\documentclass[final]{siamart250106}
\usepackage{amsmath}
\usepackage{amssymb}
\usepackage{mathtools}
\usepackage{microtype}
\usepackage{bm}
\usepackage[shortlabels]{enumitem}
\usepackage{tabularx}
\usepackage{algorithm}
\usepackage{algorithmicx}
\usepackage[noEnd,indLines,endLComment={}]{algpseudocodex}
\usepackage{overpic}
\usepackage[export]{adjustbox}
\usepackage{siunitx}
\usepackage{todonotes}
\usepackage{xspace}
\usepackage{pgfplots}
\pgfplotsset{compat=1.18}
\usepgflibrary{plotmarks}
\usepackage{cite}
\makeatletter
\newcommand{\citecomment}[2][]{\citen{#2}#1\citevar}
\newcommand{\citeone}[1]{\citecomment{#1}}
\newcommand{\citetwo}[2][]{\citecomment[,~#1]{#2}}
\newcommand{\citevar}{\@ifnextchar\bgroup{;~\citeone}{\@ifnextchar[{;~\citetwo}{]\xspace}}}
\newcommand{\citefirst}{\@ifnextchar\bgroup{\citeone}{\@ifnextchar[{\citetwo}{]\xspace}}}
\newcommand{\cites}{[\citefirst}
\makeatother

\graphicspath{{figures/}}
\makeatletter
\def\input@path{{figures/}}
\makeatother

\makeatletter
\newcommand{\multiline}[1]{%
  \begin{tabularx}{\dimexpr\linewidth-\ALG@thistlm}[t]{@{}X@{}}
    #1
  \end{tabularx}
}
\makeatother

\algnewcommand{\IIf}[2]{\State \algorithmicif\ #1\ \algorithmicthen\ #2}

\newsiamremark{remark}{Remark}
\crefname{remark}{Remark}{Remarks}
\Crefname{remark}{Remark}{Remarks}

\renewcommand{\vec}[1]{{\bm{#1}}}
\newcommand{\dom}{\Omega}
\newcommand{\bdy}{\Gamma}
\newcommand{\bdyf}{{\tilde{\bdy}}}
\newcommand{\domf}{{\tilde{\dom}}}
\newcommand{\panel}{\gamma}
\newcommand{\stripdom}{\mathcal{S}}
\newcommand{\elem}{\mathcal{E}}
\newcommand{\interface}{\bdyf}
\newcommand{\bbox}{B}

\newcommand{\uncutdom}{\dom_B}

\newcommand{\npanel}{{n_\textup{panel}}}
\newcommand{\nboxes}{{n_\textup{box}}}

\newcommand{\lap}{\Delta}
\newcommand{\part}{P}

\newcommand{\bulk}{\textup{bulk}}
\newcommand{\strip}{\textup{strip}}
\newcommand{\glue}{\textup{glue}}

\newcommand{\SL}{\mathcal{S}}
\newcommand{\DL}{\mathcal{D}}
\newcommand{\definedas}{\coloneqq}

\newcommand{\wfun}{h}

\algnewcommand{\AND}{\textbf{ and }}
\algnewcommand{\OR}{\textbf{ or }}
\algnewcommand{\NOT}{\textbf{not }}
\newcommand{\bi}{\begin{itemize}}
\newcommand{\ei}{\end{itemize}}
\newcommand{\ben}{\begin{enumerate}}
\newcommand{\een}{\end{enumerate}}
\newcommand{\be}{\begin{equation}}
\newcommand{\ee}{\end{equation}}
\newcommand{\bea}{\begin{eqnarray}} 
\newcommand{\eea}{\end{eqnarray}}
\newcommand{\ba}{\begin{align}} 
\newcommand{\ea}{\end{align}}
\newcommand{\bse}{\begin{subequations}} 
\newcommand{\ese}{\end{subequations}}
\newcommand{\bc}{\begin{center}}
\newcommand{\ec}{\end{center}}
\newcommand{\bfi}{\begin{figure}}
\newcommand{\efi}{\end{figure}}

\newcommand{\bmp}[1]{\begin{minipage}{#1}}
\newcommand{\emp}{\end{minipage}}

\newcommand{\bp}{\begin{proof}}
\newcommand{\ep}{\end{proof}}

\newcommand{\tbox}[1]{{\mbox{\tiny #1}}}

\newcommand{\C}{\mathbb{C}}

\newcommand{\R}{\mathbb{R}}

\newcommand{\bigO}{{\mathcal O}}

\newcommand{\eps}{\epsilon}

\DeclareMathOperator{\im}{Im}
\DeclareMathOperator{\re}{Re}

\DeclareMathOperator{\supp}{supp}
\newtheorem{thm}{Theorem}

\newcommand{\OB}{{\Omega_B}}

\newcommand{\tgk}{\tilde{\gamma}_k}
\newcommand{\Ek}{{E_{\rho,\tgk}}}
\newcommand{\x}{\vec{x}}
\newcommand{\y}{\vec{y}}
\newcommand{\DGL}{D_{\text{leg}}}
\newcommand{\AGL}{A_{\text{leg}}}

\title{A fully adaptive, high-order, fast Poisson solver \\ for complex two-dimensional geometries\thanks{Submitted to the editors \today.
\funding{This work was funded by the Simons Foundation.}}}

\author{Daniel Fortunato\thanks{Center for Computational Mathematics \& Center for Computational Biology, Flatiron Institute, New York, NY 10010 (\email{dfortunato@flatironinstitute.org}).} \and David B. Stein\thanks{Center for Computational Biology, Flatiron Institute, New York, NY 10010 (\email{dstein@flatironinsitute.org}).} \and Alex H. Barnett\thanks{Center for Computational Mathematics, Flatiron Institute, New York, NY 10010 (\email{abarnett@flatironinstitute.org}).}}

\headers{A fully adaptive, high-order, fast Poisson solver}{D. Fortunato, D. B. Stein, and A. H. Barnett}

\ifpdf
\hypersetup{
  pdftitle = {A fully adaptive, high-order, fast Poisson solver for complex two-dimensional geometries},
  pdfauthor = {D. Fortunato, D. B. Stein, and A. H. Barnett}
}
\fi

\begin{document}

\maketitle

\begin{abstract}
We present a new framework for the fast solution of inhomogeneous elliptic boundary value problems in domains with smooth boundaries. High-order solvers based on adaptive box codes or the fast Fourier transform can efficiently treat the volumetric inhomogeneity, but require care to be taken near the boundary to ensure that the volume data is globally smooth. We avoid function extension or cut-cell quadratures near the boundary by dividing the domain into two regions: a bulk region away from the boundary that is efficiently treated with a truncated free-space box code, and a variable-width boundary-conforming strip region that is treated with a spectral collocation method and accompanying fast direct solver. Particular solutions in each region are then combined with Laplace layer potentials to yield the global solution. The resulting solver has an optimal computational complexity of $\mathcal{O}(N)$ for an adaptive discretization with $N$ degrees of freedom. With an efficient two-dimensional (2D) implementation we demonstrate adaptive resolution of volumetric data, boundary data, and geometric features across a wide range of length scales, to typically 10-digit accuracy. The cost of all boundary corrections remains small relative to that of the bulk box code. The extension to 3D is expected to be straightforward in many cases because the strip ``thickens'' an existing boundary quadrature.
\end{abstract}



\section{Introduction}

Inhomogeneous elliptic partial differential equations (PDEs) play a central role in many areas of science and engineering, and often arise in conjunction with boundary conditions
on complicated domains.
The many fields in which this occurs include
electrostatics in the presence of a space charge,
elastostatics with a body load,
steady-state heat or chemical reaction-diffusion equations,
and (in the oscillatory case) acoustics and electromagnetics
with a distributed source.
Poisson type boundary value problems (BVPs) also arise as components of more
elaborate solvers, where they may be called a large number of times.
One example is that to solve a nonlinear elliptic BVP
a (linear) Poisson solve is needed at each quasi-Newton iteration~\cite{NonlinearBVPBook}.
A second broad example area is time-dependent solvers, in which the inhomogeneity is derived
from the solution at previous time steps.
Common applications include: (1) in computational fluid dynamics,
the pressure solve that follows each time step
in the velocity formulation for incompressible Navier--Stokes~\cite{Chorin1967} (reviewed in \cite{Guermond2006});
(2) in non-Newtonian fluids, internal stresses act as volumetric source terms
which often may combine with advection to
generate large gradients at surfaces~\cite{Stein2019}; and
(3) implicit time stepping a parabolic PDE, such as the heat or Navier--Stokes equations, using Rothe's method~\cite{Rothe1930}, which demands solvers for the
modified Helmholtz or modified Stokes (i.e., Brinkman) inhomogeneous BVPs, respectively.
In time-dependent applications, the boundary geometry may itself be evolving.

Let $\dom \subset \mathbb{R}^2$ be a simply connected domain with smooth boundary $\bdy = \partial \dom$, and let $f(\vec{x})$ for $\vec{x} \in \dom$ and $g(\vec{x})$ for $\vec{x} \in \bdy$ be smooth functions. The model inhomogeneous elliptic PDE on $\dom$ is the Poisson equation,
\begin{subequations}
\begin{alignat}{2}
\lap u(\vec{x}) &= f(\vec{x}), &&\qquad \vec{x} \in \dom, \label{eq:ipde} \\
     u(\vec{x}) &= g(\vec{x}), &&\qquad \vec{x} \in \bdy. \label{eq:ibc}
\end{alignat}
\end{subequations}
Note that the technique that we present generalizes straightforwardly to above-mentioned other PDEs, and to other boundary conditions.

Many techniques exist to discretize and solve inhomogeneous elliptic BVPs on complex geometries. Perhaps most ubiquitous are methods that operate on an unstructured volumetric mesh of the domain interior, such as finite element methods (FEMs). While they are geometrically flexible, allow for variable coefficients, and are well supported with software, the size of the resulting linear systems scales with the number of interior unknowns needed to represent $u$, which may be huge when $f$ has fine-scale variations. Furthermore, the linear systems may be ill-conditioned, and can become more so upon mesh refinement. Despite this, fast iterative methods such as multigrid~\cite{MultigridBook} and sparse direct methods~\cite{Davis2016} have made this a popular approach. However, the cost of meshing the volume, especially to high order, may be prohibitive in the setting of time-stepping with evolving geometry. This has led to FEMs based on adaptive refinement of Cartesian meshes with cut cells~\cite{Saye2017,Saye2017a}, giving decreased meshing cost. This extends ideas from immersed interface and level set methods for regular finite difference grids.

However, in the case of constant coefficients---as in \cref{eq:ipde} and many of the applications mentioned above---a potential-theory-based alternative to direct discretization allows for a much reduced number of unknowns~\cite{Mayo1992,McKenney1995}. One splits the solution $u$ as the sum of a particular solution, $v$, and a homogeneous solution, $w$, satisfying
\begin{equation}\label{eq:u_part}
\lap v(\vec{x}) = f(\vec{x}), \qquad \vec{x} \in \dom,
\end{equation}
and
\begin{equation}\label{eq:u_homo}
\begin{aligned}
\lap w(\vec{x}) &= 0,                       &&\quad \vec{x} \in \dom, \\
     w(\vec{x}) &= g(\vec{x}) - v(\vec{x}), &&\quad \vec{x} \in \bdy.
\end{aligned}
\end{equation}
The function $u = v\,+\,w$ then satisfies \cref{eq:ipde,eq:ibc}. Since the particular solution $v$ is far from unique (it need not obey any specific boundary conditions), a numerically convenient choice can be made so that evaluation of $v$ is fast and stable, as reviewed shortly.

Once $v$ is evaluated, the homogeneous BVP \cref{eq:u_homo} for $w$ must be solved; this is conveniently done using boundary element or boundary integral methods, needing only a number of unknowns sufficient to discretize the boundary and its data~\cite{Kress2014a}. This is typically orders of magnitude smaller than the system size needed with FEM.
Iterative solvers such as GMRES converge rapidly when a representation is chosen that results in a Fredholm second-kind boundary integral equation (BIE), and high-order discretizations of the boundary integral operators are available~\cite{Hao2014}. Although the resulting linear system is dense, matrix-vector products may be performed via, e.g., a fast multipole method (FMM)~\cite{Greengard1987} for the PDE fundamental solution, resulting in a solution time linear in the number of boundary unknowns. Subsequent evaluation of $w$ in the interior may then also exploit an FMM. (Note that in the non-oscillatory case the homogeneous BVP could also be solved via boundary-concentrated FEM, with a similar cost scaling~\cite{Khoromskij2003}.)

The evaluation of a particular solution obeying \cref{eq:u_part} is our main topic. There have been two main prior approaches to this:
\vspace{1em}
\begin{enumerate}
	\item One idea~\cite{Mayo1992,McKenney1995,Stein2016,Stein2017,afKlinteberg2019,Bruno2022,Stein2022b} exploits the availability of fast solvers for the Poisson problem on various simple domains. Commonly this is a uniform grid on the rectangle, with some simple boundary condition, for which there exist fast low-order finite-difference solvers via cyclic reduction~\cite{Buzbee1970} or the 2D fast Fourier transform (FFT); note that if the solution is smooth, the latter may also be used for a spectrally-accurate solution on this uniform grid. If the complex geometry $\dom$ lives within a simple domain $R$, then a particular solution on $R$---found using such a fast solver---is also a particular solution on $\dom$. (We note that fast Poisson solvers on the rectangle have recently been extended to nonuniform spectral discretizations~\cite{Fortunato2020a}, and to spheres and balls~\cite{Townsend2016,Wilber2017}, using low-rank alternating direction implicit methods.)
	\vspace{1em}
	\item Another approach~\cite{Ethridge2001,Askham2017,Anderson2023,Shen2024,Fryklund2024} finds a particular solution by convolution of $f$ with the fundamental solution (free-space Green's function), the latter being $\frac{1}{2\pi} \log \frac{1}{\|\vec{x}\|}$ for the 2D Poisson equation. This has the advantage that \textit{there is no linear solve}, merely an evaluation of a volume potential. Furthermore, the discretization of $f$ may be spatially adaptive and thus more efficient than the above FFT solvers when $f$ has fine-scale features. To evaluate this convolution in linear time, so-called ``box codes'' (or ``VFMMs'' \cite{Fryklund2022}) have been developed---FMMs specialized to an adaptive quadrature grid living on a Cartesian quadtree and reaching near-FFT speeds~\cite{Ethridge2001,Greengard1996,Askham2017}.
\end{enumerate}
\vspace{1em}
However, in both of the above approaches, the particular solution $v$ needs to be smooth enough on $\dom$ to achieve the desired order of accuracy in the overall solution $u = v\,+\,w$. There have also been two major ways to tackle this smoothness requirement:
\vspace{1em}
\begin{enumerate}
\item[(a)] The classical approach is to discretize the convolution over the domain, $v(\x) = -\int_\dom \Phi(\x,\y) f(\y) d\y$, where $\Phi$ is either the simple-domain or free-space Green's function. If done accurately, this generally gives $v$ as smooth in $\dom$ as the solution $u$ (see \cref{r:newt}). However, efficient high-order accurate approximation of the volume potential on the complex geometry is challenging. Early work in the uniform finite-difference setting extended $f$ by zero in $R \setminus \dom$ and then added careful near-boundary node corrections to recover the bulk convergence order~\cite{Mayo1992,McKenney1995}. Plain use of a box code does not solve the problem, since high accuracy would demand an excessive level of adaptive refinement towards $\bdy$.
  Recent works evaluate potentials from triangular mesh elements or irregular cut cells, to medium or high order,
  by conversion to line integrals~\cite{Anderson2023}, by density interpolation \cite{Anderson2022},
  by a two-level Ewald-type heat potential split and the nonuniform FFT \cite{Fryklund2024}, or by
  so-called {\em anti-Laplacians} (Green's 3rd identity applied elementwise) \cite{Shen2024}.
  Each of these works uses a fast algorithm to achieve $\bigO(N)$ scaling in $N$, the number of discretization
  nodes. However, of these four, in their current forms,
  only the work of Shen and Serkh \cite{Shen2024} could preserve
  linear scaling in a truly adaptive mesh.
  Furthermore, the cost of generating the needed unstructured mesh is notoriously high, especially in 3D.
	\vspace{1em}
	\item[(b)] An alternative is to use \textit{function extension} (extrapolation outside $\dom$), meaning the construction of a function $f_e$ with $f_e=f$ in $\dom$ and with some specified degree of smoothness throughout an enclosing simple domain $R$. Then $v(\x) = -\int_R \Phi(\x,\y) f_e(\y) d\y$ is a particular solution that is smooth in $\dom$, and whose evaluation (via either a fast Poisson solver or free-space box code) does not require refinement near $\bdy$. This has been used with non-adaptive spectral FFT Poisson solvers via fixed-order immersed-boundary smooth extension (IBSE) \cite{Stein2016,Stein2017}, high-order partition-of-unity extension (PUX) using radial basis functions~\cite{afKlinteberg2019,Fryklund2018,Fryklund2020}, or 1D extension along normals~\cite{Bruno2022,Epstein2022a}. In the adaptive free-space box code setting, $C^0$ extension has been done via an exterior Laplace BVP solution~\cite{Askham2017}, or at high order by an adaptive variant of PUX~\cite{Fryklund2022}. Aside from the extra computational cost (sometimes involving linear solves), the idea has two key difficulties: extrapolation is inherently ill-conditioned~\cite{Demanet2019} (becoming more so the higher the order of smoothness~\cite[Tbl.~1]{Epstein2022a}), generating large values outside $\dom$; and, close-to-touching boundaries $\bdy$ may simply not allow room for a single-valued smooth $f_e$ to exist in $R \setminus \dom$.
\end{enumerate}
\vspace{1em}
The obstacles inherent in both above approaches to the creation of a smooth $v$ motivated one of the authors recently to propose ``function intension''~\cite{Stein2022b}, namely the smooth roll-off to zero of the source term $f$ in a constant-width boundary strip region $\stripdom$ \textit{inside} $\dom$, followed by an FFT-based Poisson solve to generate a $v$ valid only in $\dom \setminus \stripdom$; a distinct particular solution is used within $\stripdom$.

\begin{figure}
  \centering
  \begin{overpic}[width=0.23\textwidth]{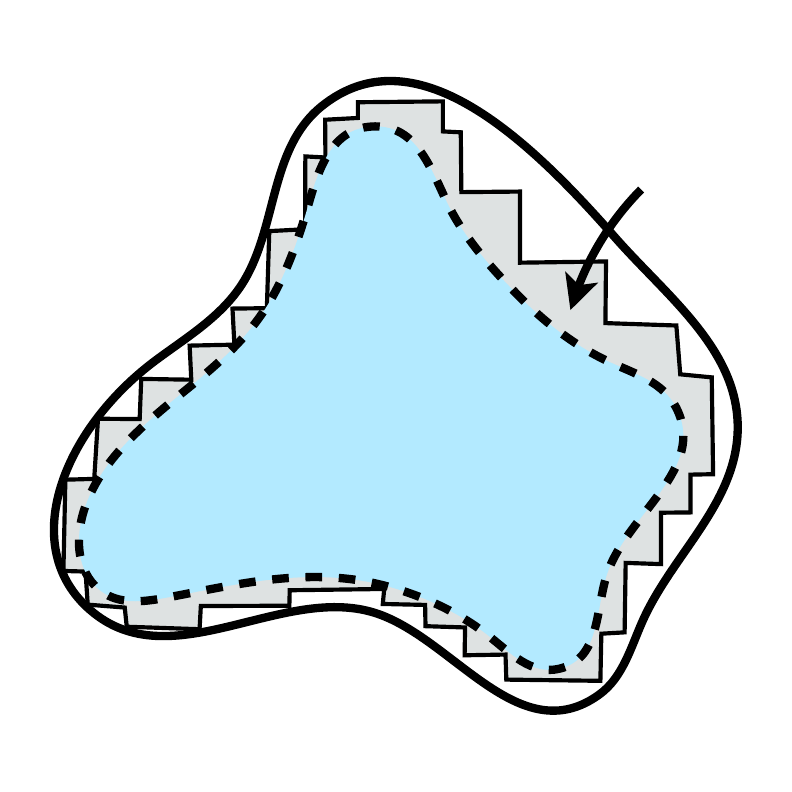}
  	\put(38,0){$v_\bulk$}
	\put(20,72){$\dom$}
    \put(81,77.3){$\dom_B$}
	\put(44,40){$\domf$}
	\put(10,9){$\bdy$}
	\put(19,28){$\bdyf$}
  \end{overpic}
  \begin{overpic}[width=0.23\textwidth]{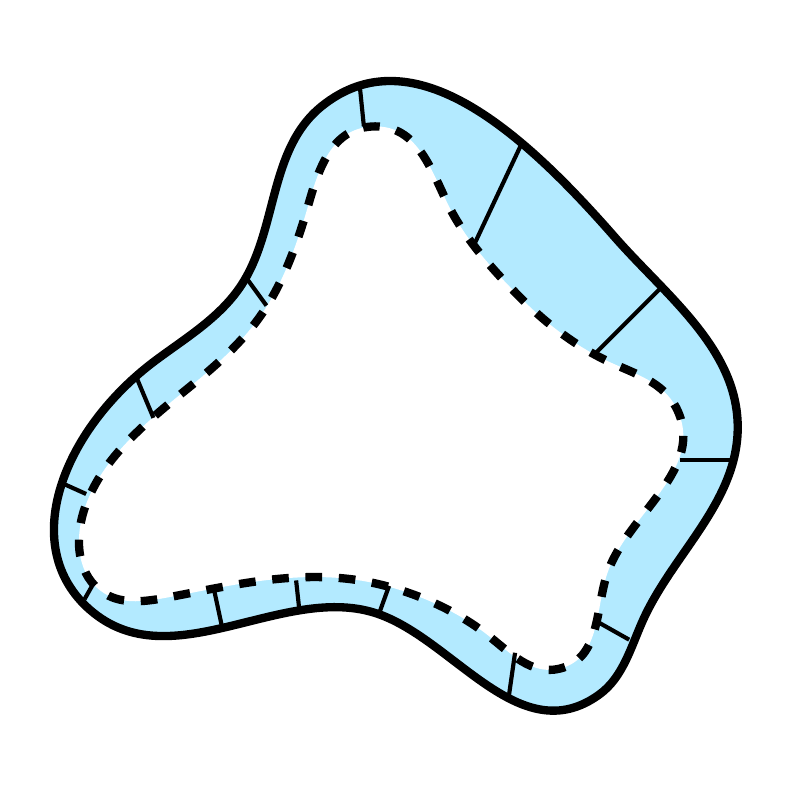}
  	\put(38,0){$v_\strip$}
	\put(67,63.5){$\stripdom$}
  \end{overpic}
  \begin{overpic}[width=0.23\textwidth]{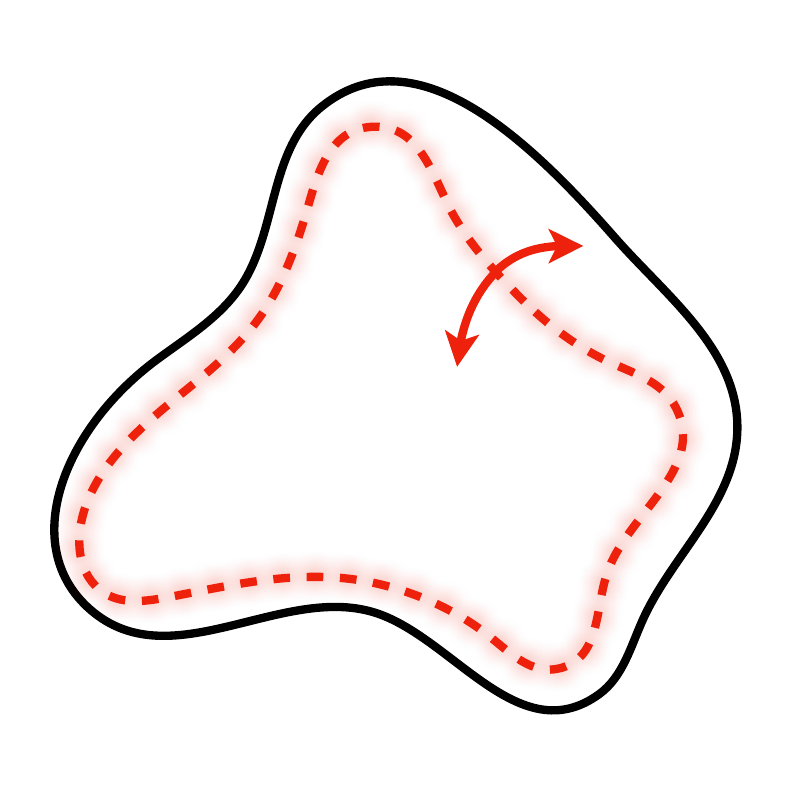}
  	\put(38,0){$v_\glue$}
	\put(25,42){\footnotesize\color{red}$\SL_\bdyf[\sigma]\!-\!\DL_\bdyf[\tau]$}
  \end{overpic}
  \begin{overpic}[width=0.23\textwidth]{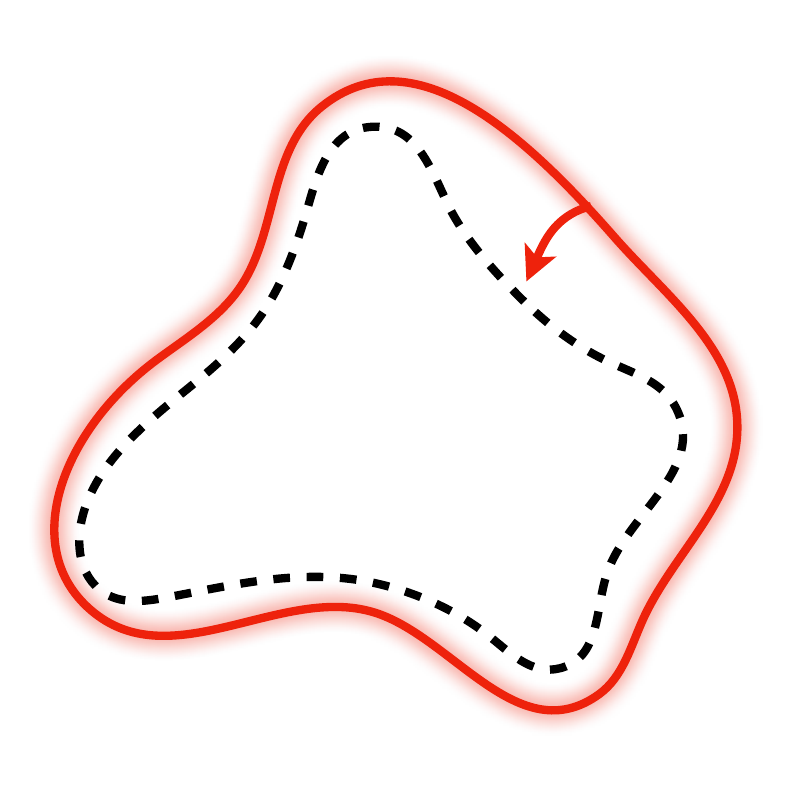}
  	\put(45,0){$w$}
	\put(77,78){\footnotesize\color{red}$\DL_\bdy[\mu]$}
  \end{overpic}
	\caption{Overview of our adaptive Poisson solver for a domain $\dom$. The left three terms comprise the particular solution $v = v_\bulk + v_\strip + v_\glue$, while $w$ is the homogeneous solution. In the left two panels, blue indicates where the particular solution is evaluated. The grey region ($\dom_B$, which includes $\domf$) created by well-separated box truncation generates the volume source for $v_\bulk$. In the right two panels, layer densities are shown in red, and evaluated throughout $\dom$. See \cref{alg:solver}.}
	\label{fig:schematic}
\end{figure}
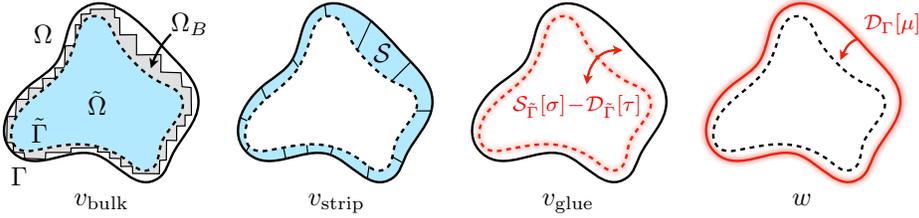

In this work, we introduce a high-order linear-scaling Poisson solver for complex geometries that is fully adaptive in handling both the inhomogeneity and the boundary, while avoiding both volume potentials on the complex geometry and function extension.
Building upon~\cite{Stein2022b}, our solver computes a particular solution $v$ using a simple decomposition of $\dom$ into two regions (see \cref{fig:schematic}):
\vspace{1em}
\begin{itemize}
\item a bulk region $\domf \subset \dom$ with boundary $\bdyf = \partial \domf$, in which we can use a free-space box-code to evaluate a particular solution $v_\bulk$ satisfying
\[
\lap v_\bulk(\vec{x}) = f(\vec{x}), \quad \vec{x} \in \domf;
\]
\vspace{-1em}
\item a thin {\em variable-width} boundary-fitted strip region $\stripdom = \dom \setminus \domf$, with boundary $\bdy \cup \bdyf$, in which we propose high-order curvilinear spectral collocation to compute a particular solution $v_\strip$ satisfying
\[
\lap v_\strip(\vec{x}) = f(\vec{x}), \quad \vec{x} \in \stripdom.
\]
\end{itemize}
As these particular solutions are piecewise defined, they generally possess jumps in both value and normal derivative across the interface $\bdyf$. We correct for these jumps by adding single- and double-layer Laplace potentials $v_\glue$ along $\bdyf$ to $v_\bulk$ and $v_\strip$, which effectively patch the solutions together\footnote{Recall that, since the PDE is second-order elliptic, values and normal derivatives are precisely the Cauchy data needed for smooth matching of solutions~\cite{EvansPDE}.} to recover a globally smooth particular solution $v$ defined everywhere inside $\dom$. Unlike in the non-adaptive work~\cite{Stein2022b} where this partitioning is defined by a constant shift in the local normal direction on $\bdy$, we seek to handle multiscale geometries and inhomogeneities which may require adaptivity. Thus, care must be taken to define the curve $\bdyf$ in a smooth geometry-aware fashion (see \cref{sec:gamma1}). This process is outlined in \cref{alg:solver}. We also improve upon~\cite{Stein2022b} by replacing the idea of ``function intension'' with ``well-separated box truncation,'' which is simpler, faster, and amenable to rigorous analysis (see \cref{t:rate}).

The paper is structured as follows. In \cref{sec:boundary_disc}, we present our standard discretization of the boundary $\bdy$ into a set of high-order panels. \Cref{sec:homo} briefly describes how a homogeneous solution $w$ may be computed using standard potential-theoretic techniques. \Cref{sec:gamma1} describes an algorithm to construct a smoothly-defined strip region in an adaptive fashion. In \cref{sec:bulk}, we describe how $v_\bulk$ can be efficiently computed as a truncated volume potential and analyze the effect that truncation has on its smoothness. \Cref{sec:strip} describes a curvilinear, composite spectral collocation method to compute $v_\strip$. In \cref{sec:glue}, we show how $v_\bulk$ and $v_\strip$ may be patched together using layer potentials to yield a globally smooth particular solution $v$. We conclude with numerical results and examples in \cref{sec:results}.

\begin{algorithm}[htb]
\caption{Adaptive solution of interior Dirichlet Poisson problem}
\begin{algorithmic}[1]
\Require{Boundary panelization of $\bdy$, inhomogeneity $f$, Dirichlet data $g$}
\Ensure{Solution $u$ to \cref{eq:ipde}--\cref{eq:ibc}}
\\[-0.5em]
\State{Construct fictitious boundary $\bdyf$ (see \cref{sec:gamma1}).}
\If{$f$ is unresolved on any elements of $\stripdom$}
    \State{Split the corresponding panels of $\bdy$ and \textbf{goto} 1}
\EndIf
\State{Construct quadtree approximation to $f$ with truncation (see \cref{sec:quadtree}).}
\State{Compute bulk solution $v_\bulk$ using a truncated box code (see \cref{sec:truncation}).}
\State{Compute strip solution $v_\strip$ using spectral collocation (see \cref{sec:strip}).}
\State{\multiline{Compute the jumps in value and normal derivative between $v_\bulk$ and $v_\strip$:
\vspace{-0.4em}
\[
\tau = v_\bulk|_\interface - v_\strip|_\interface, \qquad \sigma = \partial_\vec{n} v_\bulk|_\interface - \partial_\vec{n} v_\strip|_\interface. \vspace{-1.8em}
\]}}
\State{\multiline{Define a piecewise harmonic function to correct the jumps (see \cref{sec:glue}):
\vspace{-0.2em}
\[
v_\glue(\vec{x}) \definedas \SL_\bdyf[\sigma](\vec{x}) - \DL_\bdyf[\tau](\vec{x}). \vspace{-1.9em}
\]}}
\State{\multiline{Define the particular solution $v$:
\vspace{-0.2em}
\[
v(\vec{x}) = \begin{cases}
v_\bulk(\vec{x}) + v_\glue(\vec{x}),  &\vec{x} \in \domf, \\
v_\strip(\vec{x}) + v_\glue(\vec{x}), &\vec{x} \in \stripdom.
\end{cases} \vspace{-1.2em}
\]}}
\State{Compute the homogeneous solution $w$ by solving~\cref{eq:u_homo} (see \cref{sec:homo}).}
\State{\Return $u = v + w$}
\end{algorithmic}
\label{alg:solver}
\end{algorithm}

\section{Boundary discretization and geometry format}\label{sec:boundary_disc}

Our Poisson solver needs as input a description of the smooth domain boundary $\bdy$, a forcing function $f : \dom \to \R$, and a function $g : \bdy \to \R$ evaluating the boundary data.
We assume that the boundary is supplied as a set of disjoint panels $\{\panel_k\}_{k=1}^{\npanel}$ that resolve the boundary and such that $\bdy = \bigcup_k \panel_k$.
Specifically, the $k$th panel is described by a set of
user-supplied nodes $\{ \vec{x}_{j,k} \}_{j=1}^{p+1}$,
each node being $\vec{x}_{j,k} = (x_{j,k},y_{j,k}) \in \R^2$,
that are assumed
to be the image of the standard Gauss--Legendre nodes $\smash{\{t_j\}_{j=1}^{p+1}}$
on $[-1,1]$ under some smooth map $\vec{\Lambda}_k$.
Then $\vec{\Lambda}_k([-1,1]) = \panel_k$.
Such an input format is typical for high-order
boundary integral solvers~\cite{Wu2020}.
We typically choose the order $p$ in the range 10--20.

\begin{remark}
For the purposes of numerical tests we will need to generate resolved panel discretizations of various boundaries $\bdy$ described by an analytic or image-extracted function. This generation is common in the boundary integral equation setting and may be automatically performed in an adaptive fashion~\cite{Wu2020}.
In practice, identifying $\C$ with $\R^2$, we construct an adaptive panelization from a given parametrized curve $z(t) = x(t) + \mathrm{i} y(t)$ by numerically resolving to a specified tolerance $\eps$ a set of monitor functions: the curve $z(t)$, its parametrization ``speed'' $|z'(t)|$, and the bending energy density $|\im(z''(t)/z'(t))|^2/|z'(t)|$~\cite{Wu2020}.
For each monitor function $\zeta$ in this set and on each panel $\panel_k$, we compute the $2p+1$ Legendre coefficients $\{\hat{\zeta}_{j,k}\}_{j=1}^{2p+1}$ of $\zeta$ on $\panel_k$. Then we say that $\zeta$ is resolved on $\panel_k$ if
\[
	\sqrt{\frac{1}{p} \frac{\sum_{j=p+1}^{2p+1} \hat{\zeta}_{j,k}^2}{\sum_{j=1}^p \hat{\zeta}_{j,k}^2}} < \eps,
\]
i.e., if the tail of the Legendre coefficients has decayed to a relative tolerance of $\eps$. If $\zeta$ is not resolved on $\panel_k$, the panel is further subdivided. If all panels $\{\panel_k\}_{k=1}^{\npanel}$ are resolved to the given $\eps$, we say that the panel set resolves $\bdy$. Note that in tests of the solver, the panel nodes $\{\vec{x}_{j,k}\}_{j=1}^{p+1}$ for $1\leq k \leq \npanel$ alone are passed in to describe the geometry; the underlying parametrization is discarded.
\end{remark}

For boundary integrals with respect to arc length measure $ds$ we will need a
quadrature weight $w_{j,k}$ for each of the user-supplied set of panel nodes.
These weights are created as follows.
Let $\DGL$ be the size-$(p+1)$ square spectral differentiation matrix
that maps values to derivatives on the standard Gauss--Legendre nodes $\{t_j\}_{j=1}^{p+1}$;
the $j$th column is given by the derivative of the $j$th Lagrange basis function evaluated at the nodes.
For the $k$th panel, stacking its node coordinates as vectors and
using MATLAB-style notation
$x_{:,k} \definedas \{x_{j,k}\}_{j=1}^{p+1}$,
the speed $\|\x'(t)\|$
at the $j$th node is approximated by
\be
s_{j,k} = \sqrt{[(\DGL x_{:,k})_j]^2 + [(\DGL y_{:,k})_j]^2}.
\label{sjk}
\ee
The quadrature weight is then $w_{j,k} = W_j s_{j,k}$,
where $W_j$ are the Gauss--Legendre weights for $[-1,1]$.
Then, to high-order accuracy, for any smooth function $h:\bdy\to\R$
the change of variables from arc length to $t$ on each panel shows that
\be
\int_\bdy h(\x) ds_\x
\;\approx\;
\sum_{k=1}^\npanel \sum_{j=1}^{p+1} w_{j,k} h(\x_{j,k}).
\label{bdyquad}
\ee

The solver also needs access to the underlying arc-length parametrization induced by the set of user-supplied Gauss--Legendre panel nodes.
For each panel, say the $k$th, this is done as follows.
There exists a spectral antidifferentiation matrix $\AGL \in R^{(p+1)\times(p+1)}$ that maps derivatives to values (up to an overall constant) on the standard Gauss--Legendre nodes $\{t_j\}_{j=1}^{p+1}$; this matrix could be constructed via applying quadrature to Lagrange basis functions, but is more easily found via the pseudoinverse of $\DGL$.
Then the set $a_{j,k} = [\AGL s_{:,k}]_j$ are approximate within-panel arc-length
coordinates of the user-supplied nodes.
From this, arc-length coordinates $a(t)$ of an arbitrary $t$ in the panel may be found by polynomial interpolation from the nodes.
In particular, the panel arc length is then $a_k = a(1)-a(-1)$.
By cumulatively summing arc lengths from a fixed fiduciary panel endpoint, a global approximate arc-length parametrization of $\bdy$ is then easily built.
We will denote this by $\vec{\Lambda}(a)$, so that $\vec{\Lambda}([0,L]) \approx \bdy$ to accuracy $\eps$, where $L = \sum_{k=1}^\npanel a_k$ is the perimeter.
Note that we do not modify the user-supplied nodes (i.e., we do not reparametrize the given nodes to be Gauss--Legendre in arc length); in \cref{sec:gamma1} we will only need the $a_{j,k}$ and $a_k$ computed above.

We further require that panels are sufficiently far away from their non-neighboring panels. Specifically, the distance between a panel and any non-neighboring panel should be larger than three times the arc length of the panel \cite{Wu2020}. The given panelization is refined until this criterion is met. We use a $k$-d tree to efficiently calculate approximate panel distances~\cite{kdtree}. Finally, we require that the user-supplied panelization be level restricted, so that no two neighboring panels differ in arc length by more than a factor of two. If the given panelization does not satisfy this criterion, we refine the panelization until level restriction is satisfied.

\section{Potential theory for the homogeneous problem}\label{sec:homo}

Our main focus in the present work is the development of a fast, adaptive, and high-order accurate scheme to compute a particular solution $v$ to the inhomogeneous PDE in a complex geometry. However, we also need potential-theoretic techniques for computing the homogeneous solution $w$, which are standard and briefly described here.
The homogeneous solution $w$ to the interior Dirichlet problem \cref{eq:u_homo} is represented as the double-layer potential on $\bdy$ induced by an unknown density function $\mu$,
\be
w(\vec{x}) = \DL_\bdy[\mu](\vec{x}) := \int_\bdy \frac{\partial \Phi(\vec{x},\vec{y})}{\partial \vec{n}_\vec{y}} \mu(\vec{y}) \, ds_\y,
\label{DLP}
\ee
where $\Phi(\vec{x},\vec{y}) = \frac{1}{2\pi} \log \frac{1}{\|\vec{x}-\vec{y}\|}$ is the fundamental solution of Laplace's equation in two dimensions.
Using the jump relations of the double layer potential \cite[Thm.~6.18]{Kress2014a} leads to a second-kind integral equation for the unknown density $\mu$:
\begin{equation}\label{eq:bc_corr1}
-\frac12 \mu(\vec{x}) + \DL_\bdy[\mu](\vec{x}) = g(\vec{x}) - v(\vec{x}), \qquad \vec{x} \in \bdy.
\end{equation}
For a smooth boundary the operator \cref{DLP} has a smooth kernel,
so that a plain Nystr\"{o}m discretization~\cites[Sec.~12.2]{Kress2014a}{Hao2014} using the quadrature scheme~\cref{bdyquad} is high-order accurate. The resulting linear system for the values of $\mu$ at the set of panel nodes is well-conditioned. We solve it using GMRES, and accelerate matrix-vector products at each iteration with the 2D Laplace FMM~\cite{Greengard1987} implemented by FMMLIB2D~\cite{FMMLIB2D}. For evaluation of~\cref{DLP} close to the boundary, and for evaluation of the blocks of the system matrix between non-adjacent panels that fall sufficiently close to each other (e.g., in the case of re-entrant or close-to-touching geometries) we use the specialized panel quadrature scheme of Helsing and Ojala~\cite{Helsing2008}. In practice, we use the Helsing--Ojala scheme for any target point or target panel whose distance to the source panel is less than 1.2 times the length of the source panel.

\section{Constructing a particular solution}\label{sec:particular}

We now describe our piecewise construction of a particular solution satisfying \cref{eq:u_part}. Recall
from \cref{alg:solver}
that we form particular solutions in two regions separated by a fictitious curve $\bdyf$; see \cref{fig:schematic}.
We begin with the construction of $\bdyf$.

\subsection{Defining the fictitious curve}\label{sec:gamma1}

Given the curve $\bdy$ defined by a panelization $\{\panel_k\}_{k=1}^{\npanel}$, we aim to compute another panelized curve $\bdyf$ lying inside $\bdy$. The region between $\bdy$ and $\bdyf$ then defines the strip region $\stripdom$. To obtain a high-order accurate and scalable method, there are a number of criteria that $\bdyf$ should satisfy. The fictitious curve should be:
\vspace{0.5em}
\begin{itemize}
	\setlength\itemsep{0.1em}
	\item as smooth as the given curve $\bdy$ (or high-order accuracy may be lost);
	\item resolved using $\mathcal{O}(\npanel)$ panels (or optimal complexity may be lost);
	\item not too close to $\bdy$ (or the box code would have to adapt to overly small scales); and,
	\item not too far from $\bdy$ (or the strip region would require too many nodes in the radial direction).
\end{itemize}
\vspace{0.5em}
For simplicity, and since it induces curvilinear (as opposed to highly skew) coordinates in the strip,
we use extension in the local normal direction on $\bdy$ to define $\bdyf$.
We do not anticipate that a significant reduction in degrees of freedom is possible with a more complicated scheme.
Thus we define the fictitious curve $\bdyf$ according to a positive local width function, ${\wfun(t) : [0,L] \mapsto \R^+}$. Let $\vec{n}_{j,k}$ be the outward pointing unit normal vector at the $j$th node of panel $k$. Then, given a width function $\wfun(t)$, $\bdyf$ may be defined via perturbation in the normal direction, with each panel's nodes given by
\be
\tilde{\vec{x}}_{j,k} = \vec{x}_{j,k} - \wfun(t_{j,k}) \, \vec{n}_{j,k}
\label{eq:txjk}
\ee
for $j=1,\ldots,p+1$ and $k=1,\ldots,\npanel$, where $t_{j,k} \definedas A_k + a_{j,k}$ are the
arc-length parameters of the user-supplied nodes.
If $h$ is a smooth function (or at least as smooth as $\bdy$), then the fictitious curve $\bdyf$ will be as smooth as $\bdy$.

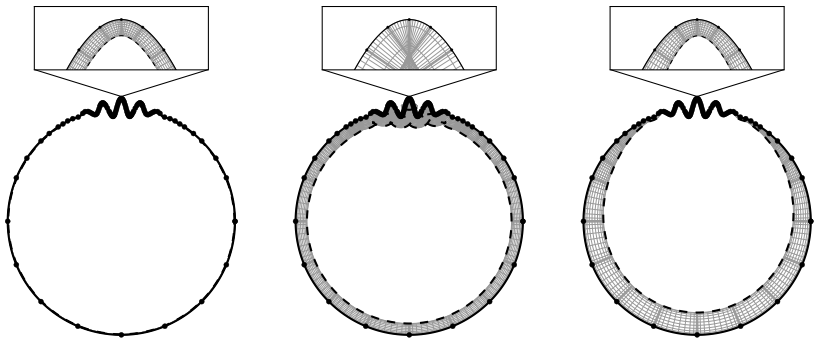
\begin{figure}[htb]
	\centering
	\begin{tikzpicture}
        \node[anchor=south west,inner sep=0] (image) at (0,0) {\includegraphics[width=0.24\textwidth]{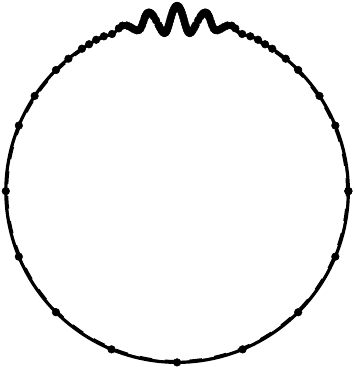}};
        \begin{scope}[x={(image.south east)},y={(image.north west)}]
            \node[anchor=south,inner sep=0,outer sep=-0.01cm] (inset) at (0.5,1.1) {\includegraphics[width=2.3cm,trim={-1cm 1cm -1cm -0.5cm},clip,frame]{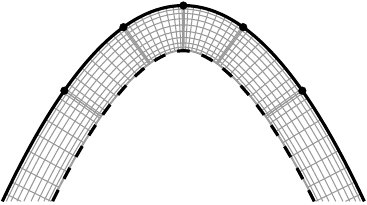}};
            \draw (inset.south west) -- (0.5,0.99);
            \draw (inset.south east) -- (0.5,0.99);
        \end{scope}
    \end{tikzpicture}
	~~~~
	\begin{tikzpicture}
        \node[anchor=south west,inner sep=0] (image) at (0,0) {\includegraphics[width=0.24\textwidth]{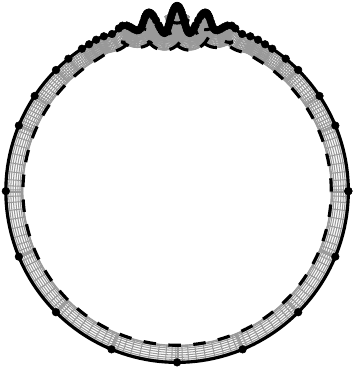}};
        \begin{scope}[x={(image.south east)},y={(image.north west)}]
            \node[anchor=south,inner sep=0,outer sep=-0.01cm] (inset) at (0.5,1.1) {\includegraphics[width=2.3cm,trim={-1cm 1cm -1cm -0.5cm},clip,frame]{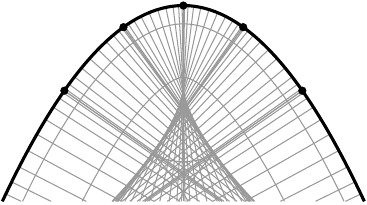}};
            \draw (inset.south west) -- (0.5,0.99);
            \draw (inset.south east) -- (0.5,0.99);
        \end{scope}
    \end{tikzpicture}
	~~~~
	\begin{tikzpicture}
        \node[anchor=south west,inner sep=0] (image) at (0,0) {\includegraphics[width=0.24\textwidth]{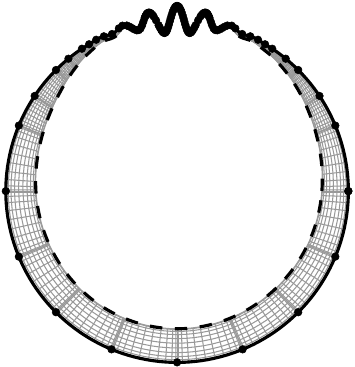}};
        \begin{scope}[x={(image.south east)},y={(image.north west)}]
            \node[anchor=south,inner sep=0,outer sep=-0.01cm] (inset) at (0.5,1.1) {\includegraphics[width=2.3cm,trim={-1cm 1cm -1cm -0.5cm},clip,frame]{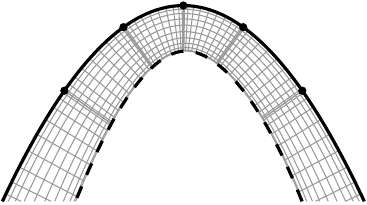}};
            \draw (inset.south west) -- (0.5,0.99);
            \draw (inset.south east) -- (0.5,0.99);
        \end{scope}
    \end{tikzpicture}
	\caption{Definition of the fictitious curve $\bdyf$ by inward normal extension. (Left) On a multiscale geometry, using a uniform strip width based on the smallest length scale results in an unnecessarily thin strip where the panel size is large. (Center) Using a width based on the largest length scale leads to self-intersection at the smallest length scales. (Right) Using a width function $h(t)$ that smoothly adapts to local panel size gives a strip that correctly resolves all geometric features.}
	\label{fig:strip_criteria}
\end{figure}

It remains to set up a width function $h$ that satisfies the above criteria.
\Cref{fig:strip_criteria} depicts the effects that different choices of $h$ can have on an adaptive geometry. Assuming that the given panelization correctly resolves all multiscale features of $\bdy$, and has been post-processed so that all non-neighboring panels are sufficiently separated,
these criteria suggest that $h$ should be proportional to the local panel size.

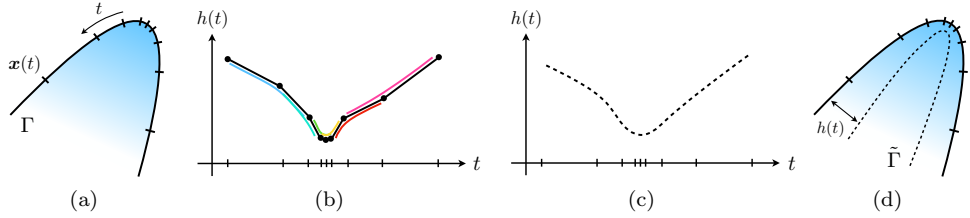
\begin{figure}[htb]
    \centering
    \begin{overpic}[width=\textwidth]{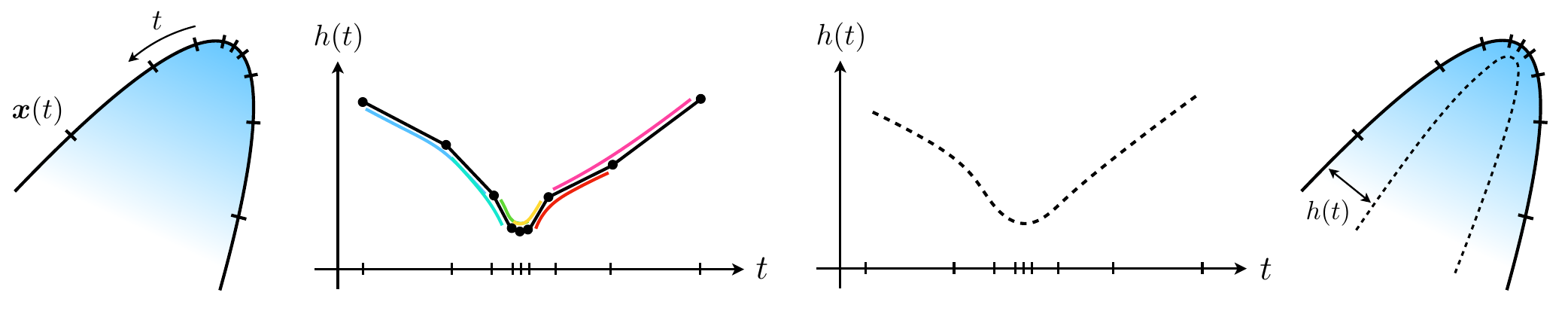}
        \put(7,-1.5){\footnotesize(a)}
        \put(32,-1.5){\footnotesize(b)}
        \put(64,-1.5){\footnotesize(c)}
        \put(89,-1.5){\footnotesize(d)}
        \put(2,6){\footnotesize$\bdy$}
        \put(90.5,2.5){\footnotesize$\bdyf$}
    \end{overpic}
    \vspace{-1em}
    \caption{Schematic for creating a smooth fictitious curve that adapts to local panel size. (a) The input panelization, given as a set of high-order Gauss--Legendre nodes sampled from the original curve $\bdy$. (b) A crude piecewise linear width function $h(t)$ is constructed from average local panel size (black circles). A rounded approximation to each kink is created, with the amount of smoothing commensurate with the local panel size (colored curves). (c) The rounded approximations are blended together into a globally smooth width function $h(t)$ using a partition of unity. (d) The original panelization is perturbed in the normal direction by $h(t)$, yielding a smooth fictitious curve $\bdyf$ that adapts to local panel size.}
    \label{fig:rounding}
\end{figure}

To define a suitable width function $h$, we begin with a crude, piecewise linear interpolation of local panel size. Let $a_k$ be the arc length of panel $k$, constructed as in \cref{sec:boundary_disc}, and define the panel endpoints $A_k \definedas \sum_{k'=1}^{k-1} a_{k'}$ so that $\vec{\Lambda}([A_k, A_{k+1}]) = \panel_k$. Define the local linear functions, $h_k^\text{lin}(t) = a_k + \frac{a_{k+1}-a_k}{a_k}(t-A_k)$. Then our starting point for a width function that adapts to local panel size is given by piecewise
linear interpolation, i.e., $h(t) = h_k^\text{lin}(t)$ for $t \in [A_k, A_{k+1}]$ and $k = 1, \ldots, \npanel$.
See~\cref{fig:rounding}(b).
(Note that, as written, the value at the endpoint $h(A_k) = a_k$ matches the panel size to its right rather than left; in practice, once $A_k$ are computed, we then replace $a_k$ with a local average of the panel sizes from the $2K$ neighboring panels centered on $A_k$, for some parameter $1 \leq K \leq 5$.)
However, this $h(t)$ is continuous but not generally $C^1$, due to kinks at endpoints.
A rounded approximation to the kink occurring between panels $\panel_{k-1}$ and $\panel_k$ at the point $A_k$ can be constructed as
\[
h_k^\text{round}(t) = a_k + \tfrac{a_k-a_{k-1}}{a_{k-1}}(t-A_k) + \left(\tfrac{a_{k+1}-a_k}{a_k} - \tfrac{a_k-a_{k-1}}{a_{k-1}}\right) r_k(t-A_k),
\]
where the middle term sets the slope for $t<A_k$ and the last term smoothly adjusts this slope to that of
$h_k^\text{lin}$ when $t>A_k$.
Here $r_k$ is a ``softplus'' (or ``smooth ReLU'') function with \emph{panel-dependent} length scale $1/\beta_k$,
\[
r_k(t) = \tfrac{1}{\beta_k} \log\left(1+e^{\beta_k t}\right),
\]
which blends from $r_k(t) \approx 0$, when $t \ll -1/\beta_k$, to $r_k(t) \approx t$, when $t \gg 1/\beta_k$.
For a length scale commensurate with local panel size, we choose $\beta_k = 2/(a_{k-1}+a_k)$.

We then combine the ``inner'' rounded approximations to each kink into smooth functions defined on each panel,
subtracting off the ``outer'' expansion $h^\text{lin}$ in the manner of matched asymptotics \cite[\S3.3.3]{logan}.
Specifically, on each panel $\panel_k$ we add together the rounded approximations from a small set of $2K$ neighboring panels, $\{\panel_{k-K}, \ldots, \panel_{k-1}, \panel_{k+1}, \ldots, \panel_{k+K}\}$, so that
\[
h_k^\text{neigh}(t) = h_{k+K+1}^\text{round}(t) + \!\!\sum_{j=k-K}^{k+K} \left( h_j^\text{round}(t) - h_j^\text{lin}(t) \right).
\]
In practice, we choose $1 \leq K \leq 5$.
The function $h_k^\text{neigh}(t)$ is locally a smooth function on panel $\panel_k$ and its adjacent $2K$ neighbors.

Finally, in order to arrive at a width function $h(t)$ that is globally smooth across all panels, we blend together the functions $\{h_j^\text{neigh}(t)\}_{j=k-K}^{k+K}$ on panel $\panel_k$ using a partition of unity, yielding
\be
h(t) = \!\! \sum_{j=k-K}^{k+K} \hat{w}_j(t) h_j^\text{neigh}(t), \qquad t \in [A_k, A_{k+1}],
\label{eq:ht}
\ee
with normalized blending functions $\hat{w}_k(t) = w_k(t) / \sum_{j=k-K}^{k+K} w_j(t)$, $w_k(t) = \rho_c(\tfrac{t-A_k}{\delta})$, and bump function $\rho_c(t)$ taken to be the prolate spheroidal wavefunction of order zero and bandwidth $c$~\cite{slepianI},
implemented by {\tt pswf.m} in Chebfun~\cite{Chebfun}. We typically choose $\delta = \tfrac{a_{k-1} + a_k}{8}$ and $c=30$.
It is \cref{eq:ht} that is used as $h(t)$ in \cref{eq:txjk}.

\subsection{The bulk problem}\label{sec:bulk}

With a suitably defined fictitious curve $\bdyf$ separating the bulk region from the strip region, we now turn to computing a particular solution $v_\bulk$ in the bulk region $\domf$, meaning that it satisfies
\begin{equation}\label{eq:bulk_problem}
\lap v_\bulk(\vec{x}) = f(\vec{x}), \quad\vec{x}\in\domf.
\end{equation}
For this, our procedure is similar to that of a box code (fast evaluation of a volume potential on a
grid adapted to resolve $f$ \cite{Greengard1996,Ethridge2000,Ethridge2001,Askham2017}),
but with a modified criterion for refinement, as well as for exclusion of well-separated
boxes.
Its result will be a $v_\bulk$ that is a high-order approximation to
\be
v_\bulk(\vec{x}) = -\int_{\Omega_B} \Phi(\vec{x}, \vec{y}) \, f(\vec{y}) \, d\vec{y},
\label{eq:vbulk}
\ee
where $\Omega_B$ is some union of boxes such that $\domf \subset \Omega_B \subset \dom$,
as sketched by the grey region with irregular boundary in the left panel of \cref{fig:schematic}.
The requirement $\Omega_B \subset \dom$ simply comes from the fact that
$f$ is not defined outside $\dom$, and,
for reasons discussed in the introduction, extending $f$ smoothly outside of
$\dom$ has difficulties that we wish to avoid.

It is clear that \cref{eq:vbulk} would satisfy \cref{eq:bulk_problem} for
{\em any} choice of source domain $\Omega_B\subset \dom$ that entirely covers $\domf$.
However, what is needed is a choice of $\Omega_B$ that results in
$v_\bulk$ being smooth on $\bdyf$ and in $\domf$.
This would indeed hold if $\Omega_B=\domf$, but would require excessive
box refinement around $\bdyf$ down to the scale of the desired tolerance $\epsilon$
(thus adding $\bigO(\epsilon^{-1} \log \epsilon^{-1})$ boxes to resolve the boundary),
hence would be very inefficient.
The same would be true at the other extreme case $\Omega_B=\dom$, with the excessive
refinement now occurring at $\bdy$.
Our choice of $\Omega_B$ lies between these two extremes, and is based on selectively
truncating boxes in the strip region.
Despite its seemingly awkward shape, we will show that it retains high-order accuracy
while inducing minimal extra box refinement.

\begin{remark}
  The known values of $f$ outside $\domf$ may also be used to obtain a globally smooth function $\tilde{f}$ through multiplication by a high-order smooth blending function which rolls off to zero in the strip region $\stripdom$, in the style of ``function intension''~\cite{Stein2022b}.
  Since $\tilde{f}$ is now less smooth than $f$, this would typically require smaller adaptive boxes in $\stripdom$, and
  hence would be less efficient than the box truncation that we propose.
\end{remark}

\subsubsection{Constructing a quadtree approximation to \texorpdfstring{$f$}{f}}\label{sec:quadtree}

The problem defines the domain $\dom$ and the source function $f$ available on $\dom$,
and we have at this point fixed the fictitious split of $\dom$ into bulk $\domf$ and strip $\stripdom$
along the curve $\bdyf$.
The user specifies the order $p$ and tolerance $\epsilon>0$.
We will take $p=16$ in all examples.
We now aim to construct a quadtree consisting of a set of boxes $\{\bbox_k\}_{k=0}^{\nboxes-1}$,
where each box is either a parent or a leaf, and
$f$ is approximated by an order-$p$ polynomial up to a tolerance $\epsilon$ on leaf boxes.

\begin{figure}[htb]
	\centering
	\vspace{1em}
	\setlength{\tabcolsep}{0.22cm}
	\newlength{\mywidth}
	\newlength{\cbheight}
	\setlength{\mywidth}{0.175\textwidth}
	\setlength{\cbheight}{0.18\textwidth}
	\begin{tabular}{c@{\hskip 0.4cm}ccccc}
		\adjustbox{valign=m}{\rotatebox{90}{\scriptsize Inhomogeneity $f$}} &
		\includegraphics[width=\mywidth,valign=m]{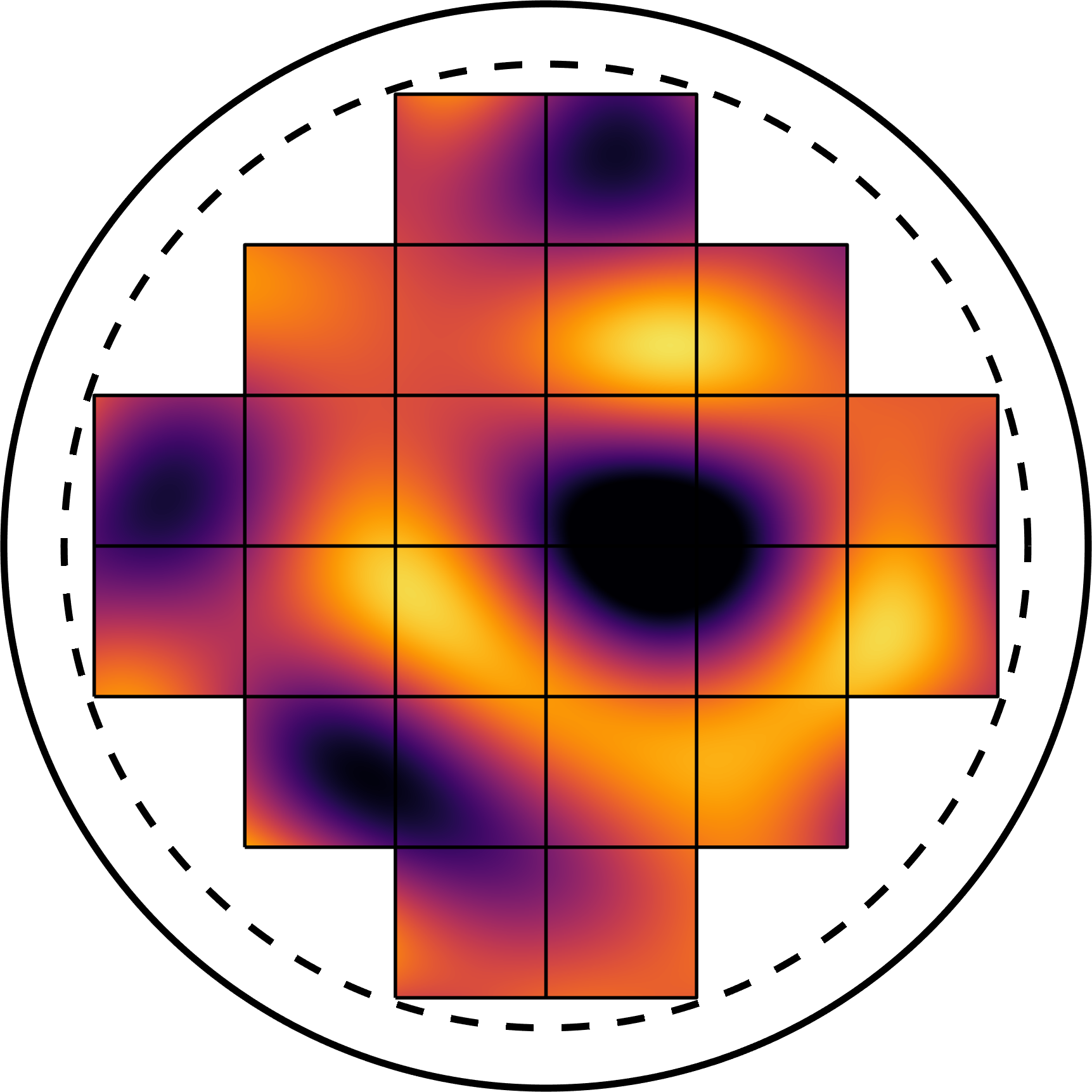} &
		\includegraphics[width=\mywidth,valign=m]{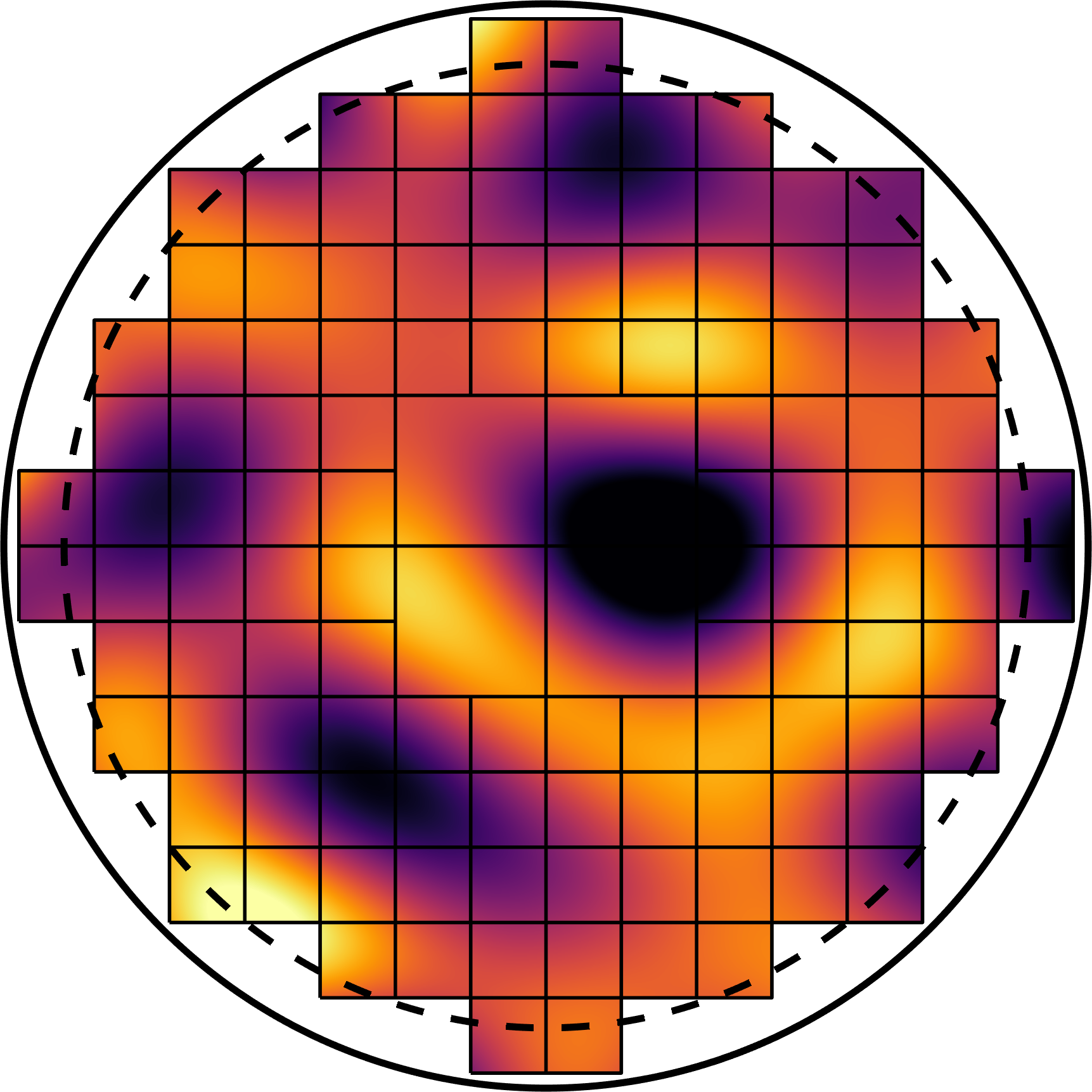} &
		\includegraphics[width=\mywidth,valign=m]{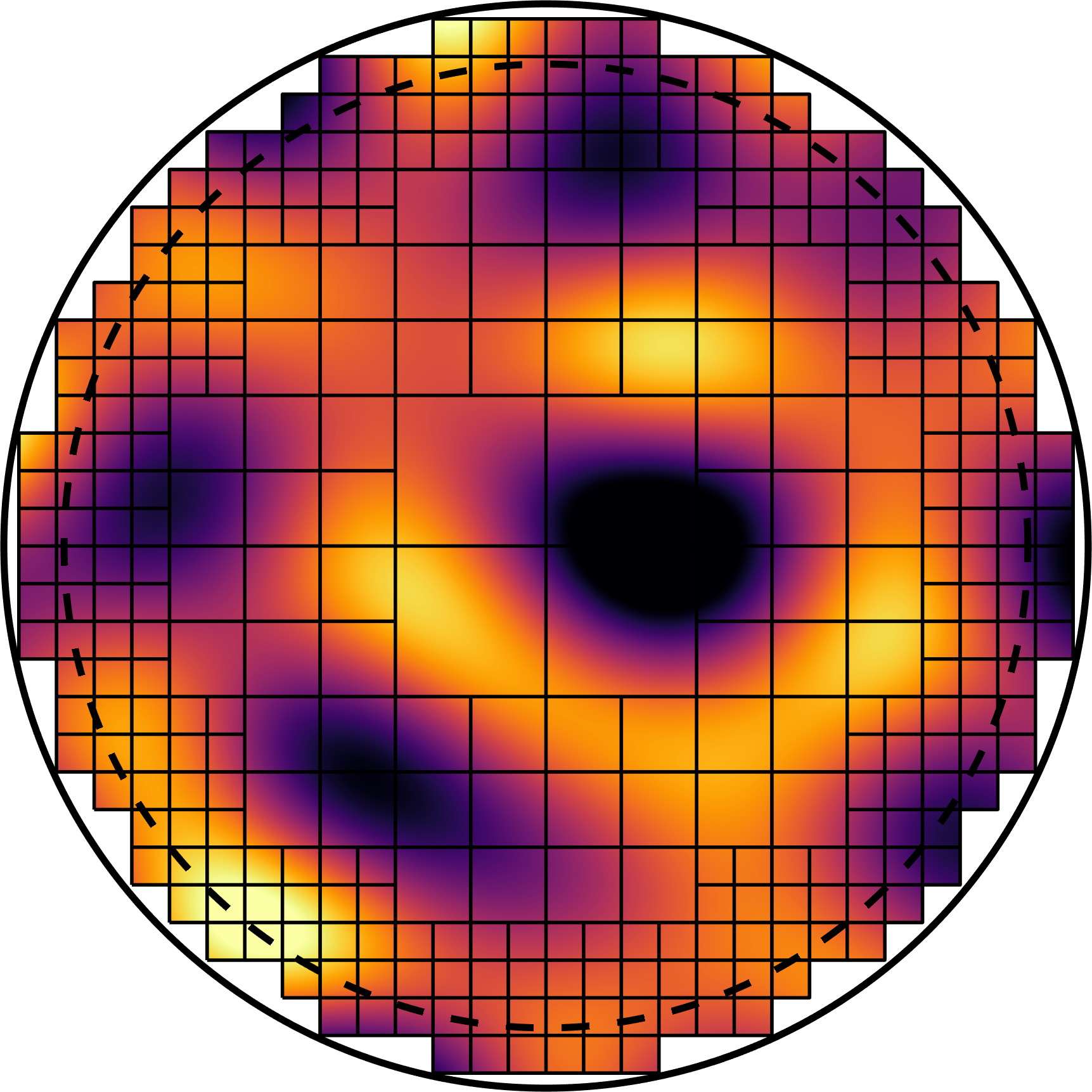} &
		\includegraphics[width=\mywidth,valign=m]{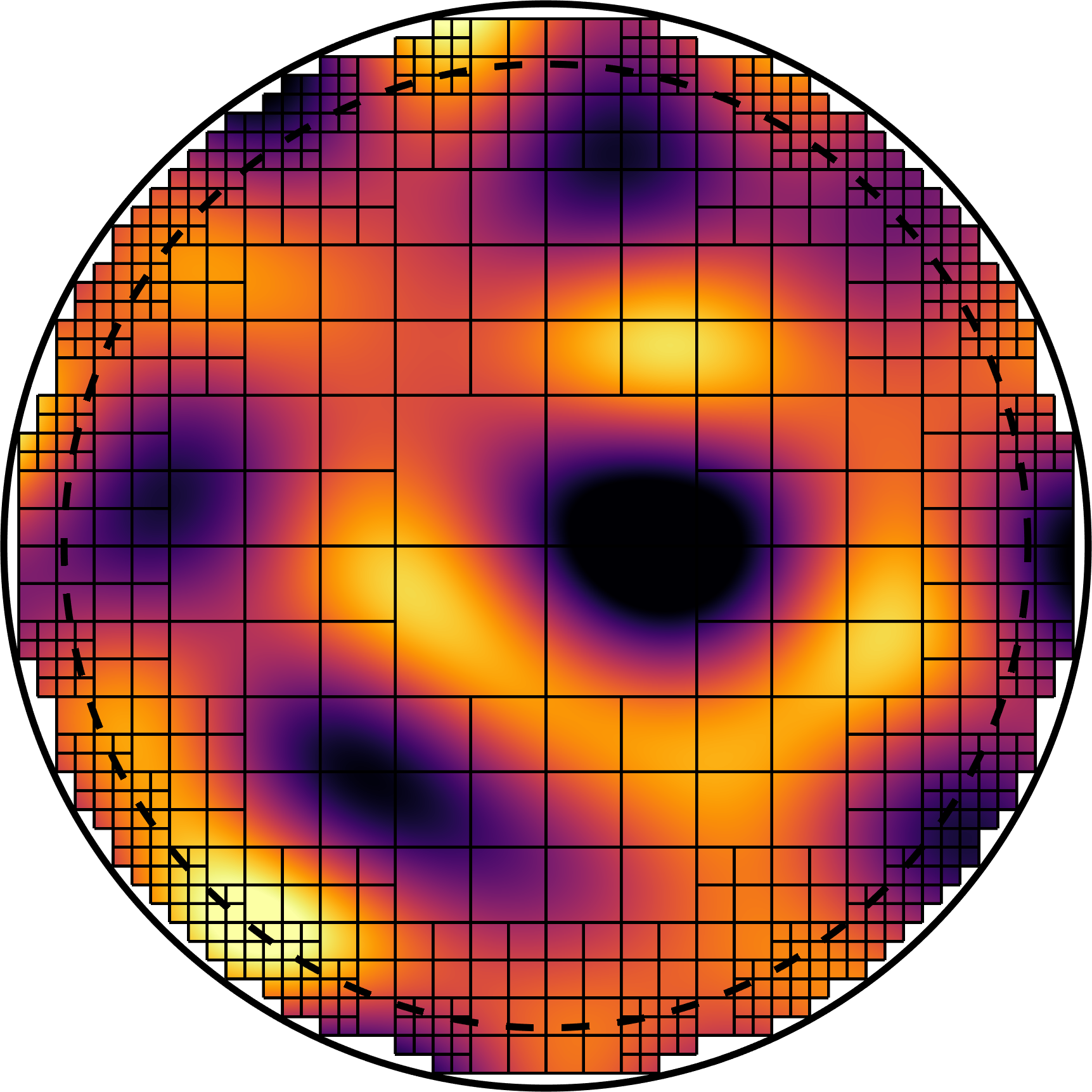} &
		\includegraphics[height=\cbheight,valign=m]{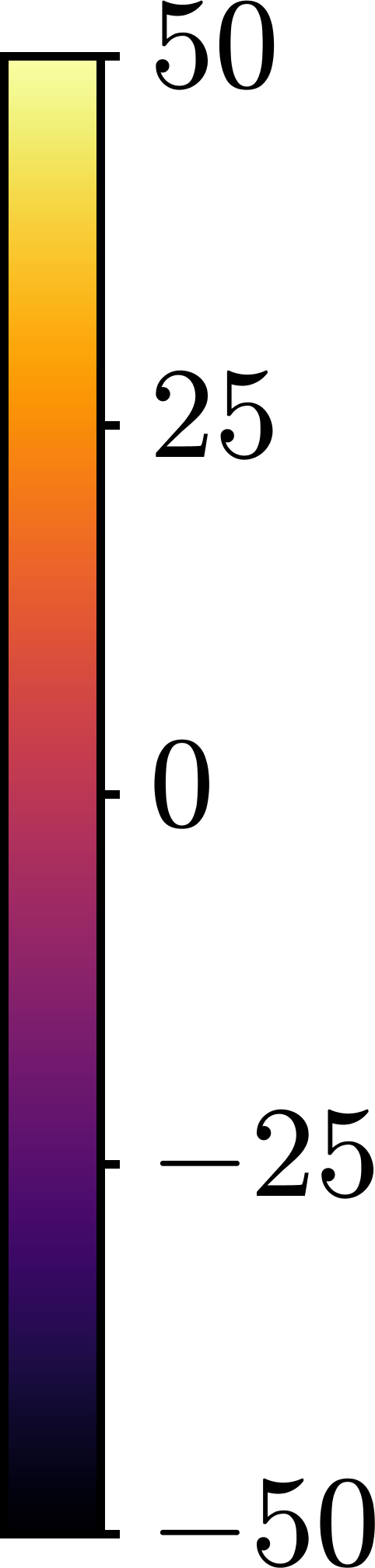} \\[1.2cm]
		\adjustbox{valign=m}{\rotatebox{90}{\scriptsize $\dfrac{\|\lap v_\bulk - f\|_\infty}{\|f\|_\infty}$}} &
		\includegraphics[width=\mywidth,valign=m]{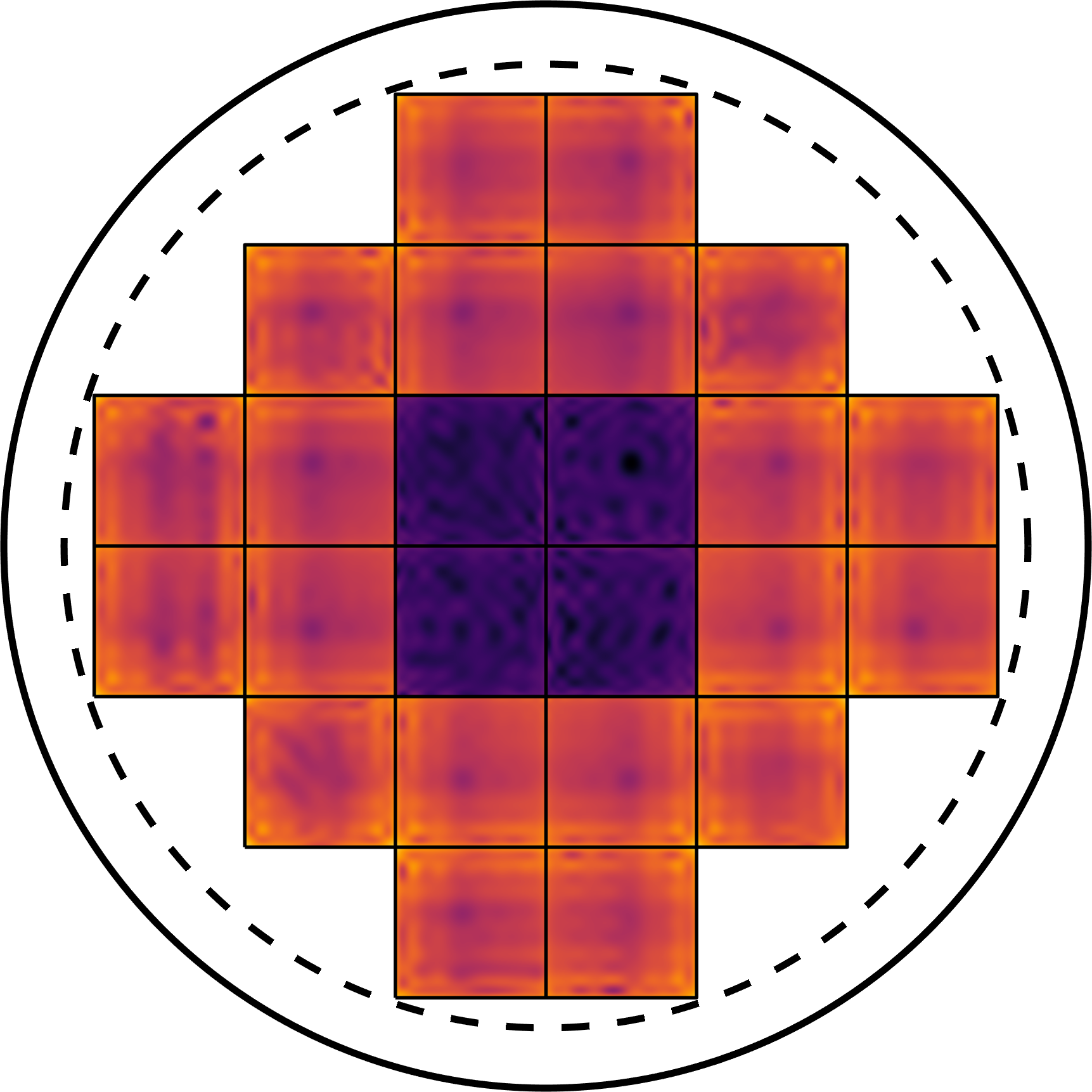} &
		\includegraphics[width=\mywidth,valign=m]{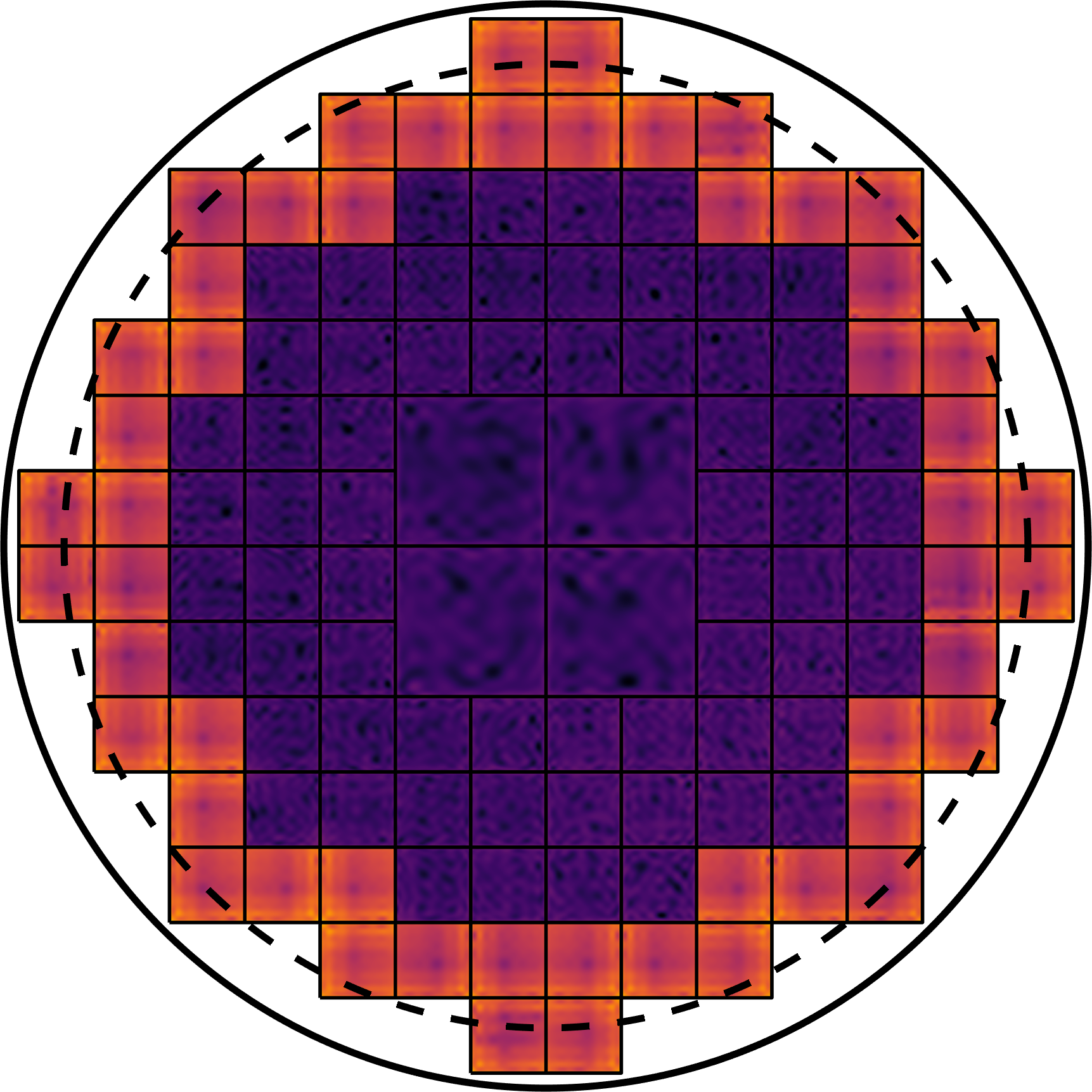} &
		\includegraphics[width=\mywidth,valign=m]{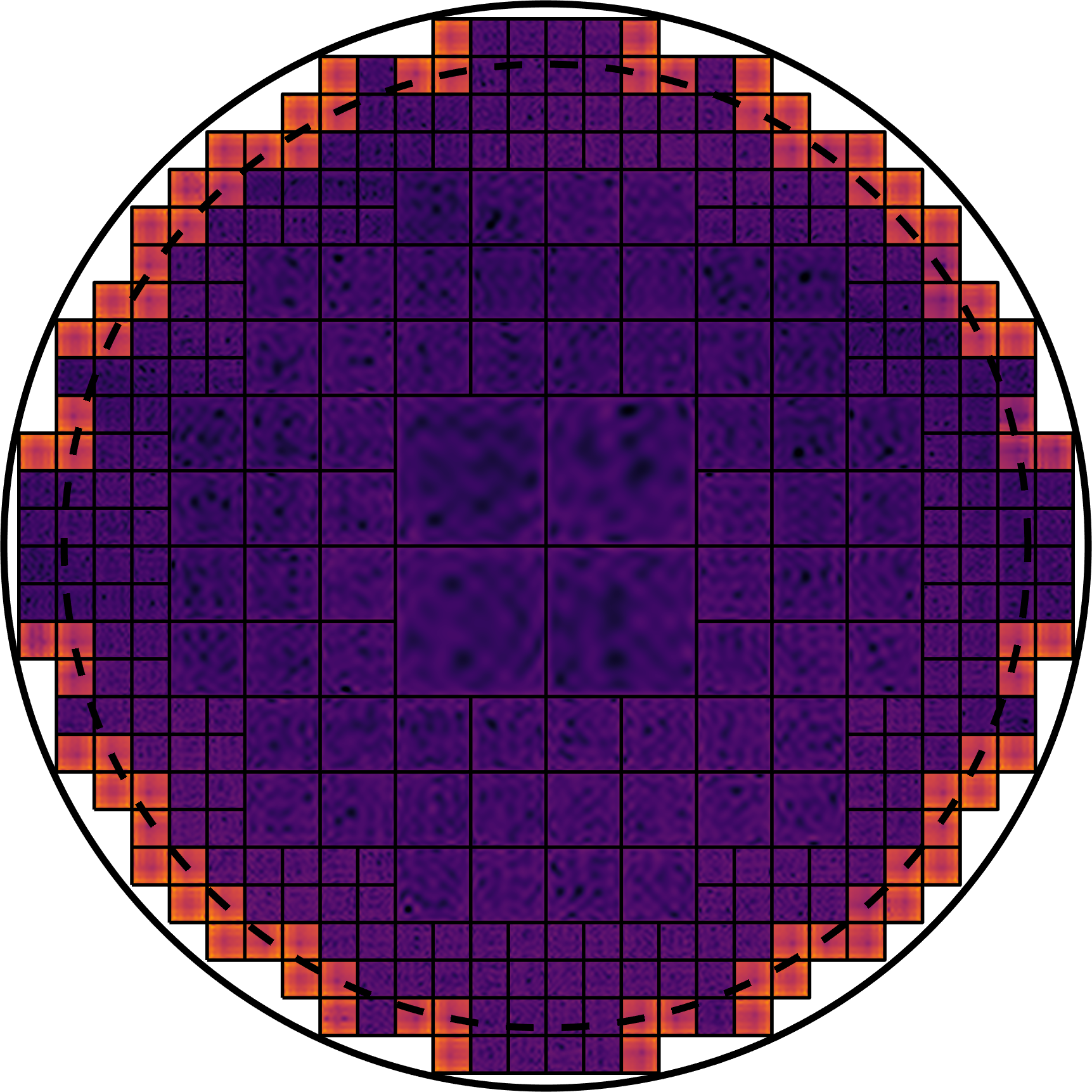} &
		\includegraphics[width=\mywidth,valign=m]{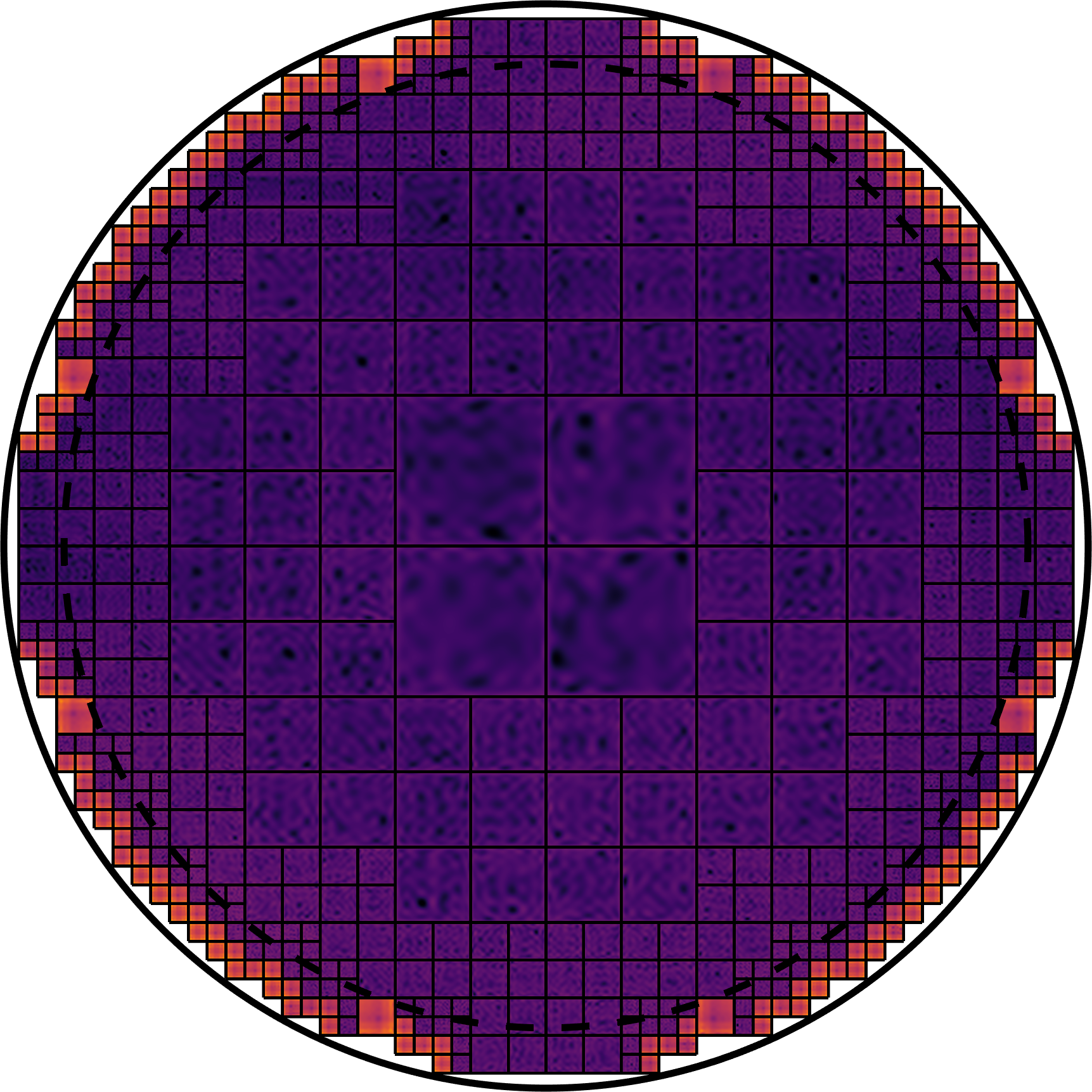} &
		\includegraphics[height=\cbheight,valign=m]{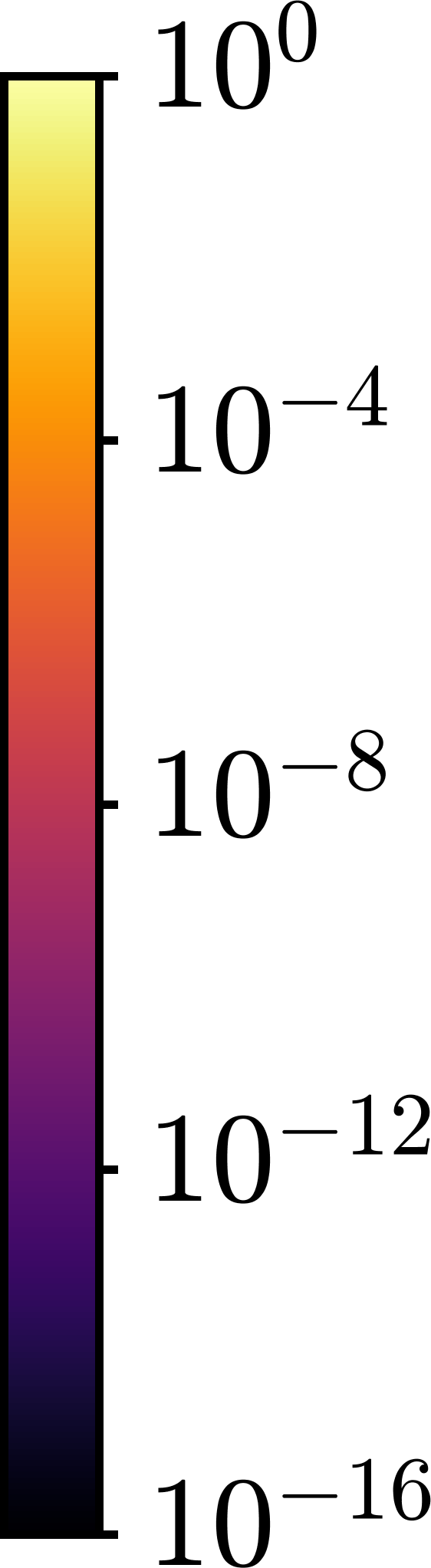} \\[1.2cm]
		\adjustbox{valign=m}{\rotatebox{90}{\scriptsize $\dfrac{\|u-u_\text{ref}\|_\infty}{\|u_\text{ref}\|_\infty}$}} &
		\includegraphics[width=\mywidth,valign=m]{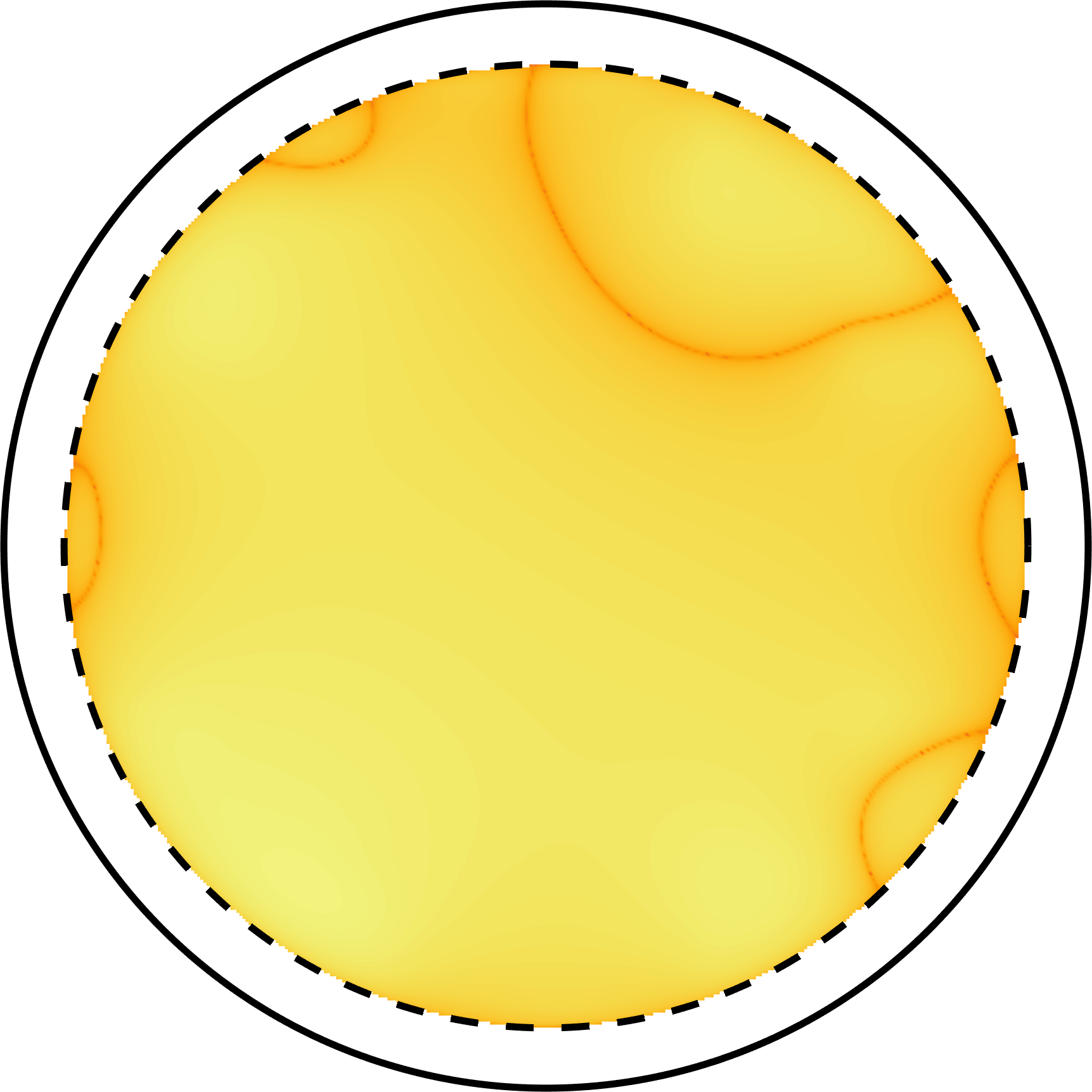} &
		\includegraphics[width=\mywidth,valign=m]{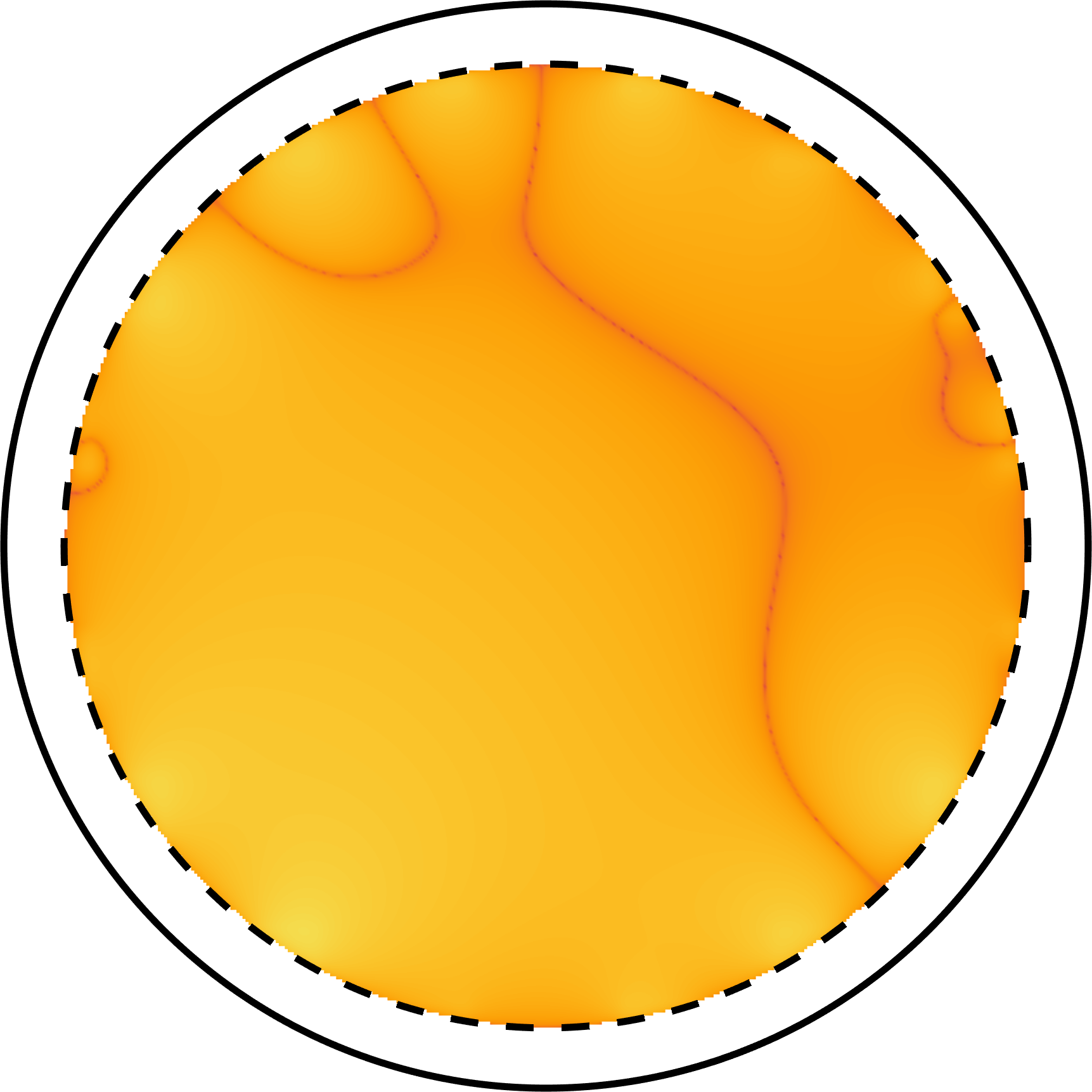} &
		\includegraphics[width=\mywidth,valign=m]{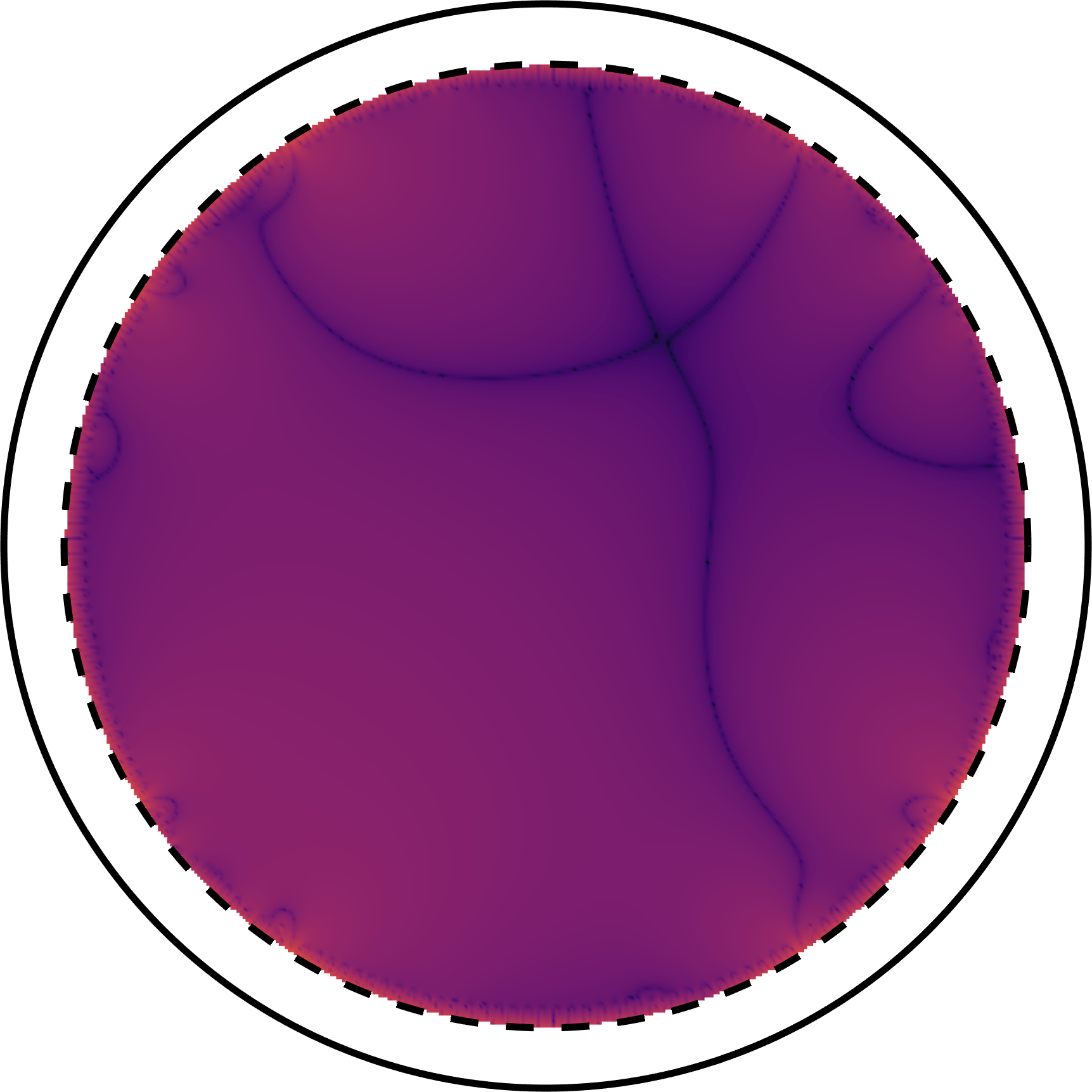} &
		\includegraphics[width=\mywidth,valign=m]{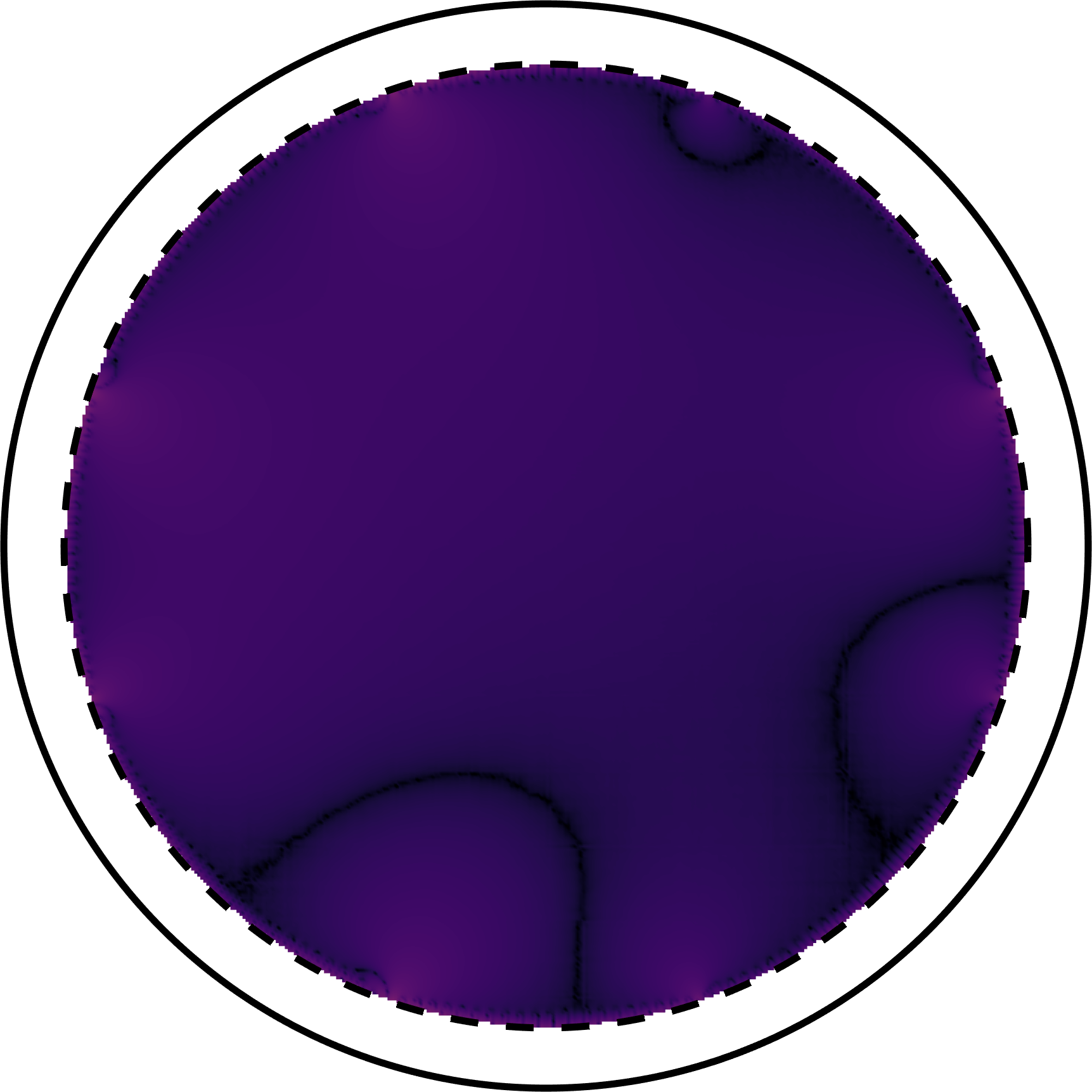} &
		\includegraphics[height=\cbheight,valign=m]{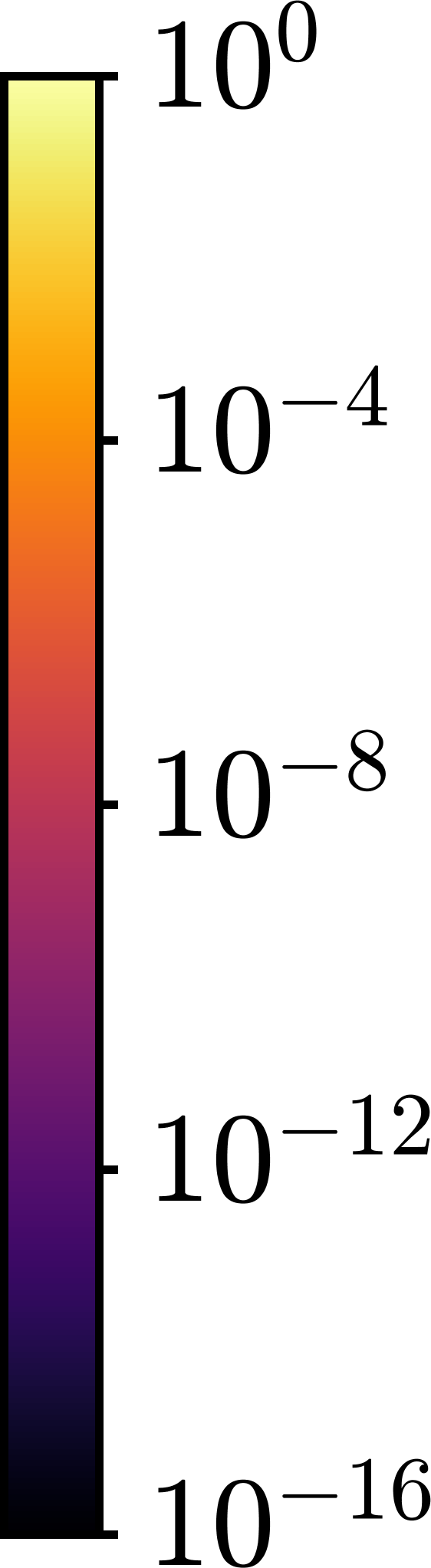}
	\end{tabular}
	\caption{Visualization of the effect of truncation when computing a volume potential. (Top row) An inhomogeneity $f$ known everywhere inside the domain (solid circle) is converted to a quadtree, with each leaf approximating $f$ to high order. Any boxes overlapping the outside of the domain are discarded, as we assume $f$ is known only inside the domain. (Middle row) The volume potential $v_\bulk$ induced by this truncated data is computed with a box code and the relative error of the residual is used as an indicator for refinement. Boxes with large residuals are subdivided and solutions are successively computed as we move column-wise to the right. (Bottom row) To visualize the effect of truncation on the final solution, a test problem is solved inside the fictitious region (dotted circle) using $v_\bulk$ as the particular solution and the computed solution $u$ is compared to a known reference solution $u_\textup{ref}$. The final quadtree suggests that boxes touching the strip region should be refined until their diameter is less than half of the strip width.}
	\label{fig:truncation}
\end{figure}

The quadtree construction process starts from a single box $B_0 \supset \dom$ containing the entire domain, and proceeds to recursively split a box $\bbox_k$ into four child boxes according to two criteria:
\ben
\item \textbf{Function resolution}. If $\bbox_k \subset \dom$, then $f$ may be fully evaluated on $\bbox_k$ and its Chebyshev coefficients $\hat{f}_{ij}^k$ computed according to the order-$p$ approximation,
\[
f(x,y) \approx \sum_{i=1}^{p+1} \sum_{j=1}^{p+1} \hat{f}_{ij}^k \, T_j^k(x) \, T_i^k(y), \qquad (x,y) \in \bbox_k,
\]
where $T_j^k$ is the $j$th Chebyshev polynomial of the first kind scaled to the domain of box $\bbox_k$.
If the coefficients of $f$ on $\bbox_k$ do not decay below the given tolerance $\epsilon$, then $\bbox_k$ is marked for refinement.
If $\bbox_k$ is entirely outside of $\dom$,
it is discarded and does not contribute.
However, if only part of $\bbox_k$ falls outside $\dom$
(i.e., the boundary $\bdy$ cuts the box), 
then extending $f$ by zero outside $\dom$ will cause its polynomial approximation to fail.
This necessitates the use of the second criterion.
\item \textbf{Separation from truncation}. If $\bbox_k$ intersects $\bdy$, then the corresponding volume potential computed on $\bbox_k$ will not be an accurate particular solution even in the part of the box
  lying in $\dom$.
  One simple option is to set the coefficients $\hat{f}_{ij}^k$ to zero in $\bbox_k$.
But how does this truncation affect the neighbors of $\bbox_k$? To examine this, we perform an experiment on a simple domain in \cref{fig:truncation}. Starting from an inhomogeneity $f$ known everywhere inside a circle, we construct a quadtree which approximates $f$ to high order and zero
any boxes which overlap the outside of the domain (\cref{fig:truncation}, top row). We then compute the volume potential $v_\bulk$ induced by this truncated data and evaluate its residual (\cref{fig:truncation}, middle row).
One might expect that the truncation generally induces {\em fake corner singularities} in the
resulting $v_\bulk$ which cause large residuals on any adjacent box---and indeed this is what is seen.
Boxes with large residual, as well as those that were cut by $\bdy$, are marked for refinement.
Moving left-to-right across the figure,
each column shows a recomputation after boxes marked for refinement are subdivided.
To examine the global effect, we solve a test problem inside $\domf$ using $v_\bulk$ as the particular solution, and measure the error against a known reference solution (\cref{fig:truncation}, bottom row).
This shows that the global error in $\domf$ is controlled by this truncated-neighbor effect.
The final quadtree suggests that any box that is well separated from truncation (i.e., having no truncated neighbors) will be accurate. Hence, for all boxes in $\domf$ to be accurate, all boxes touching the strip region $\stripdom$ should be refined until their diameter is less than half of the strip width. This is our second criterion.
\een

A quadtree construction process based on these criteria is given in \cref{alg:quadtree}.
As a consequence of the size-based splitting criterion, all truncation is pushed to the
{\em outer half} of the strip $\stripdom$.
As a final post-processing step, we perform a 2:1 balance of the resulting quadtree~\cite{deBerg2008,Malhotra2016}, so that every box is no more than one refinement level apart from its neighbors. The truncated volume potential is computed using the box code available in the \texttt{boxcode2d} library~\cite{boxcode2d}, which is called by the \texttt{treefun} package~\cite{treefun}.

\begin{remark}\label[remark]{r:keepcut}
  In practice, we choose to utilize the Chebyshev coefficients of $f$ even on cut boxes (where $f$ is taken to be extended by zero outside $\dom$), rather than setting the function to zero uniformly over the entire cut box. Though the resulting polynomial approximation of $f$ on these cut boxes is unresolved, we have found that including such cut values of $f$ can increase the overall accuracy by half a digit or more, by effectively moving the truncation location slightly further from $\bdyf$.
\end{remark}

\begin{algorithm}[htb]
\caption{Adaptive construction of truncated quadtree approximation to $f$}
\begin{algorithmic}[1]
\Require{Inhomogeneity $f$, boundaries $\bdy$ and $\bdyf$, bounding box $\bbox_0 \supset \dom$, tolerance $\epsilon > 0$}
\Ensure{Quadtree approximation to $f$ as a set of boxes $\{\bbox_k\}_{k=0}^{\nboxes-1}$}
\\[-0.5em]
\State Initialize stack of boxes $B \gets \{ \bbox_0 \}$ and $\nboxes \gets 1$
\While{$B \neq \emptyset$}
	\State Pop $\bbox_k$ from $\bbox$
    \State Determine where $\bbox_k$ lies in relation to $\bdy$ and $\bdyf$
    \State $\texttt{resolved} \gets \texttt{true}$
	\\[-0.5em]
	\LComment{Function resolution criterion}
			\State Compute Chebyshev coefficients $\hat{f}_{ij}^k$ of $f$ on $\bbox_k$ (with $f$ taken to be extended by zero outside $\dom$)
	\If{$\bbox_k \subset \dom$ \AND $\hat{f}_{ij}^k$ is \NOT resolved to $\epsilon$}
		\State $\texttt{resolved} \gets \texttt{false}$
	\EndIf
	\vspace{1em}
	\LComment{Truncation separation criterion}
	\If{$\bbox_k \cap \domf \neq \emptyset\!$ \AND $\bbox_k \cap \R^2\backslash\dom \neq \emptyset$}
		\State $\texttt{resolved} \gets \texttt{false}$
	\ElsIf{$\bbox_k \cap \stripdom \neq \emptyset$}
		\State Find nodes $\vec{x} \in \bdy$ and $\tilde{\vec{x}} \in \bdyf$ nearest to the center of $\bbox_k$
		\State Compute local strip width $s \gets \|\vec{x} - \tilde{\vec{x}}\|_2$
		\IIf{$\text{diagonal length of $\bbox_k$} > s/2$}{$\texttt{resolved} \gets \texttt{false}$}
	\EndIf
	\vspace{1em}
	\If{\texttt{resolved}}
		\State Store $\bbox_k$ as a leaf box with Chebyshev coefficients $\hat{f}_{ij}^k$
	\Else
		\State Store $\bbox_k$ as a parent box
		\State Refine $\bbox_k$ and push children onto $B$
		\State $\nboxes \gets \nboxes + 4$
	\EndIf
\EndWhile
\vspace{0.7em}
\State{\Return} $\{\bbox_k\}_{k=0}^{\nboxes-1}$
\end{algorithmic}
\label{alg:quadtree}
\end{algorithm}

\subsubsection{Smoothness of truncated volume potentials}\label{sec:truncation}

We now analyze the effect of truncating the volume potential in terms of the smoothness of the computed particular solution on each panel of $\bdyf$.
Recall from the introduction that a classical particular solution using the full solution domain
$\dom$ is
\be
v(\vec{x}) = -\int_{\dom} \Phi(\vec{x}, \vec{y}) \, f(\vec{y}) \, d\vec{y}.
\label{eq:vps}
\ee
This is smooth in the interior if $f$ is smooth, as follows from standard elliptic regularity \cite[\S6.3.2]{EvansPDE} since $\Delta v = f$ in $\dom$.

The previous subsection described our method to truncate the volume potential
to give \cref{eq:vbulk}, produced by limiting the support of the source to $\OB$, a union of boxes.
Then the difference from the classical potential \cref{eq:vps},
\[
\tilde v := v_\bulk - v,
\]
is simply the potential due to the difference $\tilde f$ in the source terms, namely
the free space convolution
\be
\tilde v = -\Phi \ast \tilde f,
\qquad \mbox{ where } \quad \tilde f := \chi_{\dom \backslash \OB} f,
\label{tildev}
\ee
where $\chi_S$ denotes the characteristic function of a set $S$,
and we note that $\supp \tilde f \subset \R^2 \backslash \OB$ and $\|\tilde f\|_{L^1(\R^2)} < \infty$.
Additional cut cells included as per \cref{r:keepcut} change $\tilde f$ but
do not change the fact that the {\em support lies outside of} $\OB$.
The proof that $\tilde v$ is smooth on $\bdyf$ will rely entirely on this fact that $\bdyf$ is
well-separated from $\supp \tilde f$; the roughness of $\tilde f$ is not relevant.

Our analysis, being based on spatial well-separation, is of a different flavor from that of
prior work such as \cite[Sec.~4]{Askham2017}, which showed that for $f$ discontinuous in a box
one expects nearly two orders (with respect to the box size)
better convergence in $v$ than in the Chebyshev representation of $f$.

\begin{figure}
\centering
\begin{minipage}{0.62\textwidth}
  \includegraphics[width=\textwidth,trim={1.5cm 0 0 0},clip]{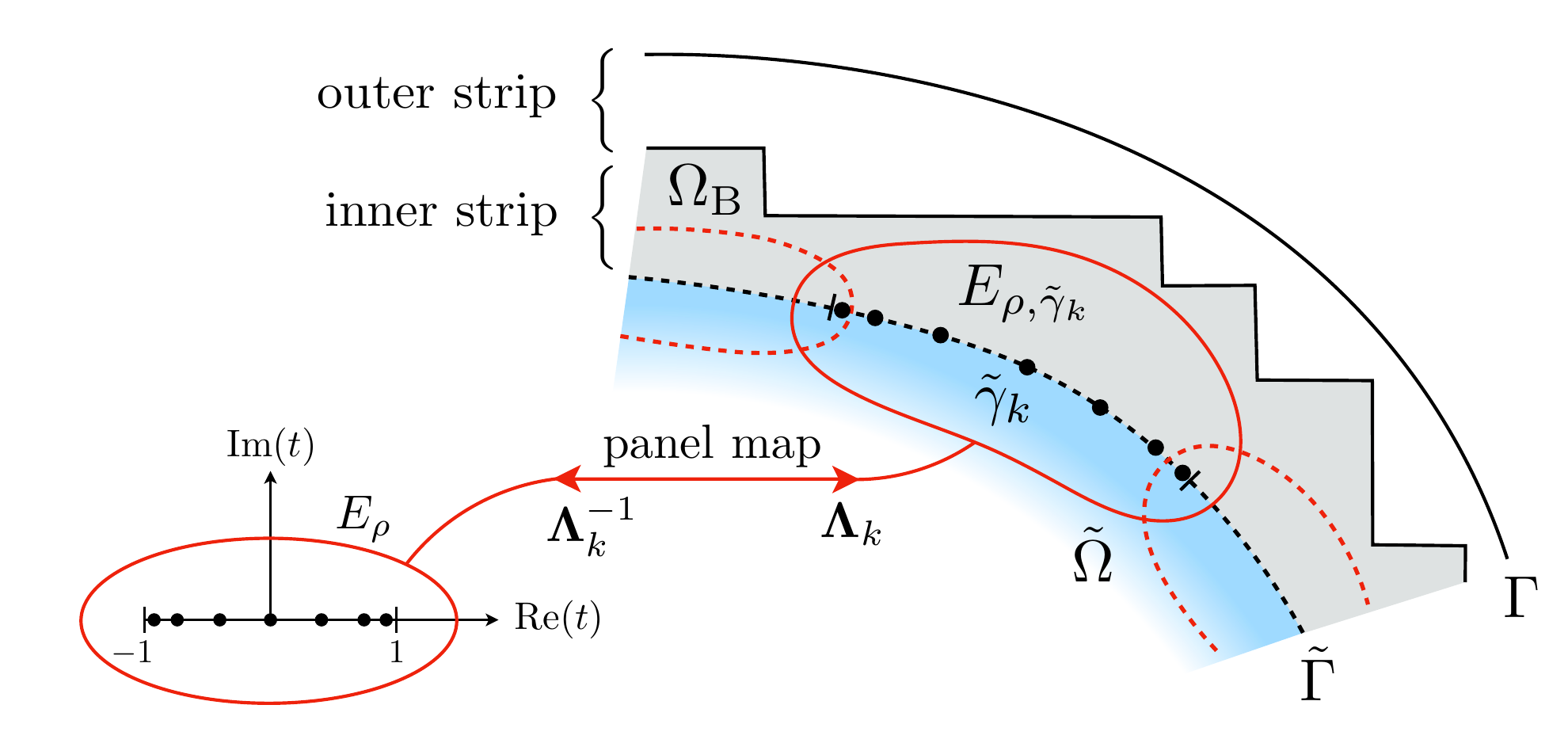}
\end{minipage}%
\hspace{0.4cm}%
\begin{minipage}{0.34\textwidth}
    \vspace{0.2cm}
    \scalebox{0.9}{%
    \begin{overpic}[width=\textwidth]{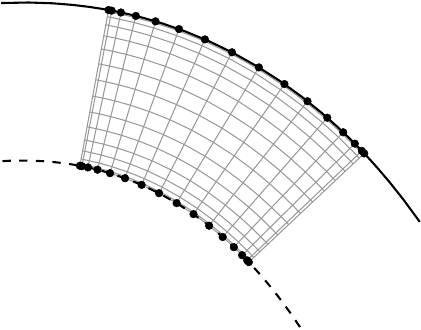}
        \put(46,44) {$\elem_k$}
        \put(1,53) {\color{gray}$\elem_{k+1}$}
        \put(75,18) {\color{gray}$\elem_{k-1}$}
        \put(61,66.5) {$\panel_k$}
        \put(35,21.5) {$\tilde{\panel}_k$}
        \put(-7,75) {$\bdy$}
        \put(-7,40) {$\bdyf$}
        \put(10.2,68) {$\vec{n}_{\mathcal{I}_k}$}
        \put(5.1,63){\tikz \draw[-latex, line width=0.15mm] (0,0) -- (-0.8,0.12);}
        \put(80.8,32.9){\tikz \draw[-latex, line width=0.15mm] (0,0) -- (-0.4,-0.4);}
        \put(80.9,41.8){\tikz \draw[-latex, line width=0.15mm] (0,0) -- (-0.4,0.4);}
        \put(87,47){$\xi$}
        \put(86,33){$\eta$}
        \put(33,-4) {(b)}
        \put(-110,-4) {(a)}
    \end{overpic}%
    }
\end{minipage}
\caption{(a) Zoom near part of the boundary, showing geometry of $k$th panel $\tilde{\panel}_k$ discretizing the fictitious curve $\bdyf$, and its bijection to the standard panel $[-1,1]$. The box code accurately computes the volume potential due to $f$ in $\uncutdom$, whose boundary is well separated from $\bdyf$. $\Ek$ is the Bernstein mapped ellipse for the $k$th panel. (b) The spectral collocation grid in the strip region. Chebyshev nodes on panels $\panel_k$ and $\tilde{\panel}_k$ are connected to form a curvilinear tensor-product grid on which differential operators are numerically constructed. Continuity and continuity of the normal derivative is enforced between elements $\elem_k$ and its neighbors $\elem_{k-1}$ and $\elem_{k+1}$.}
\label{fig:geom}
\end{figure}
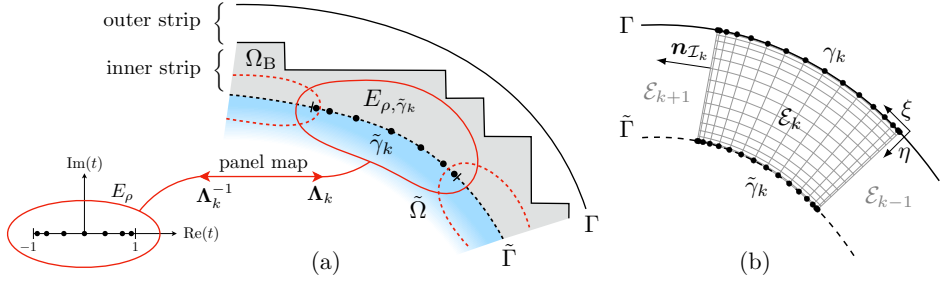

To state the result we need some definitions. Given $\rho>1$, recall (e.g.,~\cite{ATAP}) the open Bernstein ellipse for the standard interval $[-1,1]$,
\be
E_\rho \definedas \left\{(z+z^{-1})/2 \,:\, z\in\C, \; \rho^{-1} < |z| < \rho \right\}.
\label{E}
\ee
We now identify $\R^2$ with $\C$. We use $k = 1, \ldots, \npanel$ to index panels. The $k$th panel is described by a map or chart $X_p : \C \to \C$ such that $X_k([-1,1]) = \tgk$, and $X_k$ is one-to-one and analytic in some open neighborhood of $[-1,1]$. For any $\rho>1$ such that $E_\rho$ lies in this neighborhood, we define the \textit{Bernstein mapped ellipse} for the panel $\tgk$ by
\be
\Ek \definedas X_k(E_\rho).
\label{Ek}
\ee
See~\cref{fig:geom}(a).
For any function $q : \tgk \to \R$ we define its pullback $Q$
such that $Q(t) = q(X(t))$ for $-1\le t\le 1$.
Then define the degree-$n$ Chebyshev approximation $Q_n$ of any function $Q$ on $[-1,1]$
as the usual truncation of its Chebyshev expansion to degree $n$
\cite[Ch.~4]{ATAP}.
Finally, it is a useful shorthand to refer to the Chebyshev approximation of a function on
$\tgk$ as the pushforward of the Chebyshev approximation on $[-1,1]$ of its pullback,
where the pushforward of a function $Q$ simply means $Q \circ X^{-1}$.

The main result shows that, if the source truncation is entirely outside the Bernstein mapped ellipse for a panel, then $\tilde v$ and its first derivatives are analytic on that panel,
as indicated by a specific geometric convergence rate of their Chebyshev approximations.
The latter immediately implies geometric convergence of interpolants or quadrature on the panel.

\begin{thm}\label{t:rate}
  Let $\tilde v = v_\bulk - v$ be the change in particular solution \cref{eq:vbulk} from the classical
  volume potential \cref{eq:vps}.
  Fix a fictitious panel $k\in\{1,\dots,\npanel\}$.
  Let $\tilde v_n$ be the Chebyshev projection of $\tilde v$ on this panel $\tgk$,
  and similarly for $\nabla \tilde v_n$.
  Let $\rho>1$ be such that
  the panel map $X_k$ is analytic and one-to-one in $\overline{E_\rho}$, and
  $\overline{\Ek} \subset \OB$. Then
  \be
  \|\tilde v_n - \tilde v\|_{\infty,\tgk} = \bigO(\rho^{-n})
  \quad \mbox{ and } \quad
  \|\nabla \tilde v_n - \nabla \tilde v\|_{\infty,\tgk} = \bigO(\rho^{-n}),\qquad n\to\infty.
  \label{rate}
  \ee
  In particular, $\tilde v$ and $\nabla \tilde v$ are real analytic on the fictitious panel.
\end{thm}

Combining this theorem with the result that
$v$ is itself smooth on $\bdyf$ (by elliptic regularity
in the interior of $\dom$), then $v_\bulk$ is also smooth on each fictitious panel.
This is our main result for the section. It provides theoretical support for the
high-order convergence observed using the induced panelization of $\bdyf$.
Note that high-order quadtree approximation is also guaranteed, since
$v_\bulk$, being harmonic in $\domf$, must be at least as smooth
in $\domf$ as on $\bdyf$.

\begin{remark}\label{r:rho}
One can lower-bound $\rho$: the construction of $\OB$ in the previous subsection
showed that $\OB$ includes the inner half of the strip.
If strip panels are chosen
with a typical 2:1 aspect ratio, then the Bernstein ellipse preimage for $[-1,1]$ includes
$i/2$, so that $\rho > (1+\sqrt{5})/2 \approx 1.618$.
For the typically used number of upsampled fictitious panel interpolation nodes $n=24$,
the factor $\rho^{-n} \approx 10^{-5}$, implying that
at least 5 correct digits are expected (ignoring unknown prefactors).
This is rather pessimistic: our results in fact show around 10 correct digits,
we believe due to the use of \cref{r:keepcut}.
This in effect pushes the source $\tilde{f}$ out to $\bdy$, doubling the exponential
convergence rate.
\end{remark}

\begin{proof}
  We identify $\R^2$ with $\C$, and thus write $\tilde v = \re V$ for the complex logarithmic potential
  \[
  V(\vec{x}) = \frac{-1}{2\pi}\int_{\C \setminus \uncutdom} L_\vec{y}(\vec{x}) \tilde f(\vec{y}) \, d\vec{y}.
  \]
  Here $L_\vec{y}(\vec{x}) \definedas \log(\vec{x}-\vec{y})$, but with its branch cut in the $\vec{x}$ variable chosen (for each fixed $\vec{y}$) to connect $\vec{y}$ to $\infty$ while avoiding the panel Bernstein mapped ellipse. Since there is some nonzero distance between the compact sets $\supp \tilde{f}$ and $\overline{E_{\rho,\tilde\panel_k}}$, then $L_\vec{y}(\vec{x})$ is uniformly bounded over $y\in\supp \tilde{f}$ and $x\in\overline{\Ek}$.
  Combining this with $\|\tilde{f}\|_{L_1(\R^2)} < \infty$ gives $\|V\|_{\infty,\overline{\Ek}}\le M$ for some $M$. Since $L_\vec{y}(\cdot)$ is analytic in $\overline{\Ek}$, $V$ is also analytic in that set (see, e.g., \cite[Thm.~5.4]{Stein2003}).
  Then $Q(t) = V(X(t))$, the pullback of $V$ to the complex $t$-plane, being the composition of two analytic functions, is analytic in $t \in \overline{E_\rho}$ and bounded by $M$.
  Let $Q_n$ be the Chebyshev truncation of $Q$ on $[-1,1]$. Then,
  \[
  \|\tilde v_n - \tilde v\|_{\infty, \tgk} \le
  \|V_n - V\|_{\infty, \tgk}
  = \|Q_n - Q\|_{\infty, [-1,1]}
  \le  \frac{2M}{\rho-1} \rho^{-n},
  \]
  where the first inequality follows by taking the real part, and the second inequality is a standard approximation theory result~\cite[Thm.~8.2]{ATAP} for functions bounded and analytic in the Bernstein ellipse $E_\rho$. This concludes the proof for $\tilde v$.
  
Finally, since $V$ is analytic on the closed set $\overline{\Ek}$ then $V'$ must also be analytic and bounded on this set. Using $\nabla \tilde v = \re (V', iV')$, we can then apply the above argument replacing $V$ by $V'$, with some other choice for $M$, to show that the two components of $\nabla \tilde v$ obey the same result.
\end{proof}

\begin{remark}[Boundary regularity of classical Newton potential]\label[remark]{r:newt}
  In the classical potential-theory approach \cite{Mayo1992,McKenney1995}, if the
  homogeneous BVP \cref{eq:u_homo} is to have smooth data (allowing its high-order accurate
  solution), then the boundary data of the Newton potential
  \cref{eq:vps} generated by $f\in C^\infty(\dom)$ also needs to be smooth
  (along $\bdy$, since we know that it can only in general be $C^1$ in the normal direction).
  The latter is somewhat of a folk theorem.
  It certainly requires $\bdy$ to be smooth (consider a corner where
  the jump in $f$ induces a weak singularity in $v|_\bdy$).
  We do not know of literature stating the result, but note the following%
  \footnote{The argument is due to Leslie Greengard, personal communication.}.
  The Newton potential $v$ in $\R^2$ is equal to the interior
  solution (extended by zero outside of $\dom$) to $-\Delta u = f$ with $u=0$ on $\bdy$,
  minus the single-layer potential $\SL_\bdy [\partial_\vec{n} u]$.
  By regularity up to the boundary \cite[\S6.3.2, Thm.~6]{EvansPDE},
  $\partial_\vec{n} u \in C^\infty(\bdy)$, and
  since the single-layer boundary integral operator on a smooth surface $\bdy$
  is one order smoothing \cite[Thm.~7.2]{mclean}, then $v\in C^\infty(\bdy)$.
  Alternatively, proof sketches directly tackling \cref{eq:vps} exist
  \cite{tao23} (also see \cite{hu2024} which needs $\dom$ convex).
\end{remark}

\subsection{The strip problem}\label{sec:strip}

The volume potential formulation above yields a particular solution $v_\bulk$ valid everywhere inside $\domf$. We now turn to the problem of constructing a particular solution $v_\strip$ in the remaining boundary-fitted region $\stripdom$, such that
\begin{equation}\label{eq:strip_problem}
\begin{aligned}
\lap v_\strip(\vec{x}) &= f(\vec{x}), &&\quad \vec{x} \in \stripdom, \\
v_\strip(\vec{x}) &= 0, &&\quad \vec{x} \in \bdy \cup \bdyf.
\end{aligned}
\end{equation}
Note that the region $\stripdom$ is thin, as its width was chosen so that the distance between corresponding nodes on $\bdy$ and $\bdyf$ is on the order of the local panel size. Therefore, we choose to solve in this region using a multidomain spectral collocation method, with each panel creating an element that spans the entirety of the thickness of $\stripdom$ between $\bdy$ and $\bdyf$. Specifically, let $\elem_k$ be the region bounded by $\bdy$, $\bdyf$, and the straight lines connecting the left and right endpoints of panel $\panel_k$ with the left and right endpoints of panel $\tilde{\panel}_k$, respectively. See~\cref{fig:geom}(b) for reference. Denote by $\mathcal{I}_k$ the interface between elements $\elem_k$ and $\elem_{k+1}$, and $\vec{n}_{\mathcal{I}_k}$ the unit normal vector to $\mathcal{I}_k$ pointing from $\elem_k$ to $\elem_{k+1}$. (As the strip is periodic in the annular direction, we let $\mathcal{I}_{\npanel}$ be the interface between elements $\elem_{\npanel}$ and $\elem_1$, and $\vec{n}_{\mathcal{I}_\npanel}\!$ its corresponding normal vector). The multidomain boundary value problem is then formulated as
\begin{equation}\label{eq:multidomain}
\begin{aligned}
\lap v_\strip^k(\vec{x}) &= f(\vec{x}), &&\quad \vec{x} \in \elem_k, \\
v_\strip^k(\vec{x}) &= 0, &&\quad \vec{x} \in \partial \elem_k \cap \left(\bdy \cup \bdyf\right), \\
v_\strip^k(\vec{x}) &= v_\strip^{k+1}(\vec{x}), &&\quad \vec{x} \in \mathcal{I}_k, \\
\tfrac{v_\strip^k}{\partial \vec{n}_{\mathcal{I}_k}}(\vec{x}) &= \tfrac{v_\strip^{k+1}}{\partial \vec{n}_{\mathcal{I}_k}}(\vec{x}), &&\quad \vec{x} \in \mathcal{I}_k,
\end{aligned}
\end{equation}
for $k = 1, \ldots, \npanel$, where $v_\strip^k$ is the solution on element $\elem_k$ (with $v_\strip^{\npanel+1} \definedas v_\strip^1$).

We use a spectral collocation method to discretize \cref{eq:multidomain}. Let $I_\text{leg}^\text{cheb}$ be the $(p+1) \times (p+1)$ interpolation matrix which maps function values at $p+1$ Gauss--Legendre nodes to function values at $p+1$ second-kind Chebyshev nodes~\cite{ATAP}. On panel $k$, the nodes $\vec{x}^\text{cheb}_{:,k} = I_\text{leg}^\text{cheb} \vec{x}_{:,k}$ are then Chebyshev nodes. Chebyshev nodes on $\bdyf$ may be similarly defined. Letting $r^\text{cheb}$ be the order-$p$ second-kind Chebyshev nodes on $[-1,1]$, a curvilinear tensor product grid of nodes $\vec{X}_{ij,k}$ for element $\elem_k$ may be constructed as
\[
\vec{X}_{ij,k} = \left(\frac{1 + r^\text{cheb}_i}{2}\right) \tilde{\vec{x}}^\text{cheb}_{j,k} + \left(\frac{1 - r^\text{cheb}_i}{2}\right) \vec{x}^\text{cheb}_{j,k}.
\]
Let $\{\vec{\xi}_{ij} = (\xi_{ij},\eta_{ij})\}_{i,j=1}^{p+1}$ be the set of tensor-product second-kind Chebyshev nodes of order $p$ over the reference square $[-1,1]^2$. Then for each element, the nodes $\vec{X}_{ij,k}$ numerically define a mapping from the reference square $[-1,1]^2$ to $\elem_k$. That is, the coordinate mapping for each element is a function $\vec{X}_k(\xi,\eta) = (X_k(\xi,\eta), Y_k(\xi,\eta)) : [-1,1]^2 \to \R^2$ such that $\vec{X}_k(\xi_{ij}, \eta_{ij}) = \vec{X}_{ij,k}$ for $i,j = 1, \ldots, p+1$. We may numerically approximate the coordinate mapping through interpolation at the nodes via
\[
\vec{X}_k(\xi,\eta) \approx \sum_{i=1}^{p+1} \sum_{j=1}^{p+1} \vec{X}_{ij,k} \, \ell_j(\xi) \, \ell_i(\eta), \qquad (\xi,\eta) \in [-1,1]^2,
\]
where $\ell_j$ is the $j$th Lagrange polynomial associated with the second-kind Chebyshev nodes. Partial derivatives of the elemental coordinate maps, $\partial \vec{X}_k / \partial \xi$ and $\partial \vec{X}_k / \partial \eta$, may then be computed through numerical spectral differentiation~\cite{Trefethen2000}, and derivatives of the inverse coordinate mappings may be derived via the chain rule,
\begin{equation}\label{eq:inverse_coord}
\begin{aligned}
\tfrac{\partial \xi}{\partial x}  &= \hphantom{-} \tfrac{1}{J_k} \tfrac{\partial Y_k}{\partial \eta}, \qquad &\tfrac{\partial \xi}{\partial y} &= - \tfrac{1}{J_k} \tfrac{\partial X_k}{\partial \eta}, \\
\tfrac{\partial \eta}{\partial x} &= - \tfrac{1}{J_k} \tfrac{\partial Y_k}{\partial \xi}, \qquad &\tfrac{\partial \eta}{\partial y} &= \hphantom{-} \tfrac{1}{J_k} \tfrac{\partial X_k}{\partial \xi},
\end{aligned}
\end{equation}
where $J_k = \frac{\partial X_k}{\partial \xi} \frac{\partial Y_k}{\partial \eta} - \frac{\partial X_k}{\partial \eta} \frac{\partial Y_k}{\partial \xi}$ is the Jacobian of the coordinate mapping $\vec{X}_k$.

The spectral collocation method proceeds by discretizing the differential operator on each element and enforcing the PDE and boundary conditions on its interior and boundary nodes, respectively. For each $k$, denote by $v_{ij,k} \approx v_\strip^k(\vec{X}_{ij,k})$; that is, $v_{ij,k}$ is simply $v_\strip^k$ sampled on the grid of element $\elem_k$. Then we may approximate each function by
\[
v_\strip^k(\xi,\eta) \approx \sum_{i=1}^{p+1} \sum_{j=1}^{p+1} v_{ij,k} \, \ell_j(\xi) \, \ell_i(\eta), \qquad (\xi,\eta) \in [-1,1]^2,
\]
for $k = 1, \ldots, \npanel$, where we have introduced the slight abuse of notation $v_\strip^k(\xi,\eta) \allowbreak\definedas v_\strip^k(\vec{X}_k(\xi,\eta))$. Now, let $D \in \C^{(p+1) \times (p+1)}$ be the one-dimensional spectral differentiation matrix associated with the order-$p$ second-kind Chebyshev nodes on the interval $[-1,1]$~\cite{Trefethen2000}, and let $I \in \C^{(p+1) \times (p+1)}$ be the identity matrix. Then $D_\xi = D \otimes I$ and $D_\eta = I \otimes D$ are the two-dimensional differentiation matrices in the $\xi$- and $\eta$-directions on the reference square, of size $(p+1)^2 \times (p+1)^2$. Let $M[v] \in \C^{(p+1) \times (p+1)}$ denote the diagonal multiplication matrix formed by placing the entries of $v_{ij}$ along the diagonal. Using \cref{eq:inverse_coord}, one may show that differentiation matrices in the $x$- and $y$-directions are given by
\[
\begin{aligned}
D_{X_k} &= M\!\!\left[\tfrac{\partial \xi}{\partial x}\right] D_\xi + M\!\!\left[\tfrac{\partial \eta}{\partial x}\right] D_\eta, \\
D_{Y_k} &= M\!\!\left[\tfrac{\partial \xi}{\partial y}\right] D_\xi + M\!\!\left[\tfrac{\partial \eta}{\partial y}\right] D_\eta.
\end{aligned}
\]
The discrete Laplacian on element $\elem_k$ is then given by $\lap \approx (D_{X_k})^2 + (D_{Y_k})^2$.

As the multidomain formulation only couples elements to their left and right neighbors, the resulting linear system is block tridiagonal, aside from a corner block due to the periodicity of the strip region. Therefore, direct matrix inversion via block banded LU factorization---along with the Woodbury formula to correct for the corner block---takes only $\mathcal{O}(p^3\npanel)$ operations to compute the strip solutions $v_\strip^k$ for $k=1,\ldots,\npanel$.

\begin{remark}\label{rem:strip_res}
It may happen that $f(\vec{x})$ is unresolved on the strip grid induced by the given panelization of $\bdy$. To handle this case, our solver first checks if $f(\vec{x})$ is resolved on each element of the strip to the given tolerance; if $f$ is unresolved on some elements, the solver splits the corresponding panels in $\bdy$ and the algorithm restarts. This process is recursive and happens automatically at the start of \cref{alg:solver}, before quadtree construction.
\end{remark}

\subsection{Patching together \texorpdfstring{\textnormal{$v_\bulk$}}{vbulk} and \texorpdfstring{\textnormal{$v_\strip$}}{vstrip}}\label{sec:glue}

We now have particular solutions in both regions of $\dom$: $v_\bulk$ satisfying $\lap v_\bulk(\vec{x}) = f(\vec{x})$ for $\vec{x} \in \domf$, and $v_\strip$ satisfying $\lap v_\strip(\vec{x}) = f(\vec{x})$ for $\vec{x} \in \stripdom$. However, the piecewise function
\begin{equation}
v(\vec{x}) = \begin{cases}
  v_\bulk(\vec{x}), \quad \vec{x} \in \domf, \\
  v_\strip(\vec{x}), \quad \vec{x} \in \stripdom,
\end{cases}
\end{equation}
is not a globally smooth particular solution in $\dom$, since for each $\vec{y} \in \bdyf$,
\[
\lim_{\vec{x} \to \vec{y}^+} v(\vec{x}) \neq \lim_{\vec{x} \to \vec{y}^-} v(\vec{x}) \quad \text{ and } \quad \lim_{\vec{x} \to \vec{y}^+} \partial_{\vec{n}_\vec{x}} v(\vec{x}) \neq \lim_{\vec{x} \to \vec{y}^-} \partial_{\vec{n}_\vec{x}} v(\vec{x}),
\]
where superscripts of $-$ and $+$ denote limits taken from the interior and exterior of the domain, respectively. That is, the values and normal derivatives of $v_\bulk$ and $v_\strip$ do not match across the interface $\bdyf$.

Denote by $\SL_\bdyf[\sigma]$ and $\DL_\bdyf[\tau]$ the Laplace single and double layer potentials induced by the densities $\sigma$ and $\tau$, respectively, on the boundary $\bdyf$, given by
\[
\SL_\bdyf[\sigma](\vec{x}) \definedas \int_\bdyf \Phi(\vec{x},\vec{y}) \sigma(\vec{y}) \, d\vec{y}, \qquad
\DL_\bdyf[\tau](\vec{x}) \definedas \int_\bdyf \frac{\partial \Phi(\vec{x},\vec{y})}{\partial \vec{n}_\vec{y}} \tau(\vec{y}) \, d\vec{y},
\]
Such layer potentials are harmonic functions, satisfying $\lap \SL_\bdyf[\sigma] = 0$ and $\lap \DL_\bdyf[\sigma] = 0$ in all of $\R^2 \setminus \bdyf$. It can be shown that $\SL_\bdyf$ and $\DL_\bdyf$ satisfy the jump relations~\cite[Ch.~6]{Kress2014a}
\[
\begin{gathered}
\lim_{\vec{x} \to \vec{y}^+}\!\SL_\bdyf[\sigma](\vec{x}) - \!\!\lim_{\vec{x} \to \vec{y}^-}\!\SL_\bdyf[\sigma](\vec{x}) = 0, \\
\lim_{\vec{x} \to \vec{y}^+}\! \partial_{\vec{n}_\vec{x}} \SL_\bdyf[\sigma](\vec{x}) - \!\!\lim_{\vec{x} \to \vec{y}^-}\! \partial_{\vec{n}_\vec{x}} \SL_\bdyf[\sigma](\vec{x}) = -\sigma, \\
\end{gathered}
\]
and
\[
\begin{gathered}
\lim_{\vec{x} \to \vec{y}^+}\!\DL_\bdyf[\tau](\vec{x}) - \!\!\lim_{\vec{x} \to \vec{y}^-}\!\DL_\bdyf[\tau](\vec{x}) = \tau, \\
\lim_{\vec{x} \to \vec{y}^+}\! \partial_{\vec{n}_\vec{x}} \DL_\bdyf[\tau](\vec{x}) - \!\!\lim_{\vec{x} \to \vec{y}^-}\! \partial_{\vec{n}_\vec{x}} \DL_\bdyf[\tau](\vec{x}) = 0,
\end{gathered}
\]
for each $\vec{y} \in \bdyf$. The single layer potential is continuous across its boundary, with a jump in normal derivative equal to the negative of the given density. Similarly, the double layer potential has continuous normal derivative across its boundary, with a jump in value equal to the given density. Thus, we set
\[
\begin{aligned}
\tau &= v_\bulk|_\interface - v_\strip|_\interface, \\
\sigma &= \partial_\vec{n} v_\bulk|_\interface - \partial_\vec{n} v_\strip|_\interface,
\end{aligned}
\]
and define the function
\[
v_\glue(\vec{x}) = \SL_\bdyf[\sigma](\vec{x}) - \DL_\bdyf[\tau](\vec{x}),
\]
for all $\vec{x} \in \dom$. Adding this function to the piecewise-defined particular solution above results in a globally smooth particular solution to \cref{eq:ipde}, given by
\begin{equation}
v(\vec{x}) = \begin{cases}
  v_\bulk(\vec{x}) + v_\glue(\vec{x}), \quad \vec{x} \in \domf, \\
  v_\strip(\vec{x}) + v_\glue(\vec{x}), \quad \vec{x} \in \stripdom,
\end{cases}
\end{equation}
with values on $\bdyf$ defined by their limit from either side
(which are equal, to discretization accuracy).
As the layer potentials in $v_\glue$ may be rapidly evaluated using the FMM, the overall particular solution $v$ may be rapidly evaluated at any point $\vec{x} \in \dom$.

\section{Numerical results}\label{sec:results}

We now demonstrate our adaptive Poisson solver on some challenging geometries, requiring adaptivity both along the boundary to resolve geometric features and in the bulk to resolve spatial inhomogeneities. All numerical examples were run in MATLAB R2024b on single core of an M4 Max MacBook Pro with \qty{128}{\giga\byte} of memory. Our code is open source and freely available~\cite{GithubRepo}.

\subsection{A rounded raindrop}\label{sec:raindrop}

\begin{figure}[H]
	\centering
	\vspace{2em}
	\begin{tikzpicture}[overlay,xshift=1.17cm]
		\draw (1.3,4.97) -- (2,5.609);
		\draw (1.3,4.88) -- (2,4.491);
	\end{tikzpicture}%
	\begin{tikzpicture}[overlay,xshift=4.82cm]
		\draw (1.3,4.97) -- (2,5.609);
		\draw (1.3,4.88) -- (2,4.491);
	\end{tikzpicture}%
	\begin{tikzpicture}[overlay,xshift=8.5cm]
		\draw (1.3,4.97) -- (2,5.609);
		\draw (1.3,4.88) -- (2,4.491);
	\end{tikzpicture}%
	\begin{overpic}[height=5cm]{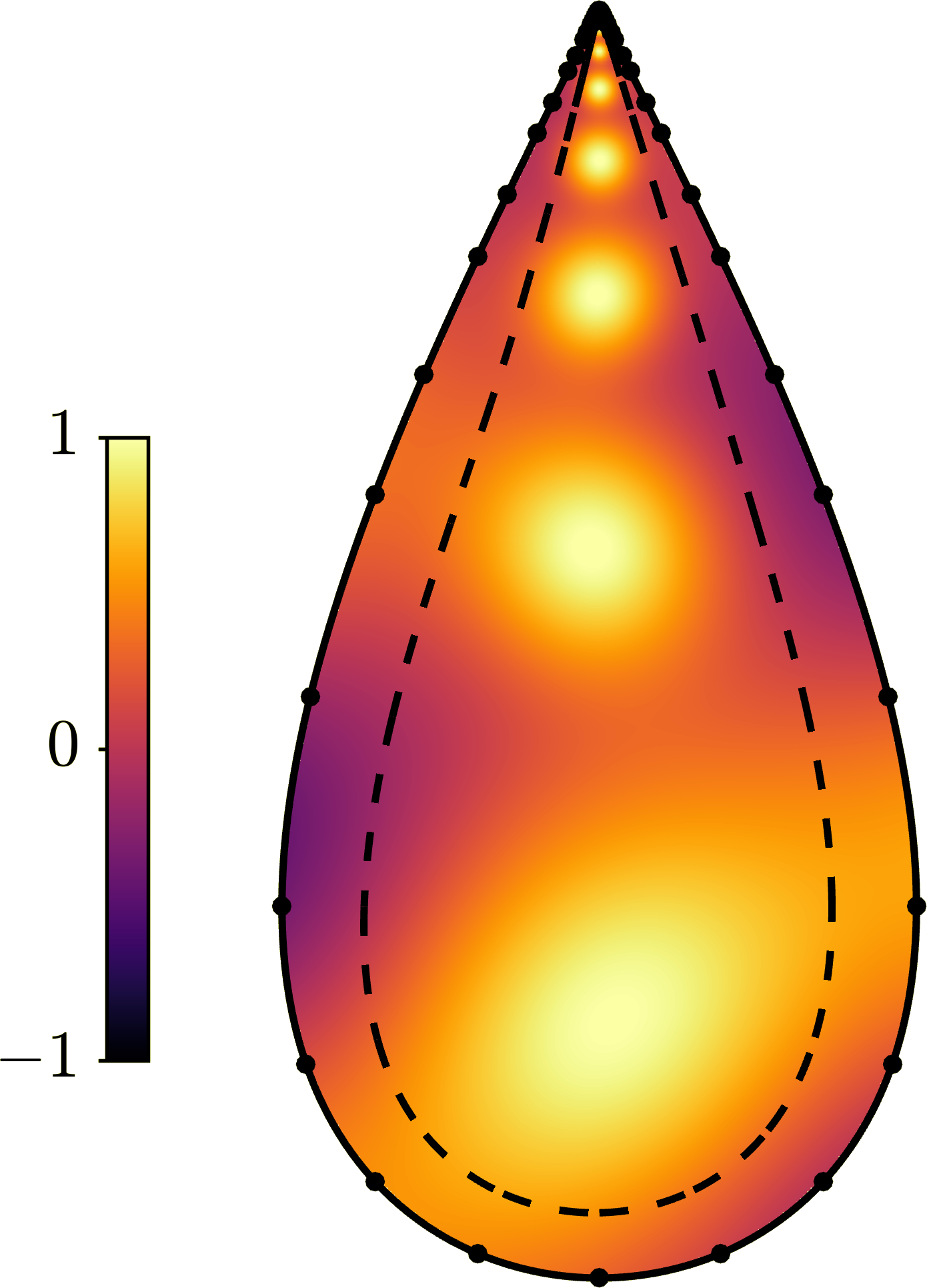}
		\put(61,90) {\includegraphics[height=1.1cm,trim={0 0.2cm 0 -0.3cm},clip,frame]{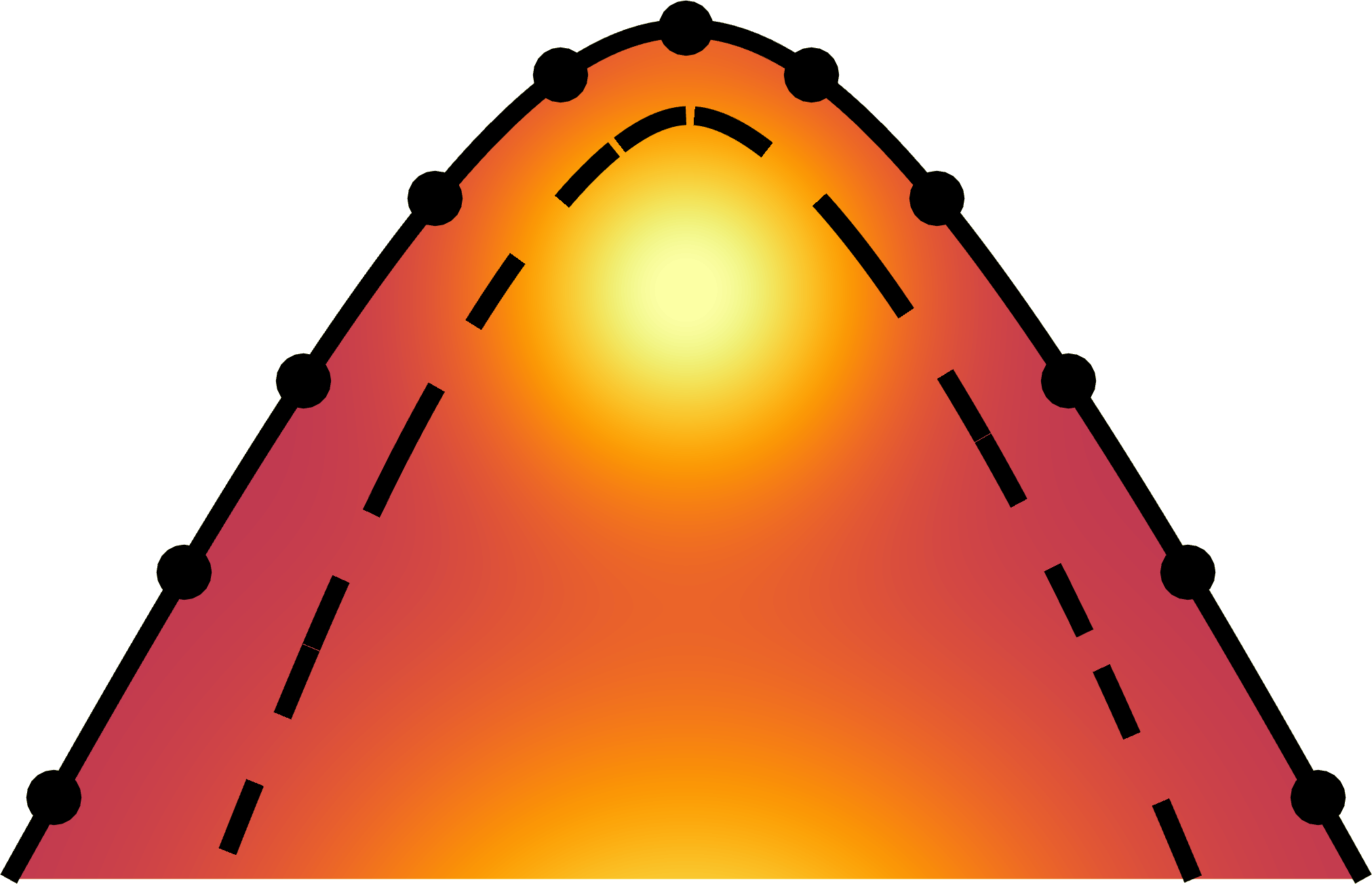}}
		\put(71,84) {\footnotesize$1000$\scriptsize$\mspace{-1mu}\times$}
	\end{overpic}
	\hspace{0.8cm}
	\begin{overpic}[height=5cm]{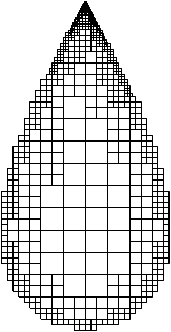}
		\put(40,90) {\includegraphics[height=1.1cm,trim={0 0.2cm 0 -0.3cm},clip,frame]{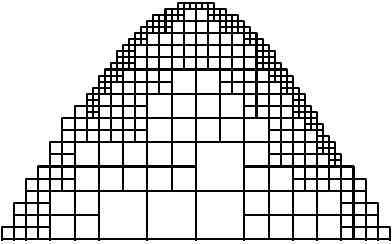}}
	\end{overpic}
	\hspace{0.8cm}
	\begin{overpic}[height=5cm]{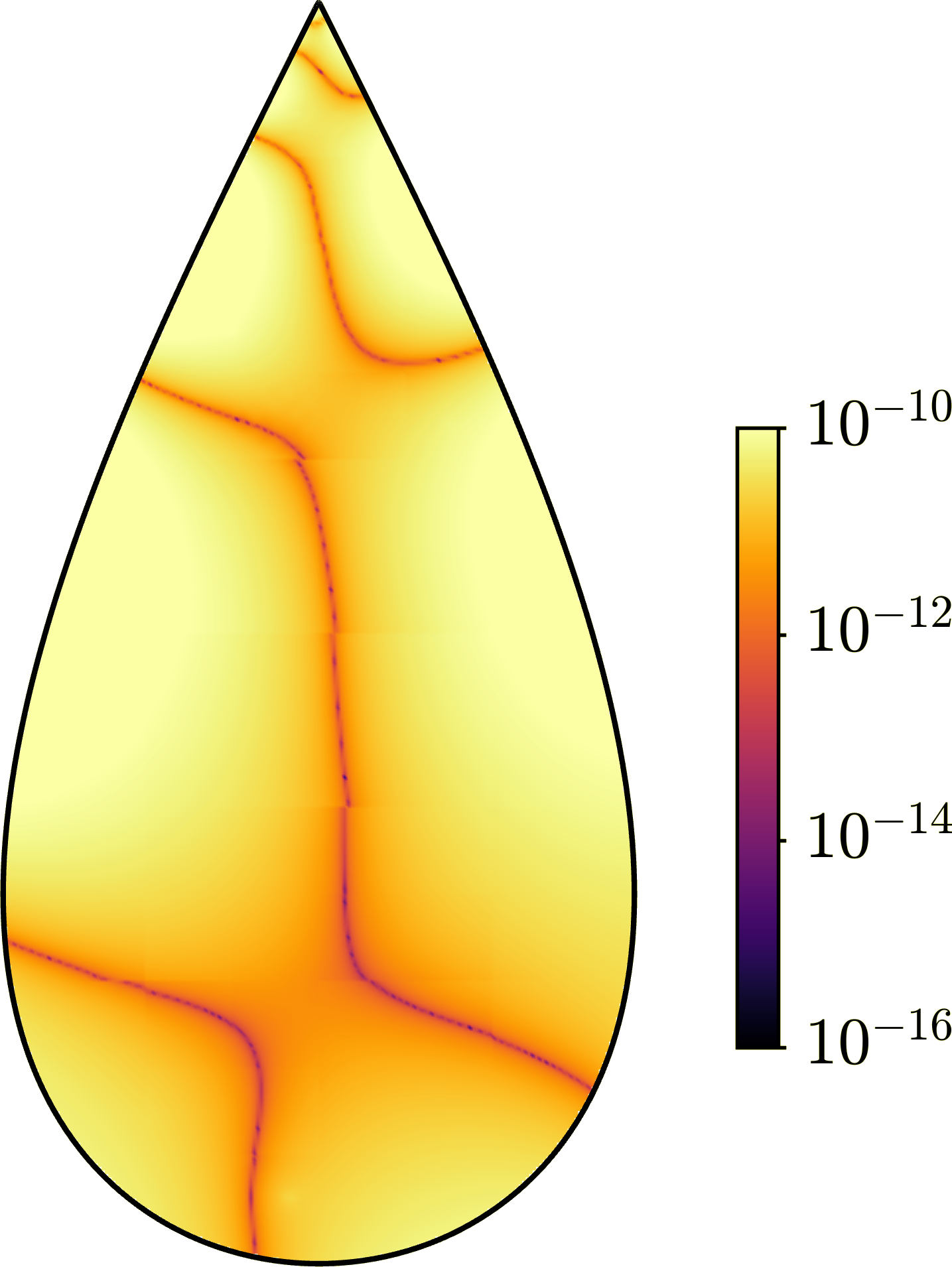}
		\put(39.5,90) {\includegraphics[height=1.1cm,trim={0 0.2cm 0 -0.3cm},clip,frame]{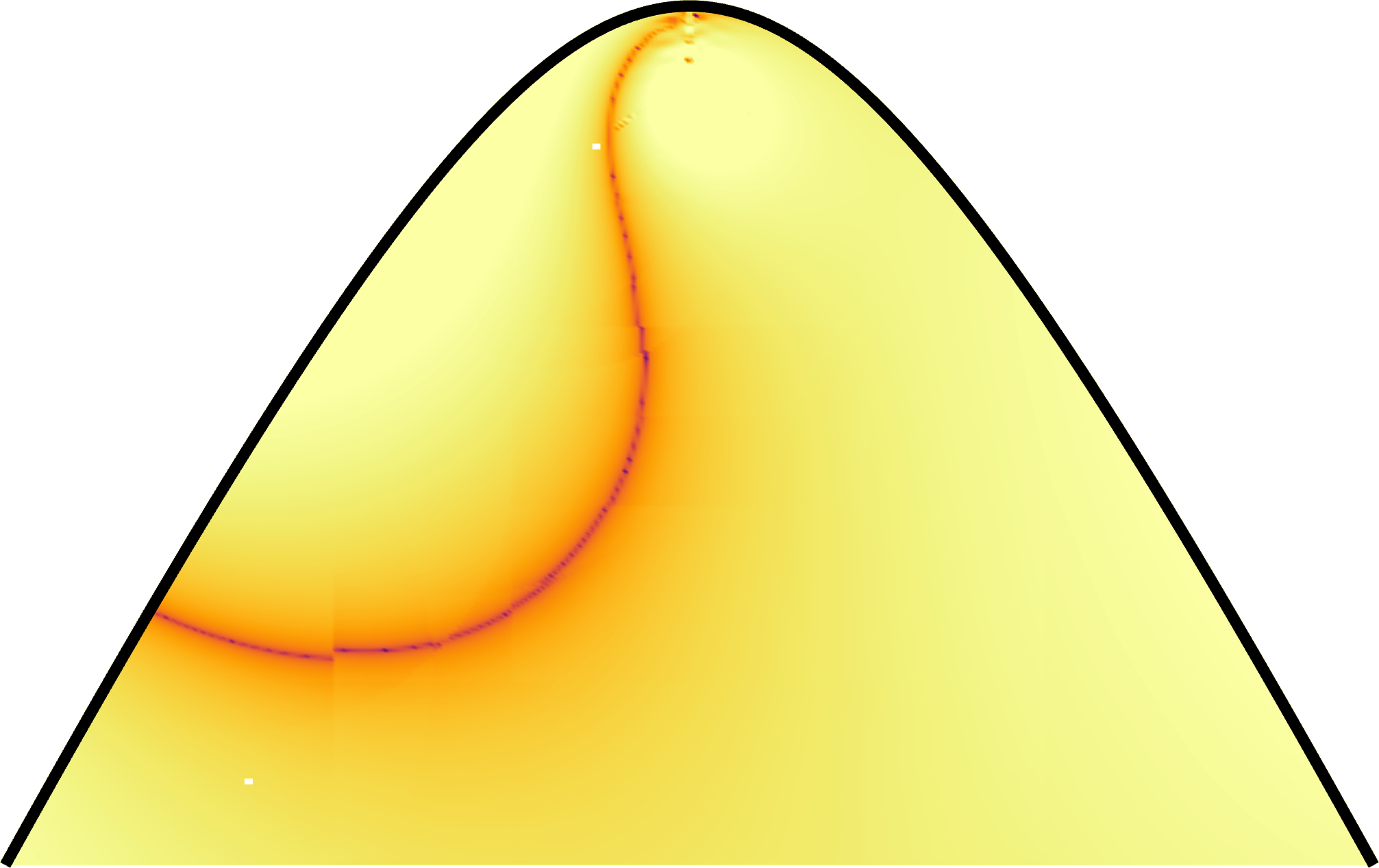}}
	\end{overpic}%
	\caption{(Left) The raindrop shape is adaptively panelized, with small panels clustering in the rounded cusp. Black circles correspond to panel endpoints. The strip region (plotted as a dashed line) conforms to this adaptive panelization. The test solution is also shown and consists of a series of Gaussians with decreasing variances clustering inside the cusp. (Center) The inhomogeneity induced by the given solution is adaptively resolved on a 16th-order background quadtree with \num{594944} unknowns and truncated according to the strip refinement criterion outlined in~\cref{sec:quadtree}. (Right) The maximum absolute error is around $10^{-10}$ over the whole domain.}
	\label{fig:raindrop}
\end{figure}
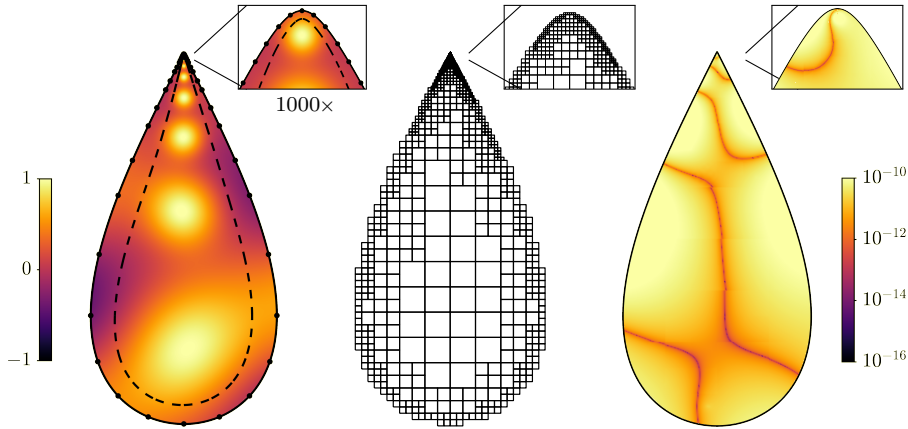

We first demonstrate our adaptive Poisson solver on a raindrop-shaped geometry shown in~\cref{fig:raindrop}, consisting of 56 panels of order 16. The corner of the raindrop is rounded to a length scale of $10^{-3}$. While the pinched end of the raindrop-like shape would force a uniform-grid method to over-refine the largest length scales, our variable-width strip region smoothly captures the transition in panel size across three orders of magnitude.

We run the solver with a requested tolerance of $10^{-10}$. Our test solution consists of a sum of Gaussians exponentially clustering into the cusp with decreasing variance, along with the smooth background function $2 x \cos{3\pi y}$ added for variation along the boundary. We analytically compute the inhomogeneity $f$ corresponding to this solution. The precomputation phase for this domain, which includes constructing the strip region and building a 16th-order quadtree to satisfy the refinement criterion, takes \qty{0.3}{\second}. The resulting quadtree possesses \num{2965} leaf nodes with \num{594944} degrees of freedom, and computation of $v_\bulk$ using a 16th-order box code~\cite{boxcode2d} takes \qty{0.06}{\second}. The strip region contains 56 elements, with each element upsampled slightly to a resolution of $16 \times 24$ to yield a strip mesh with \num{21504} degrees of freedom. Computation of the strip particular solution, $v_\strip$, takes \qty{0.3}{\second}. The homogeneous solution $w$ is computed in \qty{0.2}{\second} by solving a boundary integral equation using GMRES, with matrix-vector products accelerated by the FMM. Evaluating the solution back on all quadtree and strip nodes using the FMM takes \qty{2.7}{\second}.

\subsection{Comparison to FFT-based schemes}\label{sec:fft}

We now compare our adaptive scheme against an FFT-based uniform bulk solver. To do this, we define a series of problems with decreasing length scale $\eta$, driven by both geometry and righthand side. We generalize the raindrop shape from~\cref{sec:raindrop} by rounding the corner to a length scale of $\eta$. Such a curve is parametrized in $\C$ by
\[
z(t) = -\tfrac14 \sin{t} - \mathrm{i}\sqrt{\sin^2 \tfrac{t}{2} + \eta^2}, \qquad t \in [0, 2\pi].
\]
Its minimum radius of curvature is $R_\tbox{min} = \eta/4 + \bigO(\eta^2)$, at $t=0$. As in~\cref{sec:raindrop}, we set the righthand side to a series of Gaussians clustering into the corner with decreasing variance, with the numerical support of the narrowest Gaussian proportional to the length scale $\eta$.

\begin{figure}[H]
\definecolor{gridcolor}{rgb}{0.88,0.88,0.88}
\pgfplotsset{
    width=0.47\textwidth,
    height=0.42\textwidth,
    ylabel style={yshift=-3},
    xmin=1e0,
    xmax=1e8,
    ymin=1e-2,
    ymax=1e4,
    axis line style={line width=0.5},
    grid=major,
    grid style={gridcolor},
    ticklabel style={font=\footnotesize},
    max space between ticks=20,
    xlabel=$1/\eta$,
    xtick={1e0,1e1,1e2,1e3,1e4,1e5,1e6,1e7,1e8},
    xticklabels={$10^0$,,$10^2$,,$10^4$,,$10^6$,,$10^8$},
    every major tick/.append style={black, major tick length=0.10cm, line width=0.5},
    every minor tick/.append style={black, minor tick length=0.05cm, line width=0.5},
    legend style={at={(-0.23,1.12)},anchor=south,legend columns=-1,/tikz/every even column/.append style={column sep=0.3cm}},
    legend cell align={left}
}
\begin{tikzpicture}
\begin{loglogaxis}[name=plot1, ylabel={Time (s)}]
\addplot[color=red,mark=*,mark size=1.3,line width=1] coordinates {
    (1.600000e+00, 4.0792e-05)
    (3.200000e+00, 3.8167e-05)
    (6.400000e+00, 3.9042e-05)
    (1.280000e+01, 9.5583e-05)
    (2.560000e+01, 5.7312e-04)
    (5.120000e+01, 2.2996e-03)
    (1.024000e+02, 1.0183e-02)
    (2.048000e+02, 4.4377e-02)
    (4.096000e+02, 1.8316e-01)
    (8.192000e+02, 8.3880e-01)
    (1.638400e+03, 4.2393e+00)
    (3.276800e+03, 2.0258e+01)
    (6.553600e+03, 8.9545e+01)
};
\addplot[color=red,dashed,line width=1,domain=7.053600e+03:1e7]{2.5e-7*x^2*ln(x)};
\addplot[color=blue,mark=square*,mark size=1.2,line width=1] coordinates {
    (1.000000e+00, 2.6865e-01)
    (1.000000e+01, 4.4481e-01)
    (1.000000e+02, 6.1801e-01)
    (1.000000e+03, 8.1031e-01)
    (1.000000e+04, 1.0064e+00)
    (1.000000e+05, 1.2034e+00)
    (1.000000e+06, 1.4082e+00)
    (1.000000e+07, 1.6776e+00)
    (1.000000e+08, 1.8561e+00)
};
\end{loglogaxis}
\begin{loglogaxis}[at={($(plot1.east)+(2cm,0)$)}, anchor=west, ylabel={Memory (GB)}]
\addplot[color=red,mark=*,mark size=1.3,line width=1] coordinates {
    (1.600000e+00, 8.192000e-06)
    (3.200000e+00, 3.276800e-05)
    (6.400000e+00, 1.310720e-04)
    (1.280000e+01, 5.242880e-04)
    (2.560000e+01, 2.097152e-03)
    (5.120000e+01, 8.388608e-03)
    (1.024000e+02, 3.355443e-02)
    (2.048000e+02, 1.342177e-01)
    (4.096000e+02, 5.368709e-01)
    (8.192000e+02, 2.147484e+00)
    (1.638400e+03, 8.589935e+00)
    (3.276800e+03, 3.435974e+01)
    (6.553600e+03, 1.374390e+02)
};
\addplot[color=red,dashed,line width=1,domain=7.053600e+03:1e7]{3.6e-7*x^2*ln(x)};
\addplot[color=blue,mark=square*,mark size=1.2,line width=1] coordinates {
    (1.000000e+00, 1.279181e-02)
    (1.000000e+01, 1.202586e-02)
    (1.000000e+02, 1.440563e-02)
    (1.000000e+03, 1.687552e-02)
    (1.000000e+04, 1.955430e-02)
    (1.000000e+05, 2.194227e-02)
    (1.000000e+06, 2.424832e-02)
    (1.000000e+07, 2.688205e-02)
    (1.000000e+08, 2.927821e-02)
};
\legend{\texttt{fft2},,Proposed solver}
\end{loglogaxis}
\end{tikzpicture}
\vspace{-1em}
\caption{Performance comparison between the 2D FFT (red circles) and our adaptive solver (blue squares). We benchmark our adaptive solver on a series of teardrop problems with decreasing length scale $\eta$ and a requested error tolerance of $10^{-10}$, and record runtime and memory consumption of the entire solver. For comparison, we benchmark the 2D FFT on successively refined uniform grids using MATLAB's {\upshape\texttt{fft2}}, using ten gridpoints to resolve $\eta$; see \cref{sec:fft}.
Dashed lines indicate predicted values for the FFT. Note that our adaptive solver is solving the full Poisson problem, while the 2D FFT would only solve the bulk problem.\label{f:fft}}
\end{figure}
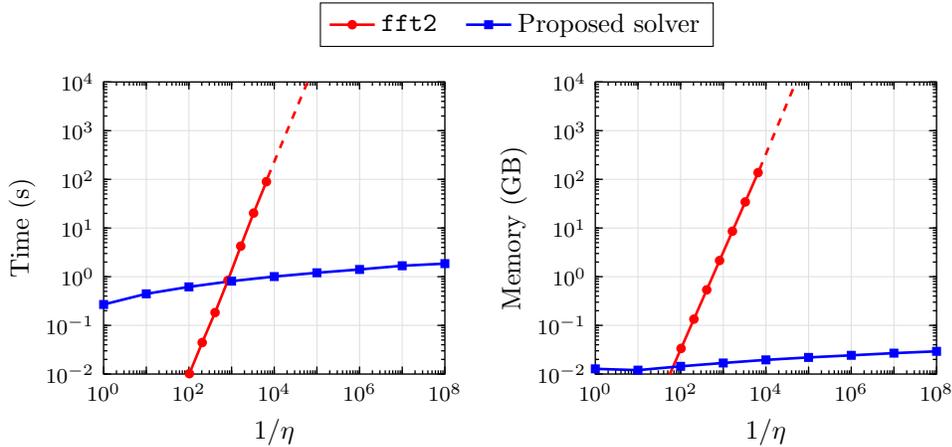

In \cref{f:fft} we vary $\eta$ from $10^0$ to $10^{-8}$ and plot the total runtime and memory consumption of our adaptive solver applied to each problem with a requested tolerance of $10^{-10}$. For comparison, we also run the 2D FFT (as implemented in MATLAB's \texttt{fft2}) on a series of successively refined $n \times n$ uniform grids and plot the same, with $n$ a power of two to allow the fastest FFTs.
The length scale parameter $\eta$ that the $n\times n$ grid for any variety of FFT-based solver
would be able to resolve cannot shrink faster than $\bigO(1/n)$.  
We choose a specific relation $\eta = 10h$, where $h=1/n$ is the grid spacing; in other words
$R_\tbox{min} = 2.5\, h$.
The figure shows the stark contrast between the runtime or memory use of our adaptive solver versus that of any nonadaptive
FFT-based solver:
runtime and memory grow very weakly with $1/\eta$ (one expects logarithmically) for the adaptive case,
while they are both $\bigO(1/\eta^2)$ for the FFT.
The upshot is that on a single shared memory node one cannot reach $\eta<10^{-4}$, whereas the adaptive solver
easily reaches $\eta=10^{-8}$ in a few seconds and a few tens of MB of memory.
Note that this $R_\tbox{min}$ is still about eight times smaller than needed in the FFT spectral solver of the second author \cite{Stein2022b}, whose Figure~7 shows that $R_\tbox{min}/h\approx 20$ is needed for 10-digit accuracy.
This prefactor would move the red FFT lines to the left, only strengthening our point. 

\subsection{A close-to-touching multiscale geometry and inhomogeneity}\label{sec:gobbler}

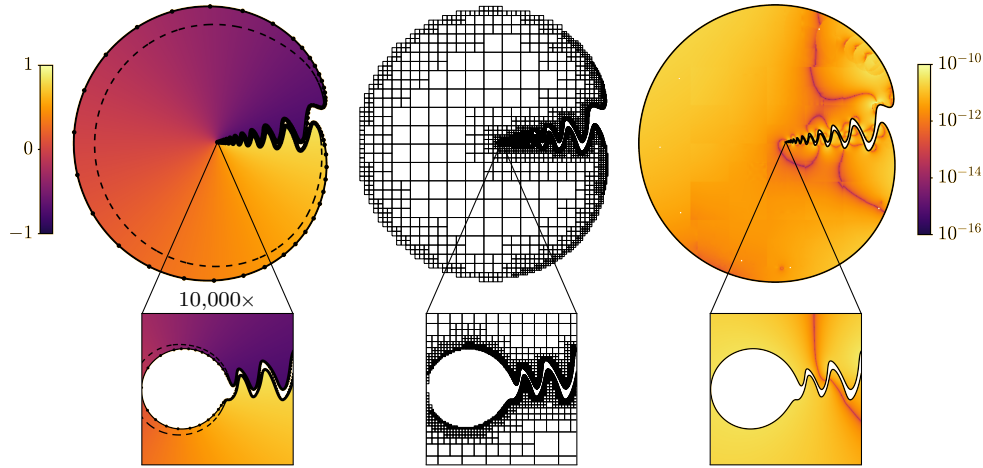
\begin{figure}[htb]
	\centering
	\begin{tikzpicture}
        \node[anchor=south west,inner sep=0] (image) at (0,0) {\includegraphics[height=3.7cm]{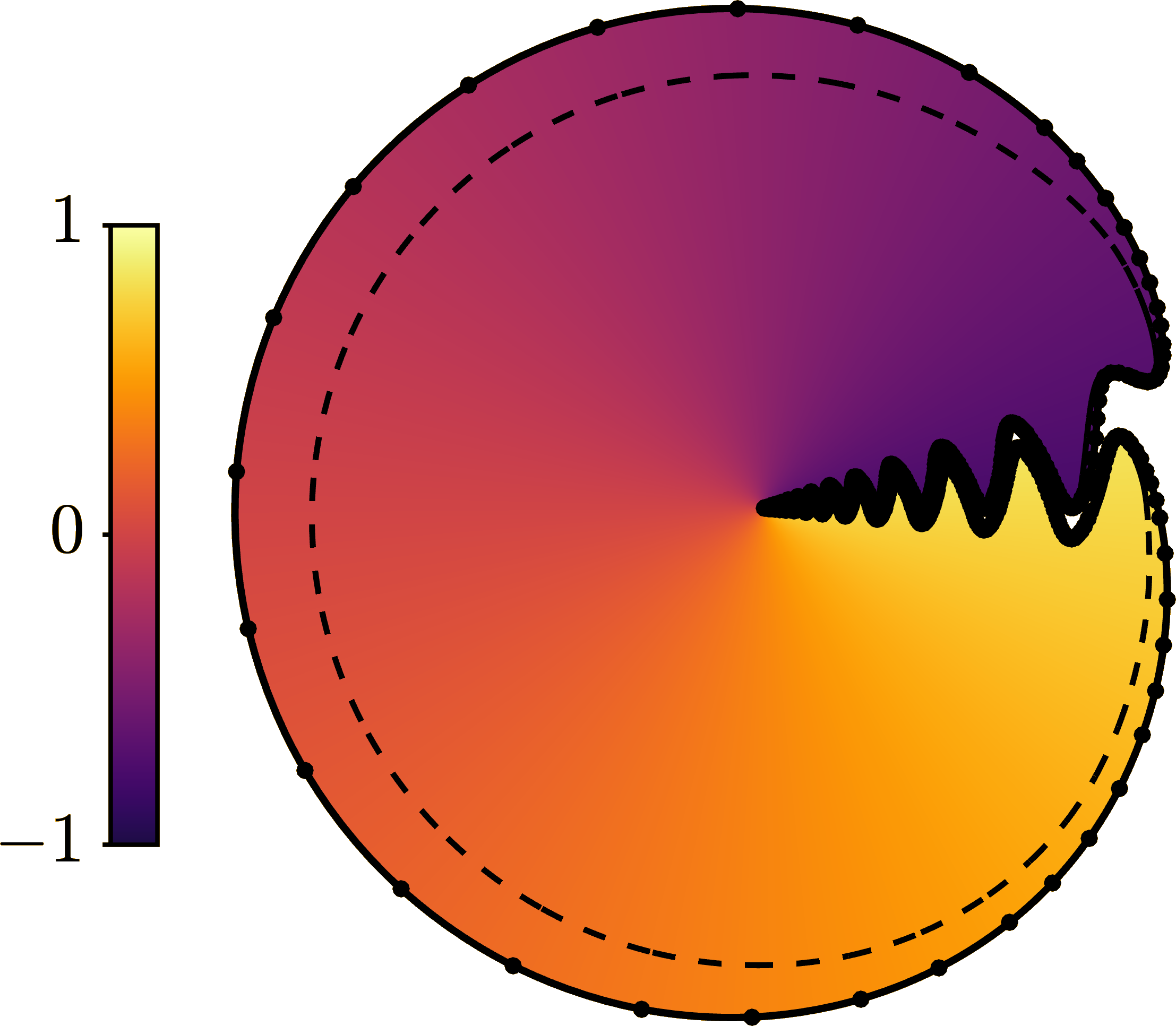}};
        \begin{scope}[x={(image.south east)},y={(image.north west)}]
            \node[anchor=north,inner sep=0,outer sep=-0.01cm] (inset) at (0.65,-0.11) {\includegraphics[height=2cm,trim={0.12cm 0cm 0.12cm 0cm},clip,frame]{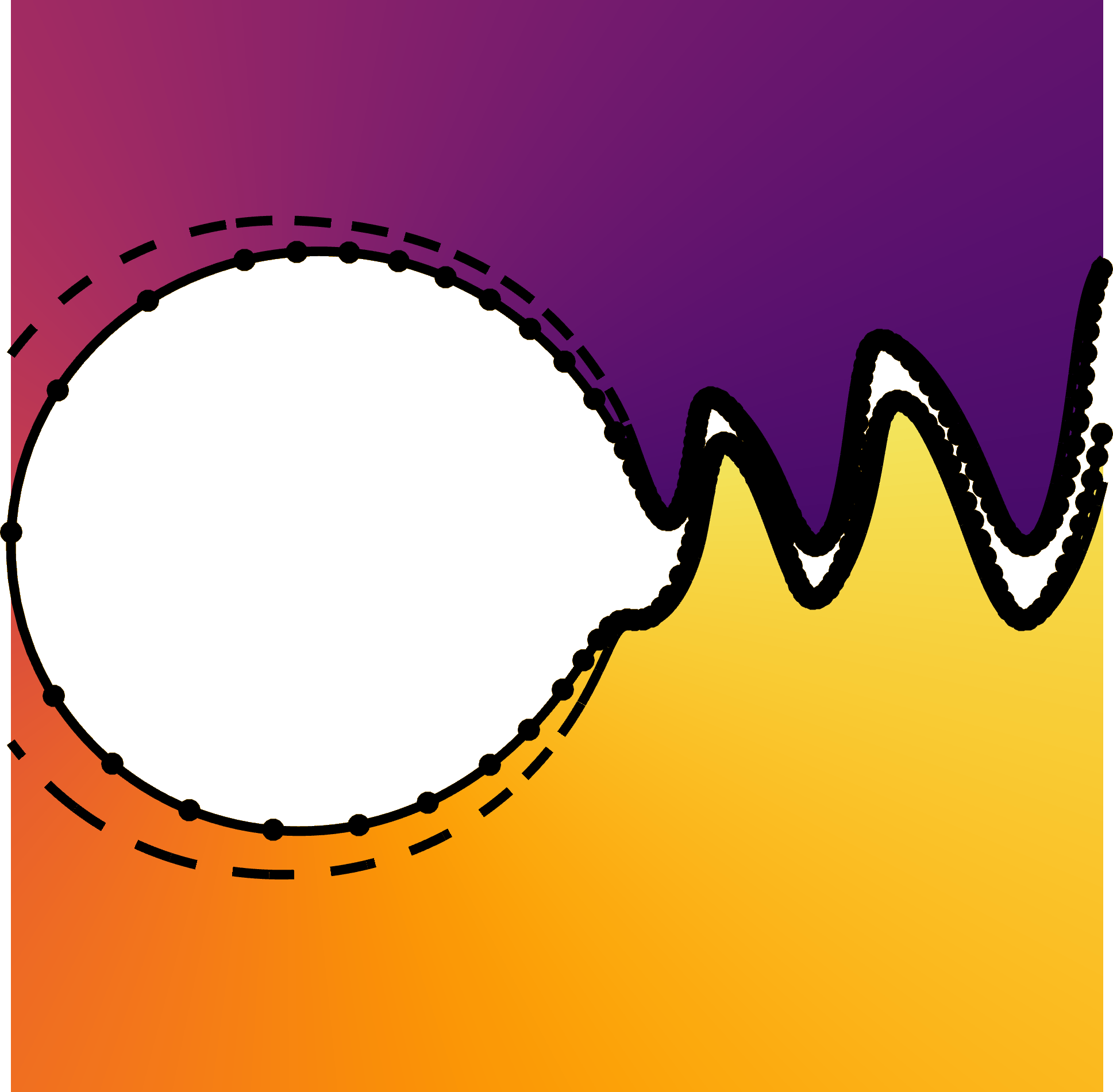}};
            \draw (inset.north west) -- (0.65,0.505);
            \draw (inset.north east) -- (0.65,0.505);
            \node[above=-0.08cm of inset] {\footnotesize\num{10000}\scriptsize$\mspace{-1mu}\times$};
        \end{scope}
    \end{tikzpicture}
    \hspace{0.17cm}
    \begin{tikzpicture}
        \node[anchor=south west,inner sep=0] (image) at (0,0) {\includegraphics[height=3.7cm]{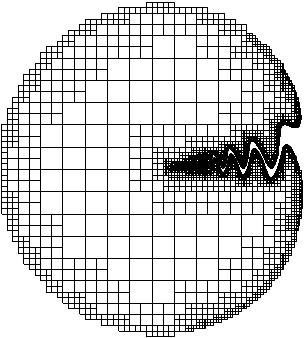}};
        \begin{scope}[x={(image.south east)},y={(image.north west)}]
            \node[anchor=north,inner sep=0,outer sep=-0.01cm] (inset) at (0.57,-0.11) {\includegraphics[height=2cm,trim={0.1cm 0.04cm 0.05cm 0.04cm},clip,clip,frame]{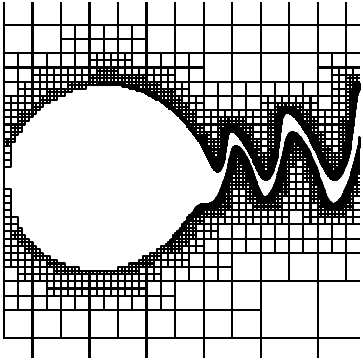}};
            \draw (inset.north west) -- (0.57,0.505);
            \draw (inset.north east) -- (0.57,0.505);
        \end{scope}
    \end{tikzpicture}
    \hspace{0.17cm}
    \begin{tikzpicture}
        \node[anchor=south west,inner sep=0] (image) at (0,0) {\includegraphics[height=3.7cm]{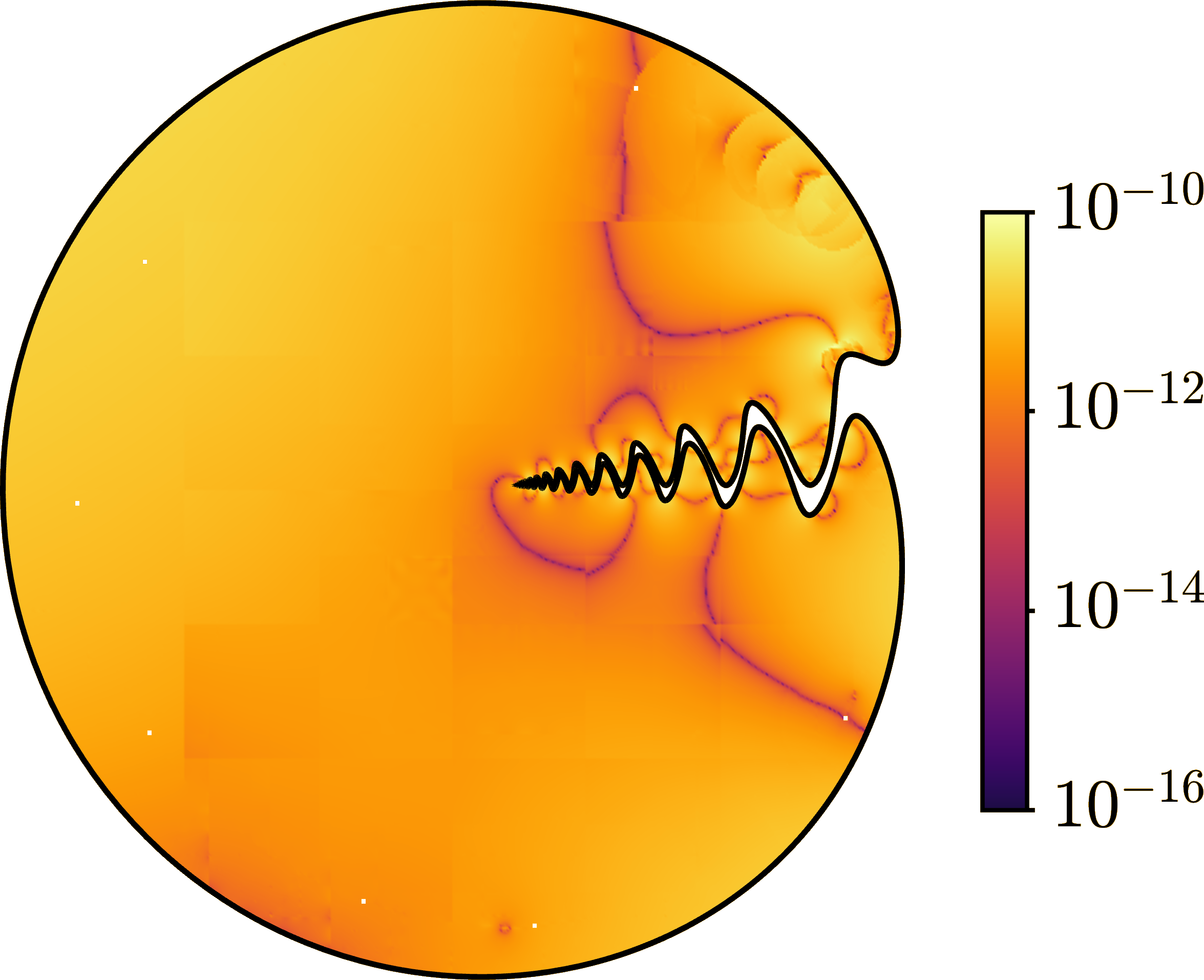}};
        \begin{scope}[x={(image.south east)},y={(image.north west)}]
            \node[anchor=north,inner sep=0,outer sep=-0.01cm] (inset) at (0.43,-0.11) {\includegraphics[height=2cm,trim={0.05cm 0 0.08cm 0},clip,frame]{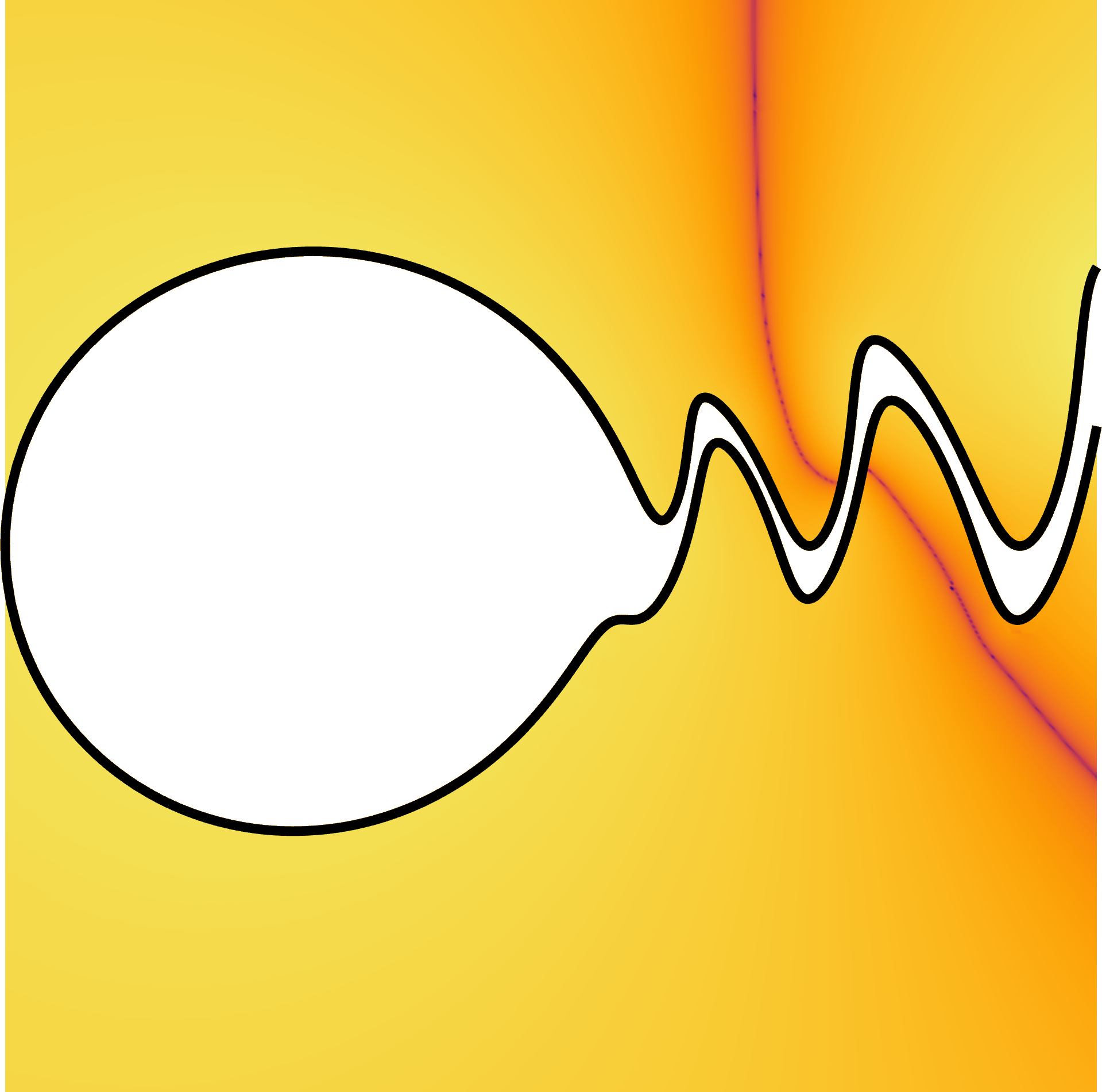}};
            \draw (inset.north west) -- (0.43,0.505);
            \draw (inset.north east) -- (0.43,0.505);
        \end{scope}
    \end{tikzpicture}
	\caption{(Left) A multiscale geometry with length scales spanning four orders of magnitude and close-to-touching ``teeth'' is adaptively panelized into \num{1970} panels. The test solution is chosen so that the inhomogeneity is given by the polar angle, with the branch cut taken between the ``teeth'' of the domain. (Center) The truncated 16th-order quadtree that resolves the inhomogeneity contains over 8 million degrees of freedom. (Right) The computed solution is accurate to about $10^{-10}$.}
	\label{fig:gobbler}
\end{figure}

We now turn to a more extreme multiscale geometry, with length scales continuously spanning four orders of magnitude and close-to-touching regions throughout all scales. This is constructed in $\C$ by spectrally-accurate blending of simple sine wave functions, followed by overall exponentiation \cite{GithubRepo}. Moreover, we prescribe the inhomogeneity $f$ to be the polar angle $\theta = \tan^{-1}(y/x)$ with the branch cut lying between the ``teeth'' of the geometry. This choice of $f$ would pose a problem for methods based on function extension, as the values of $f$ coming from the top and bottom sides conflict. The geometry and inhomogeneity are depicted in \cref{fig:gobbler}, alongside an example 16th-order quadtree and the corresponding pointwise error in the computed solution.

To further illustrate the adaptive performance of our solver, we solve a series of Poisson problems on this multiscale geometry with requested input tolerances $\eps_\text{ask} \in \{10^{-3}, \, 10^{-6}, \, 10^{-9}\}$, and measure the maximum  pointwise error $\eps_\text{get}$ in the resulting solution. \Cref{tab:gobbler_table} shows the results, along with other metrics: the polynomial degree $p \sim \log(1/\eps_\text{ask})$ used to discretize the boundary, quadtree, and strip region; the total number of degrees of freedom $N$ used to represent the solution; the time taken to set up ($T_\text{setup}$), compute a particular solution ($T_\text{part}$), compute a homogeneous correction ($T_\text{homo}$), and evaluate the solution back on the set of quadtree and strip nodes ($T_\text{eval}$); the total time ($T_\text{total}$); and the overall speed of the solver in points per second.

\begin{table}[htb]
    {\footnotesize\caption{Performance results for the adaptive Poisson solver for different requested input tolerances $\eps_\text{ask}$, applied to the multiscale geometry depicted in \cref{fig:gobbler}. We report the maximum pointwise error achieved ($\eps_\text{get}$); the polynomial degree used for the boundary, quadtree, and strip discretizations ($p$); the total number of degrees of freedom used to represent the solution ($N$); the time taken to set up the solver ($T_\text{part}$); the time taken to compute a particular solution ($T_\text{solve}$); the time taken to compute the homogeneous correction ($T_\text{homo}$); the time taken to evaluate the solution back on the set of quadtree and strip nodes ($T_\text{eval})$; the total time ($T_\text{total}$); and the speed in points per second (pps). All times are measured in seconds. Note that we choose $p$ so that $p \sim \log(1/\eps_\text{ask})$.}\label{tab:gobbler_table}
    \bgroup
    \renewcommand{\arraystretch}{1.35}
    \begin{center}
    \begin{tabular}{%
        |c%
        |S[table-format=1.1e-1,round-precision=1,round-mode=places]%
        |c%
        |S[table-format=7.0]%
        |S[table-format=1.2,round-precision=2,round-mode=places]%
        |S[table-format=1.2,round-precision=2,round-mode=places]%
        |S[table-format=1.2,round-precision=2,round-mode=places]%
        |S[table-format=2.2,round-precision=2,round-mode=places]%
        |S[table-format=2.2,round-precision=2,round-mode=places]%
        |S[table-format=6.0]%
    |} 
    \hline
    $\eps_\text{ask}$ & $\eps_\text{get}$ & $p$ & $N$ & $T_\text{setup}$ & $T_\text{part}$ & $T_\text{homo}$ & $T_\text{eval}$ & $T_\text{total}$ & {\!Speed  (pps)\!} \\
    \hline
    $10^{-3}$ & 6.588312340578689e-04 & 3 & 1327108 & 0.99337 & 0.97285 & 1.907  & 8.6756 & 12.549 & 105756 \\
    $10^{-6}$ & 8.189446986501224e-08 & 6 & 4109162 & 1.4604  & 1.4585  & 2.7812 & 21.525 & 27.225 & 150931 \\
    $10^{-9}$ & 1.839950553028658e-10 & 9 & 8418450 & 2.0387  & 2.1939  & 3.6972 & 41.992 & 49.921 & 168633 \\
    \hline
    \end{tabular}
    \end{center}
    \egroup
    }
\end{table}

\subsection{A geometry inspired by cell blebbing}\label{sec:bleb}

\begin{figure}[htb]
	\centering
  	\includegraphics[height=3.4cm]{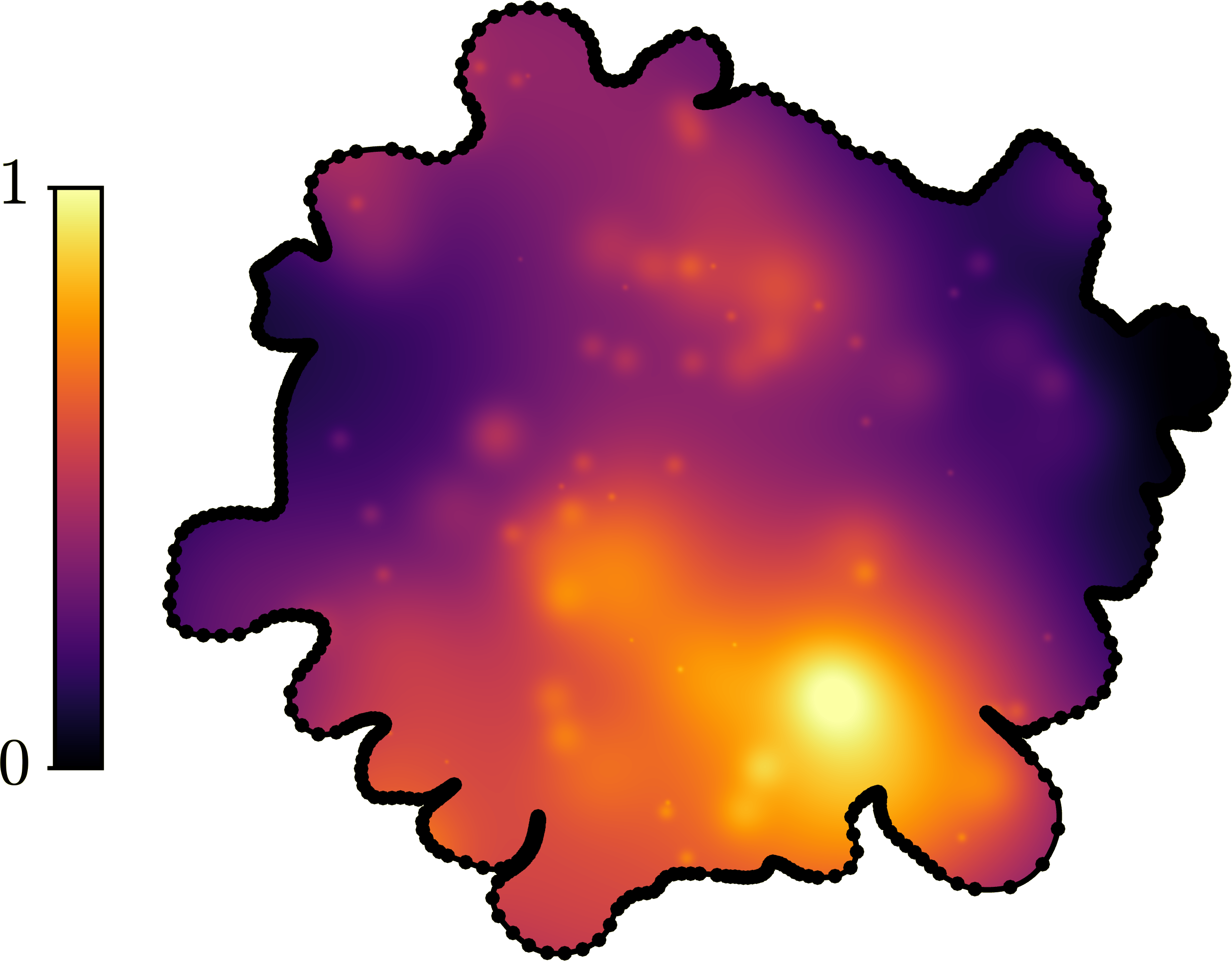}%
	\includegraphics[height=3.4cm]{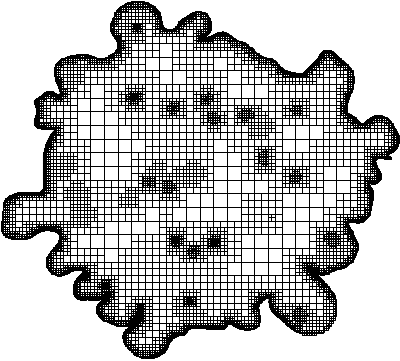}%
	\includegraphics[height=3.4cm]{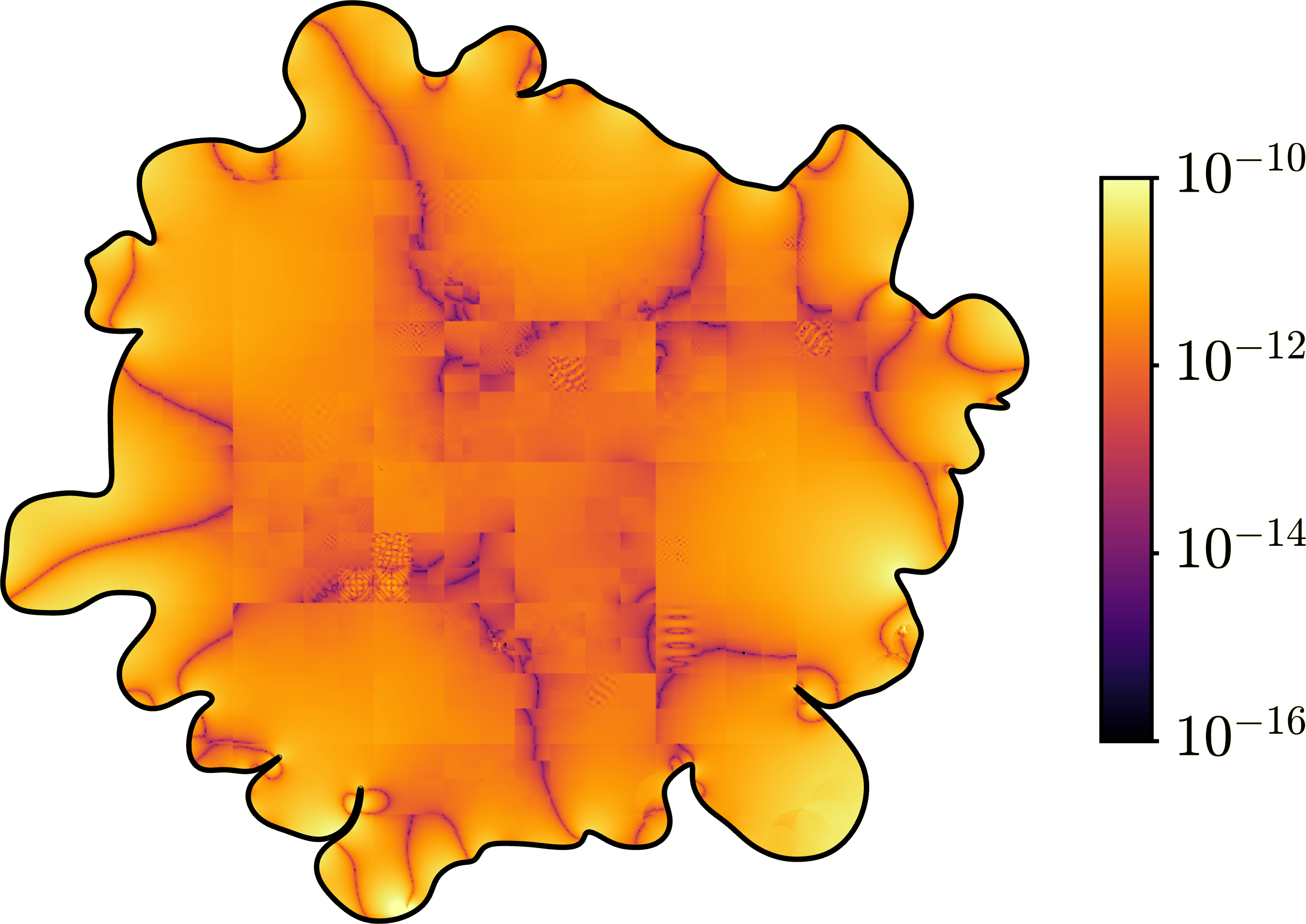}
	\vspace{-1em}
  	\caption{(Left) Computed solution to Poisson's equation on the bleb geometry, with an inhomogeneity consisting of 200 randomly placed Gaussians with variances ranging between $1$ and $10^{-5}$. (Center) The inhomogeneity is adaptively resolved on a 16th-order quadtree with 9.3 million unknowns. (Right) The maximum absolute error in the computed solution is $10^{-10}$.}
	\label{fig:bleb}
\end{figure}

We now apply our solver to a biologically-inspired geometry. Using image processing, we extract the boundary of a cell membrane undergoing blebbing from an image taken from~\cite[Fig. 6]{Charras2008}. We then fit a Fourier series to the sampled boundary and smooth it by convolving with a fixed-width Gaussian. The resulting geometry is shown in~\cref{fig:bleb} and possesses many small folds and thin filaments that would requirement very fine panels everywhere using a method based on uniform grids. We prescribe the solution inside this geometry to be 200 randomly placed Gaussians with variances ranging between $1$ and $10^{-5}$, with the inhomogeneity defined accordingly.

Due to the random placement of the Gaussians, some Gaussians occur very close to the boundary with length scales much smaller than the boundary panelization. Thus, this setup tests the case mentioned in \cref{rem:strip_res}. The solver automatically splits the boundary panels where $f$ was unresolved on the strip grid, and restarts; this process happens a three times, at which point the refined panelization induces a strip which resolved $f$ everywhere. The resulting 16th-order quadtree to resolve $f$ contains 9.3 million degrees of freedom. The maximum absolute error in the computed solution is $10^{-10}$.

\section{Conclusion}
We presented a high-order Poisson solver
that is fully adaptive with respect to both boundary geometry and forcing function.
It combines convolution with the free-space Green's function on an adaptive quadtree that resolves the forcing
function (``box code'') in the bulk, with a curvilinear spectral solver in a boundary ``strip'' region.
Layer potentials on the fictitious strip interface repair the Cauchy matching conditions
to give a particular solution for the whole domain.
The quadtree is truncated within the strip region---we prove that this maintains smoothness in the bulk---%
preventing the need for over-refinement to resolve the boundary.
For this, adapting the strip width function $h(t)$ smoothly to the local boundary panel size is crucial.
We show how our solver efficiently
handles various multiscale geometries, with features spanning up to 8 orders of magnitude,
and compare against FFT-based uniform solvers (which are impractical for 4 or more orders of magnitude).

We expect the adaptive solver presented here to naturally extend to three-dimensional problems, as the strip is based on ``thickening'' an existing boundary discretization. Defining a smooth fictitious ``shell'' in 3D from an unstructured high-order surface triangulation will require partition-of-unity smoothing based on local coordinate charts, as an arc-length parametrization is unavailable in three dimensions. We hope to develop such a solver in future work.

Many more extensions of the adaptive Poisson solver presented here are worthy of exploration. Other inhomogeneous scalar-valued PDEs with known Green's function, such as the Helmholtz or screened Poisson equations, may be solved using essentially the same piecewise representation of the particular solution~\cite{Stein2022b}. Vector-valued PDEs such as the Stokes equations or elastostatics require further development of a vector-valued strip solver. On multiply-connected domains, additional single- and double-layer corrections are needed for each interior boundary. For time-dependent problems with moving geometries or evolving inhomogeneities (e.g., as may arise in a fluid simulation), the truncated quadtree from a previous time step could be updated for the next time step by refining or coarsening only those leaf boxes which overlap the strip, reducing the setup time needed by the solver.

\section*{Acknowledgments}
We have benefited from many useful conversations with Manas Rachh, Charles Epstein, Dhairya Malhotra, Hai Zhu, Leslie Greengard, and Shidong Jiang. We thank Travis Askham for adding support for arbitrary order discretizations to the \texttt{boxcode2d} library~\cite{boxcode2d}.
The Flatiron Institute is a division of the Simons Foundation.

\bibliographystyle{siamplain}
\bibliography{references}

\begin{thebibliography}{10}

\bibitem{afKlinteberg2019}
{\sc L.~{af Klinteberg}, T.~Askham, and M.~C. Kropinski}, {\em A fast integral
  equation method for the two-dimensional {N}avier-{S}tokes equations}, J.
  Comput. Phys., 409 (2020), p.~109353,
  \url{https://doi.org/10.1016/j.jcp.2020.109353}.

\bibitem{Anderson2023}
{\sc T.~Anderson, H.~Zhu, and S.~Veerapaneni}, {\em A fast, high-order scheme
  for evaluating volume potentials on complex {2D} geometries via area-to-line
  integral conversion and domain mappings}, J. Comput. Phys., 472 (2023),
  p.~111688, \url{https://doi.org/10.1016/j.jcp.2022.111688}.

\bibitem{Anderson2022}
{\sc T.~G. Anderson, M.~Bonnet, L.~M. Faria, and C.~Pérez-Arancibia}, {\em
  Fast, high-order numerical evaluation of volume potentials via polynomial
  density interpolation}, J. Comput. Phys., 511 (2024), p.~113091,
  \url{https://doi.org/10.1016/j.jcp.2024.113091}.

\bibitem{Askham2017}
{\sc T.~Askham and A.~J. Cerfon}, {\em An adaptive fast multipole accelerated
  {Poisson} solver for complex geometries}, J. Comput. Phys., 344 (2017),
  pp.~1--22, \url{https://doi.org/10.1016/j.jcp.2017.04.063}.

\bibitem{boxcode2d}
{\sc T.~Askham, F.~Ethridge, D.~Fortunato, Z.~Gimbutas, L.~Greengard,
  M.~O'Neil, and V.~Rokhlin}, 2025,
  \url{https://github.com/flatironinstitute/boxcode2d}.

\bibitem{NonlinearBVPBook}
{\sc K.~B\"{o}hmer}, {\em Numerical Methods for Nonlinear Elliptic Differential
  Equations: A Synopsis}, Oxford University Press, 10 2010,
  \url{https://doi.org/10.1093/acprof:oso/9780199577040.001.0001}.

\bibitem{MultigridBook}
{\sc W.~L. Briggs, V.~E. Henson, and S.~F. McCormick}, {\em A Multigrid
  Tutorial, Second Edition}, SIAM, Philadelphia, PA, 2000,
  \url{https://doi.org/10.1137/1.9780898719505}.

\bibitem{Bruno2022}
{\sc O.~P. Bruno and J.~Paul}, {\em Two-dimensional {{Fourier}} continuation
  and applications}, SIAM J. Sci. Comput., 44 (2022), pp.~A964--A992,
  \url{https://doi.org/10.1137/20M1373189}.

\bibitem{Buzbee1970}
{\sc B.~L. Buzbee, G.~H. Golub, and C.~W. Nielson}, {\em On direct methods for
  solving {P}oisson's equations}, SIAM J. Numer. Anal., 7 (1970), pp.~627--656,
  \url{https://doi.org/10.1137/0707049}.

\bibitem{Charras2008}
{\sc G.~T. Charras}, {\em A short history of blebbing}, J. Microscopy, 231
  (2008), pp.~466--478, \url{https://doi.org/10.1111/j.1365-2818.2008.02059.x}.

\bibitem{Chorin1967}
{\sc A.~J. Chorin}, {\em The numerical solution of the {N}avier-{S}tokes
  equations for an incompressible fluid}, Bull. Amer. Math. Soc., 73 (1967),
  pp.~928--931.

\bibitem{Davis2016}
{\sc T.~A. Davis, S.~Rajamanickam, and W.~M. Sid-Lakhdar}, {\em A survey of
  direct methods for sparse linear systems}, Acta Numer., 25 (2016),
  pp.~383--566, \url{https://doi.org/10.1017/S0962492916000076}.

\bibitem{deBerg2008}
{\sc M.~{de Berg}, O.~Cheong, M.~{van Kreveld}, and M.~Overmars}, {\em
  Computational Geometry: Algorithms and Applications}, Springer Berlin
  Heidelberg, Berlin, Heidelberg, 2008,
  \url{https://doi.org/10.1007/978-3-540-77974-2}.

\bibitem{Demanet2019}
{\sc L.~Demanet and A.~Townsend}, {\em Stable extrapolation of analytic
  functions}, Found. Comput. Math., 19 (2019), pp.~297--331.

\bibitem{Chebfun}
{\sc T.~A. Driscoll, N.~Hale, and L.~N. Trefethen}, {\em Chebfun Guide},
  Pafnuty Publications, 2014, \url{http://www.chebfun.org/docs/guide/}.

\bibitem{Epstein2022a}
{\sc C.~Epstein and S.~Jiang}, {\em A stable, efficient scheme for
  $\mathcal{C}^n$ function extensions on smooth domains in $\mathbb{R}^d$},
  June 2022, \url{https://arxiv.org/abs/2206.11318}.

\bibitem{Ethridge2000}
{\sc F.~Ethridge}, {\em Fast Algorithms for Volume Integrals in Potential
  Theory}, PhD thesis, New York University, 2000.

\bibitem{Ethridge2001}
{\sc F.~Ethridge and L.~Greengard}, {\em A new fast-multipole accelerated
  {Poisson} solver in two dimensions}, SIAM J. Sci. Comput., 23 (2001),
  pp.~741--760, \url{https://doi.org/10.1137/S1064827500369967}.

\bibitem{EvansPDE}
{\sc L.~C. Evans}, {\em Partial Differential Equations}, vol.~19 of Graduate
  Studies in Mathematics, American Mathematical Society, Providence, RI, 1998.

\bibitem{treefun}
{\sc D.~Fortunato}, 2022, \url{https://github.com/danfortunato/treefun}.

\bibitem{GithubRepo}
{\sc D.~Fortunato}, 2022,
  \url{https://github.com/danfortunato/fully-adaptive-poisson}.

\bibitem{Fortunato2020a}
{\sc D.~Fortunato and A.~Townsend}, {\em Fast {Poisson} solvers for spectral
  methods}, IMA J. Numer. Anal., 40 (2020), pp.~1994--2018,
  \url{https://doi.org/10.1093/imanum/drz034}.

\bibitem{Fryklund2022}
{\sc F.~Fryklund and L.~Greengard}, {\em An {FMM} accelerated {P}oisson solver
  for complicated geometries in the plane using function extension}, SIAM J.
  Sci. Comput., 45 (2023), pp.~A3001--A3019,
  \url{https://doi.org/10.1137/22M153495X}.

\bibitem{Fryklund2024}
{\sc F.~Fryklund, L.~Greengard, S.~Jiang, and S.~Potter}, {\em A lightweight,
  geometrically flexible fast algorithm for the evaluation of layer and volume
  potentials}, 2024, \url{https://arxiv.org/abs/2409.11998},
  \url{https://arxiv.org/abs/2409.11998}.

\bibitem{Fryklund2020}
{\sc F.~Fryklund, M.~C.~A. Kropinski, and A.-K. Tornberg}, {\em An integral
  equation–based numerical method for the forced heat equation on complex
  domains}, Adv. Comput. Math., 46 (2020), p.~69,
  \url{https://doi.org/10.1007/s10444-020-09804-z}.

\bibitem{Fryklund2018}
{\sc F.~Fryklund, E.~Lehto, and A.-K. Tornberg}, {\em Partition of unity
  extension of functions on complex domains}, J. Comput. Phys., 375 (2018),
  pp.~57--79, \url{https://doi.org/10.1016/j.jcp.2018.08.012}.

\bibitem{FMMLIB2D}
{\sc Z.~Gimbutas and L.~Greengard},
  \url{https://github.com/zgimbutas/fmmlib2d}.

\bibitem{Greengard1996}
{\sc L.~Greengard and J.-Y. Lee}, {\em A direct adaptive {Poisson} solver of
  arbitrary order accuracy}, J. Comput. Phys., 125 (1996), pp.~415--424,
  \url{https://doi.org/10.1006/jcph.1996.0103}.

\bibitem{Greengard1987}
{\sc L.~Greengard and V.~Rokhlin}, {\em A fast algorithm for particle
  simulations}, J. Comput. Phys., 73 (1987), pp.~325--348,
  \url{https://doi.org/10.1016/0021-9991(87)90140-9}.

\bibitem{Guermond2006}
{\sc J.~Guermond, P.~Minev, and J.~Shen}, {\em An overview of projection
  methods for incompressible flows}, Comput. Meth. Appl. Mech. Engs., 195
  (2006), pp.~6011--6045, \url{https://doi.org/10.1016/j.cma.2005.10.010}.

\bibitem{Hao2014}
{\sc S.~Hao, A.~H. Barnett, P.~G. Martinsson, and P.~Young}, {\em High-order
  accurate methods for {{Nystr\"om}} discretization of integral equations on
  smooth curves in the plane}, Adv. Comput. Math., 40 (2014), pp.~245--272,
  \url{https://doi.org/10.1007/s10444-013-9306-3}.

\bibitem{Helsing2008}
{\sc J.~Helsing and R.~Ojala}, {\em On the evaluation of layer potentials close
  to their sources}, J. Comput. Phys., 227 (2008), pp.~2899--2921,
  \url{https://doi.org/10.1016/j.jcp.2007.11.024}.

\bibitem{hu2024}
{\sc J.~Hu, Y.~Huang, and J.~Lu}, {\em Boundary regularity of {R}iesz
  potential, smooth solution to the chord log-{M}inkowski problem}, 2024,
  \url{https://arxiv.org/abs/2304.14220},
  \url{https://arxiv.org/abs/2304.14220}.

\bibitem{Khoromskij2003}
{\sc B.~N. Khoromskij and J.~M. Melenk}, {\em Boundary concentrated finite
  element methods}, SIAM J. Numer. Anal., 41 (2003), pp.~1--36,
  \url{https://doi.org/10.1137/S0036142901391852}.

\bibitem{Kress2014a}
{\sc R.~Kress}, {\em Linear Integral Equations}, vol.~82 of Applied
  Mathematical Sciences, Springer, New York, NY, 3rd ed.~ed., 2014,
  \url{https://doi.org/10.1007/978-1-4614-9593-2}.

\bibitem{logan}
{\sc D.~Logan}, {\em Applied Mathematics}, Wiley, 4th~ed., 2013.

\bibitem{Malhotra2016}
{\sc D.~Malhotra and G.~Biros}, {\em Algorithm 967: A distributed-memory fast
  multipole method for volume potentials}, ACM Trans. Math. Softw., 43 (2016),
  pp.~17:1--17:27, \url{https://doi.org/10.1145/2898349}.

\bibitem{Mayo1992}
{\sc A.~Mayo}, {\em The rapid evaluation of volume integrals of potential
  theory on general regions}, J. Comput. Phys., 100 (1992), pp.~236--245.

\bibitem{McKenney1995}
{\sc A.~McKenney, L.~Greengard, and A.~Mayo}, {\em A fast {P}oisson solver for
  complex geometries}, J. Comput. Phys., 118 (1995), pp.~348--355,
  \url{https://doi.org/10.1006/jcph.1995.1104}.

\bibitem{mclean}
{\sc W.~C.~H. Mc{L}ean}, {\em Strongly elliptic systems and boundary integral
  equations}, Cambridge University Press, 2000.

\bibitem{Rothe1930}
{\sc E.~Rothe}, {\em Zweidimensionale parabolische {R}andwertaufgaben als
  {G}renzfall eindimensionaler {R}andwertaufgaben}, Math. Ann., 102 (1930),
  pp.~650--670, \url{https://doi.org/10.1007/BF01782368}.

\bibitem{Saye2017}
{\sc R.~Saye}, {\em Implicit mesh discontinuous {{Galerkin}} methods and
  interfacial gauge methods for high-order accurate interface dynamics, with
  applications to surface tension dynamics, rigid body fluid--structure
  interaction, and free surface flow: {{Part I}}}, J. Comput. Phys., 344
  (2017), pp.~647--682, \url{https://doi.org/10.1016/j.jcp.2017.04.076}.

\bibitem{Saye2017a}
{\sc R.~Saye}, {\em Implicit mesh discontinuous {{Galerkin}} methods and
  interfacial gauge methods for high-order accurate interface dynamics, with
  applications to surface tension dynamics, rigid body fluid--structure
  interaction, and free surface flow: {{Part II}}}, J. Comput. Phys., 344
  (2017), pp.~683--723, \url{https://doi.org/10.1016/j.jcp.2017.05.003}.

\bibitem{Shen2024}
{\sc Z.~Shen and K.~Serkh}, {\em Rapid evaluation of {Newtonian} potentials on
  planar domains}, SIAM J. Sci. Comput., 46 (2024), pp.~A609--A628,
  \url{https://doi.org/10.1137/22M1526666}.

\bibitem{slepianI}
{\sc D.~Slepian and H.~O. Pollak}, {\em Prolate spheroidal wave functions,
  {F}ourier analysis and uncertainty, {I}}, Bell Syst.\ Tech. J., 40 (1961),
  pp.~43--64.

\bibitem{Stein2022b}
{\sc D.~B. Stein}, {\em Spectrally accurate solutions to inhomogeneous elliptic
  {{PDE}} in smooth geometries using function intension}, J. Comput. Phys., 470
  (2022), p.~111594, \url{https://doi.org/10.1016/j.jcp.2022.111594}.

\bibitem{Stein2016}
{\sc D.~B. Stein, R.~D. Guy, and B.~Thomases}, {\em Immersed boundary smooth
  extension: {{A}} high-order method for solving {{PDE}} on arbitrary smooth
  domains using {{Fourier}} spectral methods}, J. Comput. Phys., 304 (2016),
  pp.~252--274, \url{https://doi.org/10.1016/j.jcp.2015.10.023}.

\bibitem{Stein2017}
{\sc D.~B. Stein, R.~D. Guy, and B.~Thomases}, {\em Immersed boundary smooth
  extension ({{IBSE}}): {{A}} high-order method for solving incompressible
  flows in arbitrary smooth domains}, J. Comput. Phys., 335 (2017),
  pp.~155--178, \url{https://doi.org/10.1016/j.jcp.2017.01.010}.

\bibitem{Stein2019}
{\sc D.~B. Stein, R.~D. Guy, and B.~Thomases}, {\em Convergent solutions of
  {S}tokes {Oldroyd-B} boundary value problems using the immersed boundary
  smooth extension ({IBSE}) method}, J. Non-Newtonian Fluid Mech., 268 (2019),
  pp.~56--65.

\bibitem{Stein2003}
{\sc E.~M. Stein and R.~Shakarchi}, {\em Complex Analysis}, Princeton Lectures
  in Analysis, No. 2, Princeton University Press, 2003.

\bibitem{kdtree}
{\sc A.~Tagliasacchi}, 2017, \url{https://github.com/taiya/kdtree}.

\bibitem{tao23}
{\sc T.~Tao, G.~Metafune, et~al.}, {\em Regularity of {N}ewtonian potential
  along smooth boundary}, 2023,
  \url{https://mathoverflow.net/questions/446383/regularity-of-newtonian-potential-along-smooth-boundary}.

\bibitem{Townsend2016}
{\sc A.~Townsend, H.~Wilber, and G.~B. Wright}, {\em Computing with functions
  in spherical and polar geometries {I}. {The} sphere}, SIAM J. Sci. Comput.,
  38 (2016), pp.~C403--C425, \url{https://doi.org/10.1137/15M1045855}.

\bibitem{Trefethen2000}
{\sc L.~N. Trefethen}, {\em Spectral Methods in MATLAB}, SIAM, Philadelphia,
  PA, 2000, \url{https://doi.org/10.1137/1.9780898719598}.

\bibitem{ATAP}
{\sc L.~N. Trefethen}, {\em Approximation Theory and Approximation Practice},
  SIAM, Philadelphia, PA, 2013, \url{https://doi.org/10.1137/1.9781611975949}.

\bibitem{Wilber2017}
{\sc H.~Wilber, A.~Townsend, and G.~B. Wright}, {\em Computing with functions
  in spherical and polar geometries {II}. {The} disk}, SIAM J. Sci. Comput., 39
  (2017), pp.~C238--C262, \url{https://doi.org/10.1137/16M1070207}.

\bibitem{Wu2020}
{\sc B.~Wu, H.~Zhu, A.~Barnett, and S.~Veerapaneni}, {\em Solution of
  {{Stokes}} flow in complex nonsmooth {{2D}} geometries via a linear-scaling
  high-order adaptive integral equation scheme}, J. Comput. Phys., 410 (2020),
  p.~109361, \url{https://doi.org/10.1016/j.jcp.2020.109361}.

\end{thebibliography}

\end{document}